\documentclass[12pt]{article}
\usepackage[utf8]{inputenc}
\usepackage{amsmath,amsthm,amssymb,xspace,xcolor,stmaryrd,hyperref,fullpage}
\hypersetup{colorlinks=true, linkcolor=red,citecolor=blue,urlcolor=blue}

\theoremstyle{plain}
\newtheorem{theorem}{Theorem}[section]

\newtheorem{proposition}[theorem]{Proposition}
\newtheorem{lemma}[theorem]{Lemma}
\newtheorem{corollary}[theorem]{Corollary}

\theoremstyle{definition}
\newtheorem{definition}[theorem]{Definition}

\newtheorem{remark}[theorem]{Remark}
\newtheorem{example}[theorem]{Example}

\newcommand{\N}{\mathbb{N}}
\newcommand{\A}{\mathbb{A}}
\newcommand{\PV}{\ensuremath{{\sf PV}}\xspace}
\newcommand{\bit}{\mathit{bit}}
\newcommand{\lh}{\mathit{lh}}
\newcommand{\card}{\mathit{card}}

\newcommand{\Ans}{\mathit{Answer}}

\newcommand{\PVpa}{\mathit{Th}_{\Delta^b_0(\alpha)}(\N)}
\newcommand{\PVa}{\forall\PV(\alpha)}

\newcommand{\LNP}{\mathsf{MIN}}
\newcommand{\LMIN}{\mathsf{LMIN}}

\newcommand{\Ta}{\forall\mathsf{T}^1_2(\alpha)}
\newcommand{\To}{\forall\mathsf{T}^1_2(\PV(\alpha))}
\newcommand{\Sa}{\forall\mathsf{S}^1_2(\alpha)}
\newcommand{\So}{\forall\mathsf{S}^1_2(\PV(\alpha))}
\newcommand{\BA}{\forall\mathsf{T}_2}
\newcommand{\BAa}{\forall\mathsf{T}_2(\alpha)}
\newcommand{\BAo}{\forall\mathsf{T}_2(\PV(\alpha))}

\renewcommand{\le}{\leqslant}

\renewcommand{\ge}{\geqslant}

\newcommand{\ceil}[1]{\lceil #1 \rceil}
\newcommand{\floor}[1]{\lfloor #1 \rfloor}
\newcommand{\ar}{\mathit{ar}}
\newcommand{\str}[1]{\mathcal{#1}}
\newcommand{\q}[1]{\textup{``}#1\textup{''}}

\newcommand{\fle}{\preccurlyeq}
\newcommand{\fge}{\succcurlyeq}

\newcommand{\PaOr}{\mathit{PaOr}}

\newcommand{\WPHP}{\mathit{WPHP}}
\newcommand{\PAR}{\mathit{PAR}}

\newcommand{\PHP}{\mathit{PHP}}
\newcommand{\HOP}{\mathit{HOP}}
\newcommand{\ITER}{\textit{ITER}}
\newcommand{\IND}{\textit{IND}}
\newcommand{\OPHP}{\textit{OPHP}}
\newcommand{\LPHP}{\textit{LPHP}}
\newcommand{\HAP}{\textit{HAP}}
\newcommand{\HDP}{\textit{HDP}}
\newcommand{\rPHP}{\textit{rPHP}}

\begin{document}

\title{\bf Typical forcings, NP search problems\\ and 
an extension of a theorem of Riis}

\author{ Moritz Müller\thanks{supported by  the European Research Council (ERC) under the European Unions Horizon 2020 research programme (grant agreement ERC-2014-CoG 648276 AUTAR).}\\ {\small Universitat Politècnica de Catalunya}
\\[-0.5ex] {\small C. Jordi Girona 1-3, Omega - 327, 08034 Barcelona, Spain}\\[-0.5ex] {\small \tt moritz@cs.upc.edu}}
\date{}

\maketitle
\begin{abstract} 
We define typical forcings encompassing many informal forcing arguments in boun\-ded arithmetic and give general conditions for such forcings to produce models of the universal variant of relativized $\mathsf T^1_2$.
We apply this result to study the relative complexity of total (type 2) NP search problems associated to finitary combinatorial principles. 

Complexity theory compares such problems with respect to polynomial time many-one or Turing reductions. From a logical perspective such problems are graded according to the bounded arithmetic theories that  prove their totality. The logical analogue of a reduction is to prove the totality of one problem from the totality of another. The link between the two perspectives is tight for what we call universal variants of relativized  bounded arithmetics. 
We strengthen   a theorem of Buss and Johnson (2012) that infers relative bounded depth Frege proofs of totality from polynomial time Turing reducibility. 

As an application of our general forcing method we derive a strong form of Riis' finitization theorem (1993).
We extend it by exhibiting a simple model-theoretic property that implies  independence from the universal variant of relativized $\mathsf T^1_2$ plus the weak pigeonhole principle. More generally, we show that the universal variant of relativized $\mathsf T^1_2$ does not prove (the totality of the total NP search problem associated to) a strong finitary combinatorial principle from a weak one. Being weak or strong are simple model-theoretic properties based on the behaviour of the principles with respect to finite structures that are only  partially defined.\\

\noindent{\bf Kewords:} Forcing, Bounded arithmetic, NP search problems, Proof complexity\medskip

\medskip

\noindent{\bf Declarations of interest:} none.

\end{abstract}

\newpage
\tableofcontents

\newpage

\section{Introduction}

While the method of forcing  has ``impressive success in proving independence results for set theory'' \cite[p.81]{bottom}, mathematical logic lacks general methods to prove independence of arithmetical $\Pi_1$-sentences. This lack has been pointed out by 
 Pudl\'ak in \cite{bottom}  
and repeatedly in his latest book \cite{pudlakbuch}. There he asks for a ``method that would be as powerful as forcing and work also for finite problems. [\ldots] To develop such methods is one of the principal goals in proof complexity.''\cite[p.342]{pudlakbuch} A suggestion~\cite{am} is that forcing itself could be developed to become such a method. Indeed, two landmark results of proof complexity, namely the theorems of Riis~\cite{riis,riisbrics} and Ajtai~\cite{ajtai}, have originally been proved by forcing type arguments. In contrast to set theory, however, a general theory of forcing in bounded arithmetic has not been developed.\footnote{An exception is Kraj\'i\c{c}ek's book \cite{kraforce} that follows a conceptually different set-up going back to Scott~\cite{scott}.} 
Instead, later developments ``eliminate the non-standard model theory'' \cite[p.367]{bpu} and forcing. Forcing arguments in bounded arithmetic remain largely informal and confined to the most simple kind of forcings akin to Cohen forcing in set theory.

The leading idea in  Pudl\'ak's book~\cite{pudlakbuch} or the survey~\cite{pudlakbulletin} is that the computational complexity of computational problems associated to sentences could cause independence. 
A particularly appealing instance of this idea is given by true sentences of the form $\forall x\exists y\varphi(x,y)$ where $\varphi(x,y)$ defines some polynomial time decidable and polynomially bounded relation.\footnote{All relevant technical concepts are going to be defined precisely later.} In particular, $\exists y$ is implicitly bounded, so the sentence is $\Pi_1$ in an appropriate language. The associated computational problem is the  (total) NP search problem to compute, given an input $x$, some $y$ such that $\varphi(x,y)$ is true. 
On the computational complexity side~NP search problems are compared using (polynomial time) many-one or Turing reductions and organized into various classes \cite{pls,papa}. An elegant definitorial set-up~\cite{beame} uses {\em type~2}~NP search problems where $\varphi(x,y)$ is allowed to mention  (a predicate for) an oracle $\alpha$.

By a {\em finitary combinatorial principle} we mean an existential first-order sentence $\varphi$ which is valid in the finite. Such a sentence $\varphi$ might or might not have built-in symbols, for example, an order symbol $<$ whose interpretation over universe $[n]:=\{0,\ldots, n-1\}$  is required to be the natural order.  The {\em associated type 2 NP search problem} $Q_\varphi$ asks, given $n$ (in binary) and access to an oracle~$\alpha$ that codes an (exponentially large) structure on $[n]$, to find witnesses to the existential quantifiers in $\varphi$. For example $Q_\WPHP$, for the weak pigeonhole principle $\WPHP$, asks, given $n$ and an oracle $\alpha$ coding a function $f:[n]^2\to [n]$, to find a collision of $f$. 
This problem underlies collision resistant hash functions and is thus important for cryptography (cf.~\cite{krawphp,krapud,krapolywphp,komar}). Further, Papadimitriou's seminal work~\cite{papa} identified a couple of principles $\varphi$ such that 
many natural NP search problems reduce to $Q_\varphi$.

On the logical side, there is a substantial amount of work 
 aimed at characterizing the~NP search problems which are provably total in bounded arithmetics (\cite{approx} contains a recent  survey). For example~\cite{bk}, those provably total in $\mathsf{T}^1_2$ are in the class 
 PLS  from~\cite{pls}, i.e., 
 many-one reducible to $Q_\ITER$, where $\ITER$ is the so-called {\em iteration  principle} 
 with built-in order~$<$. It is not known whether there are NP search problems outside PLS (this would imply\footnote{In fact, $\mathrm{P}\neq\mathrm{TFNP}$ seems to be much stronger than $\mathrm{P}\neq\mathrm{NP}$; see \cite{hubacek} for a recent discussion.}
  $\mathrm{P}\neq\mathrm{NP}$) but there are many such {\em type 2} problems~\cite{morioka}:\footnote{The proof given in \cite{morioka} treats  only many-one reductions. Corollary~\ref{cor:morioka} gives a stronger result.}

\begin{theorem}[Buresh-Oppenheim, Morioka 2004] \label{thm:morioka} If $\varphi$ is a finitary combinatorial principle without built-in symbols that fails in some infinite model, then $Q_\varphi$
is not Turing reducible to~$Q_\ITER$.
\end{theorem}

Papadimitriou's~\cite{papa}  principles exemplify $\varphi$ as above (see Remark~\ref{rem:papa}). 
Beame et al.~\cite{beame} showed that their associated search problems are not equivalent under Turing reductions. Equivalently~\cite{ciy}, the associated complexity classes are distinct relative to a  Cohen-generic oracle. Such oracles are produced by  forcings of the type first considered by Feferman~\cite{feferman}. We refer to~\cite{kurtz} and the references therein for more information about generic oracles. 

These oracle separations use proof  techniques underlying results stating that bounded depth Frege proofs of the propositional translation of one principle from substitution instances of another require exponential size. This  translation is a straightforwardly defined sequence of tautologies, one for each natural $n>0$, expressing totality, i.e., that $\exists y\varphi(n,y)$ is true for all oracles~$\alpha$. 
The similarity of techniques raises the suspicion that the oracle separations might follow from the proof length lower bounds. It took a while for this to be confirmed. Improving~\cite{morioka}, 
Buss and Johnson~\cite{bj}  showed:

\begin{theorem}[Buss, Johnson 2012]\label{thm:bj} Let $\varphi,\psi$ be finitary combinatorial principles. If $Q_\varphi$~is Turing reducible to $Q_\psi$, then there are quasipolynomial size bounded depth Frege proofs of~the propositional translation of $\varphi$ from substitution instances of the propositional translation~of~$\psi$. 
\end{theorem}

In fact, Buss and Johnson got {\em shallow} Frege proofs and were able to prove a partial converse (see~\cite{bj}).  Theorem~\ref{thm:bj} confirms the abovementioned suspicion. Intuitively, however, the proof length lower bounds seem to be much stronger, and it is one of the goals of the present paper to clearly confirm this intuition.

Despite these separations on the (relativized) computational complexity side, it is still open whether full relativized bounded arithmetic $\mathsf T_2(\alpha)$ has more provably total type 2 NP search problems than its second level $\mathsf T^2_2(\alpha)$. This is one of the central open problems in bounded arithmetic (e.g.\  \cite{bkt} or~\cite{approx} survey what is known). 
It is here where a general theory of forcing as Pudl\'ak asks for would be desirable. One of the most beautiful results is\footnote{The statement includes a later improvement due to Kraj\'i\v{c}ek: see \cite[Section 11.5]{krabuch}.} 

\begin{theorem}[Riis 1993]\label{thm:Briis} If $\varphi$ is a finitary combinatorial principle without built-in symbols that fails in some infinite model, then $Q_\varphi$ is not provably total in~$\mathsf T^1_2(\alpha)$.
\end{theorem}

This holds for~$\mathsf S^2_2(\alpha)$ by known conservativity, but fails for~$\mathsf T^2_2(\alpha)$ and~$\WPHP$~\cite{mpw}. 
Riis' original proof~\cite{riis,riisbrics} used a variant of ``the first forcing argument in the context of weak arithmetic'' \cite[p.278]{krabuch} due to Paris and Wilkie~\cite{pw}. These forcings are essentially different from Feferman's forcing mentioned above: the latter
 expands the standard model by an unbounded set while the former expand a nonstandard model by a bounded~set.

\paragraph{Results}
We consider {\em universal variants} of bounded arithmetics and especially the theories $\So,\To,\BAo$ in the language $\PV(\alpha)$ that contains a symbol for every polynomial time algorithm with oracle $\alpha$. They are defined using the same (induction or) minimization schemes as the usual bounded arithmetics $\mathsf S^1_2(\PV(\alpha)), \mathsf T^1_2(\PV(\alpha)),\mathsf T_2(\PV(\alpha))$ but have as base theory $\PVa$, the theory of all universal sentences true in the standard model for all oracles $\alpha$. 
Adding $\PVa$ harmonizes the computational and logical approach to type 2 NP search problems in that the logical notion of consequentiality over various theories coincides with natural notions of reductions. 
 In particular, a type~2 NP search problem $\varphi(x,y)$ is Turing reducible to another $\psi(u,v)$ if and only if $\varphi(x,y)$ is a {\em consequence} of $\psi(u,v)$ over $\So$. This means, roughly, that the totality of $\varphi(x,y)$ is provable in~$\So$ plus the totality of $\psi(u,v)$ for all oracles  that are polynomial time computable relative to $\alpha$.
This follows from known witnessing theorems, the contribution here consists mainly in spelling out the right definitions.
Indeed,
we give  a quite simple proof of

\begin{theorem}\label{thm:bjstrong} 
 Let $\varphi,\psi$ be finitary combinatorial principles. If  $Q_\varphi$ is a consequence of~$Q_\psi$ over $\BAo$, then there are quasipolynomial size bounded depth Frege proofs of the propositional translation of $\varphi$ from substitution instances of the propositional translation~of~$\psi$. 
\end{theorem}

By the equivalence of Turing reducibility and consequentiality over $\So$, this result strengthens  Theorem~\ref{thm:bj} by replacing $\So$ by $\BAo$. Thereby, it confirms the abovementioned intuition that the oracle separations of~\cite{beame} seem to be much weaker than the corresponding proof length lower bounds.

As already mentioned, progress to understand the relative complexity of type 2 NP search problems is hindered by our lack of general methods to prove independence from relativized bounded arithmetics. Here we describe a general forcing method to prove independence from $\To$, situated in the framework of \cite{am}. Many, mostly informal, forcing type arguments in bounded arithmetic  use what we call {\em typical} forcings with {\em typical graded} forcing frames. We prove a general  theorem stating that under a series of simple technical conditions such forcings produce models of $\To$. We stress that this result refers to arbitrary forcings not necessarily of the Cohen type. We refrain from reproducing this rather technical statement here and refer to Theorem~\ref{thm:T12}.  It is meant as a contribution to Pudl\'ak's question in the relativized setting. Our main result, described next, is obtained as an application.

We first reexamine Riis' Theorem~\ref{thm:Briis} in the light of 
Theorem~\ref{thm:T12}.
 We give a new proof using a 
 natural forcing whose conditions are partial oracles that code  {\em partial} structures on~$[n]$ that  embed into an infinite model where $\varphi$ fails and hence do not verify the principle. The generic $\alpha^N$ then codes a {\em total} structure on $[n]$ that falsifies $\varphi$. It is straightforward to verify the conditions of Theorem~\ref{thm:T12} for this forcing, so we get a slight strengthening of Theorem~\ref{thm:Briis} with $\To$ replacing $\mathsf T^1_2(\alpha)$.
This yields the

\begin{corollary}\label{cor:morioka}  If $\varphi$ is a finitary combinatorial principle without built-in symbols  that fails in some infinite model, then 
$Q_\varphi$ is independent from $Q_\ITER$ over $\To$.
\end{corollary}

Being {\em independent} just negates being a consequence. Recalling the relation of Turing reducibility and $\So$, we see that the corollary strengthens  Theorem~\ref{thm:morioka} in that it replaces $\So$ by $\To$. 

Our main interest are finitary combinatorial principles without built-in symbols.
 Theorem~\ref{thm:Briis} suggests to study their relative strength over $\To$.  We aim at a  model-theoretic property implying independence of $Q_\varphi$ from~$Q_{\tilde\varphi}$ over $\To$. Note this is stronger than refuting Turing reducibility.
For example, Buss et al.
~\cite[Theorem~10]{bkt} proved this for~$\tilde\varphi=\WPHP$ and $\varphi=\HOP$.  The {\em Herbrandized ordering principle} $\HOP$ states, roughly, that  partial orders have minimal elements. The proof uses quite involved combinatorics specifically tailored for $\HOP$. Nevertheless, the authors point out that
the proof ``relies on the fact that the injective $\WPHP$ is very over-determined, in the sense that even relatively small subsets of the $a^2$ pigeons must already contain a collision.''~\cite{bkt} 
This hints at the possibility that there is a more general theorem, one concerning independence from ``very over-determined'' principles. We  formalize and quantify the determinacy of a principle and then prove such a general result. 
This is done again by forcing with partial structures.
The cited comment means that the $\WPHP$ is verified in `small' partial structures, i.e., with only a small fraction of function values defined. We shall call such principles {\em weak} in distinction from {\em strong} ones and prove our main result:

\begin{theorem}\label{thm:main} If $\varphi$ is a strong finitary combinatorial principle without built-in symbols and~$\tilde\varphi$ is  a weak finitary combinatorial principle, then
$Q_\varphi$ is independent  from $Q_{\tilde\varphi}$ over $\To$.
\end{theorem}

We view Theorem~\ref{thm:main} as an extension of Riis' Theorem~\ref{thm:Briis} because its proof extends our proof of Theorem~\ref{thm:Briis} which we consider natural and intuitive. Taking $\WPHP$ for $\tilde\varphi$, it gives a simple model-theoretic criterion, namely being strong, for independence from $\To$ ``plus $\WPHP$''. 
We check it applies to many of the commonly studied principles (Section~\ref{sec:disc}), and, in particular, to $\HOP$. 
Compared to~\cite{bkt} our proof is different. First, it does not rely on the already mentioned witnessing theorem for~$\mathsf T^1_2(\alpha)$ by PLS \cite{bk}. Second, it has to sidestep the amplification of failure of $\WPHP$ (cf.~\cite[Section~2]{thapen1}) since this is not available for general weak~$\tilde \varphi$. However, the combinatorial core of the argument is `the same' and isolated as the Core Lemma~\ref{lem:dense2}. Our forcing set-up interprets it as a density argument.

\medskip

We would like to emphasize the comparative simplicity of our proofs of the mentioned results. The proof of Theorem~\ref{thm:bjstrong} proceeds by an intuitive model-theoretic argument followed by an application of the standard propositional simulation. This is technically much simpler than the more direct and quite elaborate construction of  propositional proofs in \cite{bj}. 

The proof of Theorem \ref{thm:main} is a straightforward application of our general forcing Theorem~\ref{thm:T12}. 
Intuitively speaking, the combinatorics needed to fuel the forcing argument are akin to those one would aim at when trying to refute Turing reducibility. The surplus value added by the forcing machinery then consists in strengthening the independence from $\So$ to $\To$.
We hope this can make a point in favor of further developing the general theory of forcing in bounded arithmetic. 

\section{Universal variants of bounded arithmetics}\label{sec:pv}

Usually the bounded arithmetic $\mathsf{S}^1_2$ is written in Buss' language and shown to have a conservative extension that proves  Cook's theory $\PV$~\cite{cookpv}, a theory having symbols for all polynomial time functions. One can add a predicate $\alpha$ and show $\mathsf{S}^1_2(\alpha)$ has a conservative extension $\mathsf{S}^1_2(\PV(\alpha))$ containing $\PV(\alpha)$, Cook's theory for functions computed in polynomial time with oracle $\alpha$. The universal variants $\forall\mathsf{S}^1_2$ and $\So$ use instead~$\forall\PV$ and~$\forall\PV(\alpha)$, respectively, the true universal  theories of polynomial time (with oracle $\alpha$). Basic lemmas concerning bounded arithmetics carry over to the universal variants without surprises, and we sketch the deve\-lop\-ment  only insofar as we shall need it or insofar it allows for a smooth introduction of notations and concepts used later on. 

This section has preliminary character.
Section~\ref{sec:univar} defines universal variants of bounded arithmetics in the languages $\PV$ and~$\PV(\alpha)$. Section~\ref{sec:aux} discusses auxiliary theories in the language $\PV\cup\{\alpha\}$, leading to a useful technical lemma (Lemma~\ref{lem:PVpa}). We prove it via a detour in propositional logic in Section~\ref{sec:prop}, thereby recalling the Paris-Wilkie translation. 
Section~\ref{sec:subst} treats substitutions of formulas for  oracles, and Section~\ref{sec:deforacle} spells out how to define oracle computations and prove conservativity of $\So$ over $\Sa$.

\subsection{Definitions and notations}\label{sec:univar}

A {\em language $L$} is a set  of function and relation symbols $S$ each having an arity $\ar(S)\in \N$. We view constants as nullary function symbols. 
Writing a formula $\varphi$ or a term $t$ as~$\varphi(\bar x)$ or $t(\bar x)$ means that the free variables of $\varphi$ or $t$ are among those in the tuple $\bar x$.
The interpretation of $S\in L$
in an $L$-structure~$\str A$ with universe $A$ is denoted by superscript~$S^A$. The interpretation of a term $t(x_0,\ldots,x_{r-1})$ is denoted $t^A$, a function from $A^r$ into $A$. Often we do not notationally distinguish between $\str A$ and $A$, or omit the superscript when $\str A$ is clear from context.
An $L$-formula {\em with parameters from} $\str A$ is a formula in the language obtained from $L$ by adding every $a\in A$ as a constant. Such formulas are interpreted in $\str A$ understanding 
that the new constants are interpreted by themselves.

\medskip

The language $\PV$ contains the binary relation symbol~$<$ and a function symbol for every  polynomial time Turing machine. We consider every such machine to take as inputs $\bar n\in\N^r$ for some fixed $r\in\N$ which is the arity of its symbol in $\PV$.
The {\em standard $\PV$-model} has universe $\N$ and interprets these symbols by the function computed by the machine, and $<$ by the natural order. We denote the standard model also by $\N$ and do not distinguish notationally between a symbol in $\PV$ and its interpretation in~$\N$. We let
$$
\forall\PV
$$
denote  the set of universal $\PV$-sentences which are true in the standard $\PV$-model $\N$. This theory goes back to 
DeMillo and Lipton \cite{millolipton}.

\medskip

To fix some notation we list some functions in $\PV$. 
It contains the {\em smash} $n\#m:=2^{|n|\cdot|m|}$ where the {\em length} $|n|:=\ceil{\log (n+1)}$ is the length of the binary expansion of $n$ (except that $|0|=0$). 
 We have a binary $\bit(i,n)$
with $n=\sum_{i<|n|}\bit(i,n)\cdot 2^i$ and $\bit(i,n)=0$ for $i\ge |n|$. We think of a number with binary expansion $100110$ as coding the string 
$01100$. We have a binary function in $\PV$ that maps a string to an initial segment of a given length. More precisely, $\PV$ contains a function mapping  
$(n,j)$ with $j<|n|$ to 
$$\textstyle
 n_{<j}:=2^j+\sum_{i<j}\bit(i,n)\cdot 2^i.
$$
 Every $n$ {\em codes} the set $\{i\in\N \mid \bit(i,n)=1\}$ of cardinality $\card(n)$. We also write~$x\in y$ for $\bit(x,y){=}1$.
For every finite sequence $(n_0,\ldots, n_{k-1})\in\N^k$ there is a unique $n\in\N$ such that $\lh(n)=k$ and 
$(n)_i=n_i$ for $i<\lh(n)$ and $(n)_i=0$ for $i\ge \lh(n)$. Here, $\lh(n)$ is a unary function in $\PV$ and  $(n)_i$ is a binary function in $\PV$ applied to $(n,i)$. There is a $k$-ary $t_k\in\PV$ such that $t_k(n_0,\ldots, n_{k-1})=n$; we write $\langle n_0,\ldots, n_{k-1}\rangle$ instead $t_k(n_0,\ldots, n_{k-1})$. Further,
$$
(n)_{<j}
$$
is the code of  $(n_0,\ldots, n_{\min\{k,j\}-1})$.
We assume that for some constant $c>0$
\begin{equation}\label{eq:seqbound}
\textstyle k<|n| < c\cdot  (1+\sum_{i<k}|(n)_i|).
\end{equation}


\medskip 

Let 
$\alpha$ be a unary relation symbol. For a structure $M$ not interpreting $\alpha$ we let $(M,\alpha^M)$ denote its expansion interpreting $\alpha$ by $\alpha^M\subseteq M$. In particular, $(\N,\alpha^\N)$ with
$\alpha^\N\subseteq\N$ is the expansion of the standard $\PV$-model $\N$ which interprets $\alpha$ by $\alpha^\N$. 
This structure has an expansion~$\langle\N, \alpha^\N\rangle$  interpreting the language $\PV(\alpha)$ which 
extends $\PV\cup\{\alpha\}$ by adding a symbol for every polynomial time oracle Turing machine. 
Such a symbol is interpreted in~$\langle\N, \alpha^\N\rangle$ by the function the machine computes when given oracle~$\alpha^\N$. 

\begin{definition}
The theory $\PVa$ is the set of  universal $\PV(\alpha)$-sentences which are true in $\langle\N, \alpha^\N\rangle$  for every $\alpha^\N\subseteq\N$.
\end{definition}

We use standard notations for formula classes. The {\em existential} or {\em universal closure} of a formula is the sentence obtained by existentially or universally quantifying its free variables.
For a set of formulas~$\Phi$ we let $\exists\Phi$ ($\forall\Phi$) be the closure of $\Phi$ under existential (universal) quantification $\exists x$ ($\forall x$).
A formula in a language containing  $\PV$ is {\em bounded} if it is obtained from atomic formulas by Boolean combinations and {\em bounded quantifiers} $\exists x{<}t, \forall x{<}t$ 
where $t$ is a $\PV$-term not containing~$x$. The {\em sharply bounded} formulas are similarly defined but allow only
{\em sharply bounded quantifiers} $\exists x{<}|t|, \forall x{<}|t|$.
We shall always indicate the language in the notation:
 the set of sharply bounded formulas in one of the lan\-guages~$\PV,\PV\cup\{\alpha\}$ or~$\PV(\alpha)$ is denoted by $\Delta_0^b,\Delta_0^b(\alpha)$ and $\Delta_0^b(\PV(\alpha))$ respectively. 
  Closing under positive Boolean combinations, sharply bounded quantification and existential (non-sharply) bounded quantifica\-tion~$\exists x{<}t$ defines the sets
$\Sigma^b_1,\Sigma^b_1(\alpha)$ and~$\Sigma^b_1(\PV(\alpha))$, respectively. 
The sets of bounded formulas are denoted  $\Sigma_\infty^b,\Sigma_\infty^b(\alpha)$ and $\Sigma_\infty^b(\PV(\alpha))$. 

\medskip

The following is easy to see.



\begin{lemma}\label{lem:QE} Every $\Delta_0^b(\PV(\alpha))$-formula is $\PVa$-provably equivalent to 
some quantifier free $\PV(\alpha)$-formula; hence,
 $\PVa$ proves every $\forall\Delta_0^b(\PV(\alpha))$-sentence
 which is true in~$\langle\N,\alpha^\N\rangle$ for every $\alpha^\N\subseteq\N$. 
Analogous statements hold for $\forall\PV$.
\end{lemma}


Let $\Phi$ be a set of formulas.
The {\em minimization scheme}~$\LNP(\Phi)$ and  {\em length minimization scheme}
$\LMIN(\Phi)$ contain, respectively, for every~$\varphi(y,\bar x)\in\Phi$ the universal closure of
\begin{eqnarray*}
&&\varphi(x,\bar x)\to \exists y{\le} x\ \big(\varphi(y,\bar x)\wedge \forall z{<}y \ \neg\varphi(z,\bar x)\big),\\
&&\varphi(x,\bar x)\to \exists y{\le} x \ \big(\varphi(y,\bar x)\wedge \forall z{<}y\big (|z|{<}|y|\to \neg\varphi(z,\bar x)\big)\big).
\end{eqnarray*}

We introduce notation for  {\em universal} variants of some relativized bounded arithmetics, namely those that are going to play a role later on:
\begin{equation}\label{eq:thys}
\begin{array}{lcl}
\BAo&:=&\PVa\cup\LNP(\Sigma^b_\infty(\PV(\alpha))),\\
\To&:=&\PVa\cup \LNP(\Sigma^b_1(\PV(\alpha))),\\
\So&:=&\PVa\cup \LMIN(\Sigma^b_1(\PV(\alpha))).
\end{array}
\end{equation}

The theories $\BA,\forall\mathsf T^1_2,\forall\mathsf S^1_2$ are similarly defined in the language $\PV$ using $\forall\PV$ in place of ~$\PVa$. 
Buss' original theories $\mathsf T_2, \mathsf T^1_2,\mathsf S^1_2$ (in the language $\PV$) are similarly defined but using 
a subset of $\forall\PV$ based on Cook's theory~\cite{cookpv} (cf.~\cite{krabuch}).

\subsection{An auxiliary theory}\label{sec:aux}

Let the theory
$$
\PVpa
$$
consist of all $\forall\Delta_0^b(\alpha)$-sentences which are true in $(\N,\alpha^N)$ for every $\alpha^N\subseteq\N$. Then
define $$\BAa,\Ta,\Sa$$ as in \eqref{eq:thys} but in the language $\PV\cup\{\alpha\}$ and using $\PVpa$ instead $\PVa$.
These definitions look less natural than their analogues in the language $\PV(\alpha)$ but, in fact,  the theories are not really different:
 
\begin{proposition}\label{prop:cons} $\So$ is conservative over $\Sa$, in fact, every model $(M,\alpha^M)$ of~$\Sa$ 
has a unique expansion to a model $\langle M,\alpha^M\rangle$ of $\So$; conversely, every model of $\So$ has this form. 

The same holds for~$\Ta$ and~$\To$, as 
well as for $\BAa$ and $\BAo$
\end{proposition}
 
We give a proof in Section~\ref{sec:deforacle}. We feel $\PVpa$ is the right analogue of $\PVa$ or $\forall\PV$ in the language $\PV\cup\{\alpha\}$ because the analogue of Lemma~\ref{lem:QE} fails. Indeed:

\begin{proposition}\label{prop:ugly}
$ \PVpa$ is not equivalent to a universal theory.
 \end{proposition}
 
 \begin{proof} Let $M$ be a proper elementary extension of the standard $\PV$-model $\N$. 
We claim that  for every $\alpha^M\subseteq M$ all universal sentences of $\PVpa$ are true in $(M,\alpha^M)$.
Equivalently, 
every quantifier free $(\PV\cup\{\alpha\})$-formula~$\varphi(\bar x)$ which is satisfiable in $(M,\alpha^M)$ is also satisfiable in~$(\mathbb N,X)$ for some $X\subseteq\N$.

\medskip

To see this, let $\bar a$ be a tuple from $M$ such that $(M,\alpha^M)\models\varphi(\bar a)$.  The formula $\varphi(\bar x)$ is a Boolean combination of its atomic subformulas of the form 
$$
t(\bar x){<}s(\bar x),\ t(\bar x){=}s(\bar x), \ \alpha(t(\bar x))
$$ 
for certain $\PV$-terms~$t(\bar x),s(\bar x)$. The formula  $\varphi(\bar a)$ is truth functionally satisfied when its atomic subformulas are assigned their truth values in~$(M,\alpha^M)$. Since $M$ is an elementary extension of $\N$,
 there exists a tuple~$\bar n$ from~$\N$ satisfying in~$\N$ the same inequalities (and equalities) of terms appearing in $\varphi(\bar x)$
   as $\bar a$ does in~$M$.  We can thus choose  $X\subseteq\N$ that contains~$t^\N(\bar n)$ if and only if $\alpha^M$ contains~$t^M(\bar a)$; here,~$t(\bar x)$ ranges over the terms appearing in~$\varphi(\bar x)$.  Then~$(\N,X)$ and $\bar n$ give the same truth assignment to the atomic subformulas in~$\varphi(\bar x)$ as $(M,\alpha^M)$ and $\bar a$. Hence, $(\N,X)\models\varphi(\bar n)$, and our claim is proved.

\medskip

It now suffices to show $(M,\alpha^M)\not\models\PVpa$ for $\alpha^M:=\N\subseteq M$. Indeed, any nonstandard $a\in M\setminus\N$ falsifies (plugging for $x$)
$$
\alpha(0)\wedge\forall y{<}|x|(\alpha(y)\to\alpha(y{+}1))\to \alpha(|x|)
$$
in $(M,\alpha^M)$. But the universal closure of this formula is in $\PVpa$. 
 \end{proof}
 
The definitions of $\Sa,\Ta,\BAa$ are robust with respect to these issues:

\begin{lemma}\label{lem:PVpa}
$\forall\PV\cup\LMIN(\Sigma^b_1(\alpha))$ proves $ \PVpa.$
\end{lemma}

We give a proof in the following section via a detour in propositional logic.


\subsection{Propositional logic and simulation}\label{sec:prop}

 Propositional formulas (in negation normal form) are built from literals and constants $0,1$
 using $\vee,\wedge$. A literal is a constant, a variable $X$ or a negated variable $\neg X$.  For a formula $F$ we let $\neg F$ be obtained by swapping $\wedge/\vee$ and  $0/1$ and $X/\neg X$. {\em Depth 0} formulas are literals; {\em depth $d+1$} formulas are depth $d$ formulas or  disjunctions or conjunctions thereof.

We fix a Frege system: a set of finitely many sound inference rules such that any formula~$F$ which is a logical consequence of a set of formulas $\Gamma$, has a Frege proof from $\Gamma$. This is a sequence of formulas ending with $F$ such that all formulas are either from $\Gamma$ or follow from earlier formulas by an inference rule. See \cite[Section~4.4]{krabuch} for precise definitions.
 A depth~$d$  Frege proof is one that contains only depth~$d$ formulas. 
 
Formulas $F$, finite sets of formulas $\Gamma$, Frege proofs $\pi$ etc. are coded by (binary strings coded by) natural numbers. The {\em size} of these objects is the length of the coding number.

The {\em Paris-Wilkie translation} maps  $\Sigma^b_\infty(\alpha)$-sentences $\varphi$ with parameters from $\N$ to 
  propositional formulas $\langle\varphi\rangle$
  in variables written  $\langle \alpha(m)\rangle,m\in\N$. Atomic sentences without~$\alpha$ are mapped to $0$ or $1$ according to their truth value (in the standard $\PV$-model); atoms~$\alpha(t)$ for a closed $\PV$-term $t$ (without variables and with parameters from $\N$) are mapped 
  to~$\langle\alpha(m)\rangle$ for $m:=t^\N$ the value of $t$ in the standard $\PV$-model. Recursively, define $\langle\neg\varphi\rangle:=\neg\langle\varphi\rangle$,
   $\langle\varphi\wedge\psi\rangle:=\langle\varphi\rangle\wedge\langle\psi\rangle$, 
and $
\textstyle\langle\forall y{<}t\ \psi(y)\rangle:=\textstyle\bigwedge_{m<t^{\N}} \langle\psi(m)\rangle
$.

A (partial) assignment $A$ 
{\em agrees with} $\alpha^\N\subseteq\N$ if for all $m\in\N$, $A$ is either undefined on the variable $\langle\alpha(m)\rangle$ or maps it to the truth value of $\alpha(m)$ in $(\N,\alpha^\N)$. If such $A$ is defined on all variables of $\langle\varphi(\bar n)\rangle$, then
\begin{equation}\label{eq:Aalpha}
A\models\langle\varphi(\bar n)\rangle\Longleftrightarrow(\N,\alpha^\N)\models\varphi(\bar n).
\end{equation}

For every fixed $\varphi(x_0,\ldots, x_{k-1})\in \Sigma^b_\infty(\alpha)$ the formulas $\langle\varphi(\bar n)\rangle$ have constant depth and quasipolynomial size. More precisely,
 there is $d\in\N$ such that for all $\bar n=(n_0,\ldots,n_{k-1})\in\N^k$ the formula
$\langle\varphi(\bar n)\rangle$ has depth  $d$ and size  at most $ 2^{(1+\sum_{i<k}|n_i|)^{d}}$.

\begin{proof}[Proof of Lemma~\ref{lem:PVpa}.] 
We formalize \eqref{eq:Aalpha} for  fixed $\varphi(\bar x)\in \Delta_0^b(\alpha)$. There is a $\PV$-function
mapping $\bar n$ to (the code of)~$\langle\varphi(\bar n)\rangle$. We write $\langle\varphi(\bar x)\rangle$ for this function. We code
 assigments~$A$ by sequences with entries $(A)_i$ of the form $\langle m_i,b_i\rangle$ meaning $A(\langle\alpha(m_i)\rangle)=b_i$. 
Choose a  $\Delta_0^{b}(\alpha)$-formula $ \q{z \textit{ agrees with }\alpha}$ 
with the obvious meaning.
Choose a quantifier free $\PV$-formula $\q{z\textit{ is defined on }y}$ defining  (in the standard $\PV$-model) the pairs $(A,F)$ 
of assignments $A$ and formulas $F$ such that $A$ is defined on every variable appearing in $F$. 
Choose a quantifier free $\PV$-formula $\textit{Sat}(x,y)$ defining  the pairs $(A,F)$ of assignments
$A$ and formulas $F$ such that $A$ is defined on $F$ and satisfies~$F$.
Thus
 \eqref{eq:Aalpha} means that $\PVpa$ proves
 \begin{equation}\label{eq:fAalpha}
\q{z \textit{ agrees with }\alpha}\wedge\q{z \textit{ is defined on }\langle\varphi(\bar x)\rangle}\to\big( \textit{Sat}(z,\langle\varphi(\bar x)\rangle) \leftrightarrow \varphi(\bar x) \big).
\end{equation}

\noindent{\em Claim 1.}   $\forall\PV$ proves \eqref{eq:fAalpha}.
\medskip

\noindent{\em Claim 2.}  $\forall\PV\cup\LMIN(\Sigma_1^b(\alpha))$ proves $\exists z \big(\q{z \textit{ agrees with }\alpha}\wedge\q{z \textit{ is defined on }\langle\varphi(\bar x)\rangle}\big).
$ \medskip


We omit the straightforward proofs.
Let $\forall\bar x\varphi(\bar x)\in\PVpa$. By \eqref{eq:Aalpha}, $\langle\varphi(\bar n)\rangle$ is a tautology for every tuple $\bar n$ from $\N$.
In other words, $\forall\PV$ contains the universal closure of
$$
\q{z \textit{ is defined on }\langle\varphi(\bar x)\rangle}\to\textit{Sat}(z,\langle\varphi(\bar x)\rangle).
$$
This and the two claims  imply that 
$\forall\PV\cup\LMIN(\Sigma_1^b(\alpha))$ proves $\varphi(\bar x)$.
\end{proof}


The propositional simulation of $\mathsf{T}_2(\alpha)$ extends to its universal variant $\BAa$: 

\begin{proposition}\label{prop:simulation} Let $\varphi(x_0,\ldots,x_{k-1})\in\Sigma^b_\infty(\alpha)$. If $\BAa$ proves $\varphi(\bar x)$,
then there is $d\in\N$ such that for every $\bar n=(n_0,\ldots,n_{k-1})\in\N^k$ there is a size $2^{(1+\sum_{i<k}|n_i|)^{d}}$ depth $d$ Frege proof of~$\langle\varphi(\bar n)\rangle$. 
\end{proposition}

 \begin{proof} If $\BAa\vdash\varphi(\bar x)$, then $\mathsf{T}_2(\alpha)\vdash(\psi(\bar x,\bar y)\to\varphi(\bar x))$
 for some $\forall\bar x\bar y\psi(\bar x,\bar y)\in\PVpa$. 
By the usual propositional simulation (see~\cite[Corollary 9.1.4]{krabuch}),  for all $\bar n,\bar m$ there is
a constant depth quasipolynomial (in $\bar n,\bar m$) Frege proof of 
$\neg\langle\psi(\bar n,\bar m)\rangle \vee\langle \varphi(\bar n)\rangle$. Choose the all 0 tuple $\bar 0$ for $\bar m$
 and note, as in the previous proof, that $\langle\psi(\bar n,\bar 0)\rangle$ has size polylogarithmic in $\bar n$. 
By \eqref{eq:Aalpha}, $\langle\psi(\bar n,\bar 0)\rangle$ is a tautology, so has a constant depth  proof of size exponential in  $|\langle\psi(\bar n,\bar 0)\rangle|$, so quasipolynomial in~$\bar n$. Modus ponens gives a proof 
 of $\langle \varphi(\bar n)\rangle$.
\end{proof}

 \subsection{Oracle substitutions}\label{sec:subst}
 
We need some notation for substitutions of oracles by formulas: given a $(\PV\cup\{\alpha\})$-formula~$\chi$ and a $\PV(\alpha)$-formula $\psi(u,\bar y)$ let the $\PV(\alpha)$-formula 
\begin{equation}\label{eq:subpsi}
\chi[\alpha/\psi(\cdot,\bar y)]
\end{equation}
be obtained from $\chi$ by replacing each atomic subformula of the form $\alpha(t)$ for some $\PV$-term~$t$ by $\psi(t,\bar y)$. As usual we silently assume that bounded variables in $\chi $ are suitably renamed to become distinct from those in
$\bar y$. 

Of particular interest is the substitution of the oracle by a set polynomial time computable in it. We use special, suggestive notation in this case: for $f(u,\bar y)\in\PV(\alpha)$ we write
\begin{equation}\label{eq:subfct}
\chi[\alpha/f_{\bar y}^{-1}(0)]:=\chi[\alpha/f(\cdot , \bar y){=}0].
\end{equation}

Let $(N,\alpha^N)$ be a $(\PV\cup\{\alpha\})$-structure. A set $A\subseteq N$ is 
{\em $\Delta_1^b(\alpha)$-definable in $(N,\alpha^N)$} if there are $\Sigma^b_1(\alpha)$-formulas with parameters from $N$ defining $A$ and its complement $N\setminus A$.

\begin{lemma} \label{lem:subst} Let $\mathsf T$ be $\Sa,\Ta$ or $\BAa$ and  $(N,\alpha^N) $ be a model of $\mathsf T$. If $A\subseteq N$ is $\Delta_1^b(\alpha)$-definable in $N$, then
$(N,A)\models\mathsf T$.
\end{lemma}

\begin{proof} Consider the case $\mathsf T=\Sa$, the others are similar. 
It is easy to check that $(N,A)$ satisfies $\LMIN(\Sigma_1^b(\alpha))$. That it also satisfies $\PVpa$
then follows from Lemma~\ref{lem:PVpa}.
%
 \end{proof}

 \subsection{Defining oracle computations}\label{sec:deforacle}
 
We think of an oracle computation as a binary decision tree whose inner nodes are labeled by queries to the oracle and whose leaves are labeled by the output. The tree is potentially huge but implicitly feasible in the sense that there is a polynomial time function $t(\bar x,z)$ computing the output or the next query from the input $\bar x$ and the answers $z$ obtained sofar. The answers are coded by the bits of the number $z$, the most significant one not being used.
As a convention, we shall code queries by odd numbers and outputs by even numbers. 


\begin{definition}\label{df:tree}
Let $M$ be a model of $\forall\PV$ and let $t(\bar x,z),h(\bar x)$ be definable functions in~$M$. Then
 $t(\bar x,z)$ is a {\em decision tree (of height at most $h(\bar x)$)  in $M$} if $M$ satisfies
the universal closure of~\eqref{eq:oraclefct} (and~\eqref{eq:qucomp}):
\begin{eqnarray}\label{eq:oraclefct}
i{<}|z|{-}1\wedge t(\bar x,z_{<i})\text{ is even}\ \to\ t(\bar x,z){=}t(\bar x, z_{<i}),\\\label{eq:qucomp}
 h(\bar x){\le}|z|\ \to\ t(\bar x, z)\text{ is even}.
\end{eqnarray}

For $\PV$-terms  $t(\bar x,z),h(\bar x)$ 
we say $t(\bar x,z)$ is a {\em decision tree (of height $h(\bar x)$)} if this holds in all $\PV$-models, that is, if $\forall\PV$  proves~\eqref{eq:oraclefct} (and  \eqref{eq:qucomp}).

Let $\alpha^M\subseteq M$ and $t(\bar x,z)$ be a decision tree in $M$. Then $c\in M$  is a {\em sequence of $\alpha^M$-answers  to $t$ on $\bar a$} if 
$(M,\alpha^M)\models \Ans^\alpha_t(\bar a,c)$ where
\begin{equation}\label{df:answeralpha}
\Ans^\alpha_t(\bar x,z):=\forall i{<}|z|{-}1\ \big( \bit(i,z){=}0\ \leftrightarrow\ \neg\alpha(\lfloor t(\bar x,z_{<i})/2\rfloor)\wedge  t(\bar x,z_{<i})\text{ is odd} \big).
\end{equation}
If additionally, $t(\bar a,c)$ is even (in $M$), we call $c$ {\em  complete}.
\end{definition}

Note that $\Ans^\alpha_t$ is $\Delta^b_0(\alpha)$ if $t$ is a $\PV$-term.


 

\begin{lemma}\label{lem:falphaf} For every  $f(\bar x)\in\PV(\alpha)$ there are $t(\bar x,z),h(\bar x)\in\PV$
such that $t(\bar x,z)$ is a decision tree of height $|h(\bar x)|$
 and $\PVa$ proves 
\begin{equation}\label{eq:falphaf}
f(\bar x) {=}y\ \leftrightarrow\ \exists z{<}h(\bar x) \big(\Ans^\alpha_t(\bar x,z)\wedge t(\bar x,z){=}2y\big)
\end{equation}
\end{lemma} 

\begin{proof} 
Let $f$ correspond to the oracle machine $\A$. Choose $t(\bar x,z)\in\PV$ representing  the following algorithm: on input $(\bar x,z)$   run $\A$ on $\bar x$ answering queries by $\bit(0,z),\bit(1,z),\ldots$ until either $\A$ halts with result $y$ or asks the $|z|$-th query $y$ (hence $\bit(|z|-1,z)=1$ is not used to answer queries); in the first case output $2y$ and in the second output $2y+1$. Choose~$h(\bar x)\in\PV$ such that~$|h(\bar x)|$  is bigger than the number of steps taken by $\A$ on $\bar x$.

Let $g(\bar x)\in\PV(\alpha)$correspond to the oracle machine that on $\bar x$ runs $\A$ and outputs 
the string of oracle answers.
Then $\PVa$ contains the universal closure of
$$
\Ans^\alpha_t(\bar x,g(\bar x))\wedge t(\bar x,g(\bar x)){=}2f(\bar x),
$$
so $f(\bar x) {=}y$ implies the r.h.s.\ of \eqref{eq:falphaf}. The converse is clear by Lemma~\ref{lem:QE}.
\end{proof}

\begin{proof}[Proof of Proposition~\ref{prop:cons}.]  We only prove the first statement. Let $(M,\alpha^M)\models \Sa$.
The theory
$\Sa$ proves that the r.h.s.\ of \eqref{eq:falphaf}
defines (the graph of) a function. We may hence
define the expansion $\langle M,\alpha^M\rangle$ according to the equivalences~\eqref{eq:falphaf}. Uniqueness is clear by Lem\-ma~\ref{lem:falphaf}.  By standard means (see e.g.\ \cite[Theorem~1.3.3.3]{busshand}) the same lemma implies that~$\Sigma^b_1(\PV(\alpha))$-formulas are equivalent to $\Sigma^b_1(\alpha)$-formulas, provably in the theory $\mathsf T$ consisting of $\Sa$ plus the universal closures of the  equivalences~\eqref{eq:falphaf}. Since $\mathsf T$ holds in $\langle M,\alpha^M\rangle$, we can infer
$\mathsf{LMIN}(\Sigma^b_1(\PV(\alpha)))$  from $\mathsf{LMIN}(\Sigma^b_1(\alpha))$. 
Further, $\Sigma^b_1(\alpha)$-formulas are $\Sa$-provably equivalent to 
{\em strict} $\Sigma^b_1(\alpha)$-formulas (see e.g.~\cite[Lemma~5.2.14]{krabuch}), i.e.,  formulas obtained from $\Delta_0^b(\alpha)$-formulas by bounded existential quantification. It follows that quantifier free $\PV(\alpha)$-formulas are $\mathsf T$-provably equivalent
to $\forall\Delta_0^b(\alpha)$-formulas. Thus $\langle M,\alpha^M\rangle\models\PVa$, and we conclude $\langle M,\alpha^M\rangle\models\So$.

Conversely, if $N\models\So$, then its $\PV\cup\{\alpha\}$-reduct models $\Sa$ (Lemma~\ref{lem:QE}). Thus $N$ equals the unique expansion of this reduct to an $\So$-model.
 \end{proof}
 
  
 \section{NP search problems and propositional proofs}
 
 A (type 2) NP search problem is given by a polynomial time decidable (with oracle $\alpha$) and polynomially bounded relation $R(x,y)$ such that for every $x$ there exists $y$ with $R(x,y)$. The computational task is, given $x$ (and oracle $\alpha$), to compute some $y$ with $R(x,y)$. The set of such problems (without oracle) has been introduced to complexity theory in \cite{megiddopapa}. 
 Many natural such problems ask to find a certain configuration in an exponentially large first-order structure given by an oracle $\alpha$. For example, the {\em ($n^2$ to $n$) weak pigeonhole principle $\WPHP$} asks given $n$ (in binary) and an oracle $\alpha$ coding a function $f$ from $[n]^2$ into $[n]$ to find an assignment to $x,y,x',y'$ in $[n]$ satisfying
\begin{equation}\label{eq:wphpformula}
(f(x,y){=}f(x',y')\wedge \neg x{=}x')\ \vee\ (f(x,y){=}f(x',y')\wedge \neg y{=}y').
\end{equation}

Similarly, there is a  type 2 NP search problem $Q_\varphi$ associated to every existential sentence~$\varphi$ which is valid in the finite: given~$n$ (in binary) and an oracle $\alpha$,  witness the existential quantifiers in the structure coded by~$\alpha$ on $[n]$.

Section~\ref{sec:searchprbl} formally defines type 2 NP search problems, Turing and many-one reductions, and characterizes Turing reducibility by $\So$-provability (Proposition~\ref{prop:Tred}). This characterization is essentially known. It is one of our main motivations to study universal variants of bounded arithmetics.
 Section~\ref{sec:coding} discusses two ways 
how to encode finite structures by finite sets, the {\em unary} and the {\em binary} encoding.
Section~\ref{sec:FOsearch}  then formally defines the search problems $Q_\varphi$ above. Finally, we derive Theorem~\ref{thm:bjstrong} in Section~\ref{sec:bj}.
 
\subsection{NP search problems}\label{sec:searchprbl}

Formally, we identify a {\em type 2 NP search problem} with a 
$\Delta_0^b(\alpha)$-for\-mu\-la~$\varphi(x, y)$
such that for some $\PV$-term $t( x)$ the following are true in $(\N,\alpha^{\N})$ for all $\alpha^M\subseteq\N$:
\begin{eqnarray}\label{eq:total}
&&\forall  x\exists  y\varphi( x, y),\\\label{eq:tbd}
&&\forall x\forall y(\varphi(x,y)\to y{<}t(x)) .
\end{eqnarray}
If $\varphi( x, y)$ is a $\Delta_0^b$-formula we speak of a {\em type 1 NP search problem}. We refer to \eqref{eq:total} as the {\em totality} and to \eqref{eq:tbd} as the {\em boundedness} of $\varphi(x,y)$.
We shall discuss examples in Section~\ref{sec:comb}. 
Being {\em solvable in polynomial time} means that there is $f(x)\in\PV(\alpha)$ such that 
$\forall x \varphi(x,f(x))$ is true in $\langle\N,\alpha^\N\rangle$ for all $\alpha^\N\subseteq\N$. By Lemma~\ref{lem:QE} this means that $\PVa$ proves $\varphi(x,f(x))$.

This terminology follows \cite{thakra,thask}
 except that we allow only a unary predicate $\alpha$ instead of an arbitrary finite first-order language. 
 Section~\ref{sec:coding}  spells out how $\alpha$ can code such languages. 
 Our choice allows technically simple definitions of reductions:

 \begin{definition}\label{df:Tred} Let $\varphi(x, y)$ and $\psi( u, v)$ be type 2 NP search problems. We say
 $\varphi( x, y)$ is {\em (polynomial time) Turing reducible} to $\psi( u, v)$ if there are $f(z,x,v),g(x,v),h(x,v)\in\PV(\alpha)$
 and  $\PV$-terms $q(x),s(x)$ such that the universal closure of
 \begin{equation}\label{eq:Tred}
 \begin{split}
&g(x,w){<}s(x)\ \wedge  
\Big(\forall i{<}|q(x)| \ \psi\big(g(x,(v)_{<i}),(v)_{i}\big)\big[\alpha/f(\cdot,x,(v)_{<i}){=}0\big]
\to\ \varphi(x,h(x,v))\Big).
\end{split}
 \end{equation}
 is true in $\langle\N,\alpha^N\rangle$ for all $\alpha^N\subseteq\N$.
If $q(x)$ is the constant $\ell\in\N$, we speak of Turing reducibility {\em  with $\ell$ queries}; and for $\ell=1$ we speak of {\em many-one} reducibility.
 \end{definition}

  Intuitively,  \eqref{eq:Tred} interprets $v$ as a sequence $(v)_0,(v)_1,\ldots$ of answers to oracle queries to~$\psi( u, v)$. The $i$-th query 
  is given by some instance $u_i$ and oracle $\beta_i$ computed by $g$ and~$f$ from the input $x$ and answers $(v)_0,\ldots, (v)_{i-1}$ obtained sofar, namely $u_i:=g(x,(v)_{<i})$ and $\beta_i:=f(\cdot,x,(v)_{<i})^{-1}(0)$. Finally,~$h$ returns some solution to $\varphi(x,y)$. 
  For simplicity,
the formalization \eqref{eq:Tred} assumes Turing reductions to always make the same number  of queries, namely  $|q(x)|$, independently from the answers obtained.  The first conjunct ensures  that
the size $|u_i|$ of the queries is bounded by $|s(x)|$, hence the whole computation runs in time polynomial in $|x|$. This first conjunct can be omitted 
 in case $q(x)$ is constant. In particular,  $\varphi( x, y)$ is many-one reducible to $\psi( u, v)$ if and only if there are $f,g,h\in \PV(\alpha)$ such that for all $\alpha^\N\subseteq \N$, $\langle\N,\alpha^\N\rangle$ satisfies the universal closure of
\begin{equation}\label{eq:manyone}
  \psi(g(x),v) \big[\alpha/f_x^{-1}(0)\big]
 \to\ \varphi(x,h(x,v)).
\end{equation}

The following notation mimics notations  like $\WPHP[\PV(\alpha)]$ familiar from the literature.
For a type 2 NP search problem $\varphi(x,y)$, let
\begin{equation}\label{eq:subPV}
\varphi[\PV(\alpha)]:=\Big\{  \forall x\forall \bar z\exists y\ \varphi(x,y)[\alpha/f^{-1}_{\bar z}(0)]\mid f(u,\bar z)\in\PV(\alpha) \Big\}.
\end{equation} 

 \begin{definition}\label{df:conseq} Let $\varphi(x, y)$ and $\psi( u, v)$ be type 2 NP search problems and $\mathsf T$ a theory. 
 We say $\varphi( x, y)$ is a {\em consequence of  $\psi( u, v)$ over $\mathsf T$} if
 $ \mathsf T\cup \psi[\PV(\alpha)]$ proves $ \exists y\varphi(x,y).$ 
 
 Otherwise
we say  $\varphi( x, y)$ is {\em independent over $\mathsf T$ from $\psi( u, v)$}.
 \end{definition}

We are not aware of a reference for this notion of consequence. It is a natural logical analogue of the complexity theoretic notion of reducibility. The mode of speech follows Hanika \cite[Definition~4.4]{hanikathesis} whose notion is  weaker in that  $\mathsf T$ is given only one sentence from~$\psi[\PV(\alpha)]$ when asked to prove  $ \exists y\varphi(x,y)$. 
We state the following only for the universal variants of bounded arithmetics that we explicitly defined but it is clear from the proof that it holds for other universal variants as well.

\begin{proposition}\label{prop:trans} Let $\mathsf T$ be  $\So,\To$ or $\BAo$. Consequentiality over~$\mathsf T$ is transitive as a relation over type 2 NP search problems.
\end{proposition}

\begin{proof} Let $\mathsf T'$ be $\Sa,\Ta$ or $\BAa$ if  $\mathsf T$ is $\So,\To$ or $\BAo$, respectively. Suppose $\varphi( x, y)$, $\psi(u,v)$, $\chi( z, w)$ are type 2 NP search problems and $\psi(u,v)$ is a consequence of 
$\chi( z, w)$ over $\mathsf T$ and $\varphi( x, y)$ is a consequence of $\psi( u, v)$ over $\mathsf T$. We have to show that
$\varphi( x, y)$ is a consequence of $\chi( z, w)$ over $\mathsf T$.

Let a model of $\mathsf T\cup\chi[\PV(\alpha)]$ be given. By Proposition~\ref{prop:cons} it has the form $\langle N,\alpha^N\rangle$ for~$(N,\alpha^N)\models\mathsf T'$. For contradiction, assume 
$\langle N,\alpha^N\rangle\models\forall y\neg\varphi(n,y)$ for some $n\in N$. Then $\langle N,\alpha^N\rangle\models\forall v\neg\psi(m,v)[\alpha/f_{\bar a}^{-1}(0)]$ for certain  $f(z,\bar z)\in\PV(\alpha)$ and $\bar a,m$ from~$N$. Let~$A$ be the set defined by $f(z,\bar a){=}0$ in $\langle N,\alpha^N\rangle$ and note it is $\Delta^b_1(\alpha)$-definable in  $(N,\alpha^N)$ by Lemma~\ref{lem:falphaf}.
Then $(N,A)\models \forall v\neg\psi(m,v)$ and $(N,A)\models\mathsf T'$ by Lemma~\ref{lem:subst}, so  
$\langle N, A\rangle\models\mathsf T$ by Proposition~\ref{prop:cons}. We are left to show $\langle N,A\rangle\models\chi[\PV(\alpha)]$.

For contradiction, assume  there are 
$g(u,\bar u)\in\PV(\alpha)$ and $\bar b,k$ from $N$ such that
$\langle N,A\rangle\models\forall w\neg\chi(k,w)[\alpha/g_{\bar b}^{-1}(0)]$. If $\theta(u,\bar u)$ is a $(\PV\cup\{\alpha\})$-formula equivalent to $g(u,\bar u){=}0$ in $\langle N,A\rangle$, then $(N,A)$ satisfies $\forall w\neg\chi(k,w)[\alpha/\theta(\cdot,\bar b)]$ and thus 
\begin{equation}\label{eq:theta}
\langle N,\alpha^N\rangle\models\forall w\neg\chi(k,w)\big[\alpha/\theta(\cdot,\bar b)[\alpha/f_{\bar a}^{-1}(0)]\big].
\end{equation} 

To get the desired contradiction it suffices to choose $\theta$ such that
$\theta(u,\bar u)[\alpha/f_{\bar z}^{-1}(0)]$ is $\PVa$-provably equivalent to $\ell(u,\bar u,\bar z){=}0$ for some $\ell\in\PV(\alpha)$. Indeed, 
then \eqref{eq:theta} gives  $\langle N,\alpha^N\rangle\models\forall w\neg\chi(k,w)\big[\alpha/\ell_{\bar b,\bar a}^{-1}(0)]$, contradicting $\langle N,\alpha^N\rangle\models\chi[\PV(\alpha)]$.

We choose $\PV$-terms $t(u,\bar u,z),h(u,\bar u)\in\PV$ according Lemma~\ref{lem:falphaf} for $g(u,\bar u)$ and take  the r.h.s. of \eqref{eq:falphaf} for $\theta(u,\bar u)$. Then $\theta(u,\bar u)[\alpha/f_{\bar z}^{-1}(0)]$ is $\Sigma^b_1(\PV(\alpha))$. But the leading  $\exists z{<}h(u,\bar u)$ can be eliminated by replacing $z$ by $s(u,\bar u,\bar z)$ for a suitable $s\in \PV(\alpha)$. The resulting $\Delta^b_0(\PV(\alpha))$-formula is $\PVa$-provably equivalent to  $\ell(u,\bar u,\bar z){=}0$ for a suitable  $\ell\in\PV(\alpha)$.
\end{proof}

The following proposition characterizes natural reducibilities by consequentiality over the universal variants $\PVa$ and $\So$. 
The interesting directions from right to left follow from known witnessing theorems. In \cite[Fact~4.6]{hanikathesis} and \cite[Proposition~7.1]{pudlakbulletin} proofs appear for  $\mathsf{S}^1_2(\alpha)$ and only one member of $\psi[\PV(\alpha)]$. 
 The converse directions from left to right are easy given the definition of~$\PVa$.

\begin{proposition} \label{prop:Tred}  Let $\varphi(x, y)$ and $\psi( u, v)$ be type 2 NP search problems. 
\begin{enumerate}\itemsep=0pt
\item[(a)] $\varphi( x, y)$ is  Turing reducible to $\psi( u, v)$ if and only if $\varphi(x,y)$ is a consequence of $\psi( u, v)$ over $\So$.
\item[(b)] There is $\ell\in\N$ such that $\varphi( x, y)$ is  Turing reducible to $\psi( u, v)$ with $\ell$ queries if and only if 
$\varphi(x,y)$ is a consequence of $\psi( u, v)$ over $\PVa$.
\item[(c)]  $\varphi( x, y)$ is solvable in polynomial time if and only if $\PVa$ proves $\exists y\varphi(x,y)$, if and only if
$\So$ proves  $\exists y\varphi(x,y)$.
\end{enumerate}
\end{proposition}

\begin{proof} We first prove (a). For the  direction from left to right,
assume $\varphi( x, y)$ is  Turing reducible to $\psi( u, v)$. By Lemma~\ref{lem:QE}, $\PVa$ proves
the $\Delta_0^b(\PV(\alpha))$-formula   \eqref{eq:Tred}.  We argue in $\So\cup\psi[\PV(\alpha)]$ that there exists $v$ satisfying
$$
\forall i{<}|q(x)| \ \psi\big(g(x,(v)_{<i}),(v)_{i}\big)\big[\alpha/f(\cdot,x,(v)_{<i}){=}0\big].
$$
Let $\chi_0(x,v)$ be obtained by replacing $\forall i{<}|q(x)|$ by $\forall i{<}\lh(v)$, and 
consider the formula
$$
\chi_1(x,w):=\exists v\big( \lh(v){=}|q(x)|{-}|w|\wedge  \chi_0(x,v) \big).
$$

Let $t(x)$ witness boundedness \eqref{eq:tbd} of $\psi(u,v)$. We can assume~$t(x)$ is non-decreasing, i.e., $\forall\PV$ proves $(x{\le} x'\to t(x){\le} t(x'))$. By \eqref{eq:seqbound} the quanti\-fier~$\forall i{<}\lh(v)$ can be sharply bounded in$v$ and  $\exists v$ can be bounded by~$c\#(q(x)\# t(s(x)))$ for a  suitable  $c\in\N$. Hence, $\chi_1(x,w)$ is 
$\forall\PV$-provably equivalent to a $\Sigma_1^b(\PV(\alpha))$-formula.
Since trivially $\chi_1(x,q(x))$,  $\LMIN(\Sigma_1^b(\PV(\alpha)))$ gives  a minimal length $w$ with $\chi_1(x,w)$. Then $|w|=0$ because each answer sequence  $v$ can be prolongued by any  $v'$
with $\psi\big(g(x,v),v'\big)\big[\alpha/f(\cdot,x,v){=}0\big]$; and such~$v'$ exists by $\psi[\PV(\alpha)]$.

\medskip

For the  direction from right to left, assume $\So\cup\psi[\PV(\alpha)]$ proves $ \exists y\varphi(x,y)$. 
By compactness there are $\ell\in\N$ and $f_0(z,\bar z_0),\ldots, f_{\ell-1}(z,\bar z_{k-1})\in\PV(\alpha)$ and a quantifier free $\PV(\alpha)$-formula $\chi(\bar x)$ such that $\forall\bar x\chi(\bar x)\in\PVa$ and such that $\mathsf{S}^1_2(\PV(\alpha))$ proves
\begin{eqnarray*}
&&\exists y\bar w\forall\bar v\ \theta(x,y,\bar w,\bar v),\textup{ where}\\
&&\quad\bar w:= \bar x \ u_0\bar z_0\cdots u_{\ell-1}\bar z_{\ell-1}, \\
&&\quad\bar v:=v_0\cdots v_{\ell-1},\\
&&\quad\theta:=\textstyle
\neg \chi(\bar x)\vee\bigvee_{i<\ell} \neg\psi(u_i,v_i)\big[\alpha/f_i(\cdot,\bar z_{i}){=}0  \big]\vee\varphi(x,y).
\end{eqnarray*}

By a well-known witnessing argument (see \cite[Theorem~7.3.3]{krabuch}) a witness tuple $y\bar w$ is computable from $x$ by a polynomial time counterexample computation~\cite{counter}: a polynomial time Student computes  a candidate $y^0\bar w^0$ and sends it to a computationally unbounded Teacher; Teacher answers with
a counterexample $\bar v^0$, i.e., such that $\neg\theta(x,y^0,\bar w^0,\bar v^0)$; then Student computes another candidate~$y^1\bar w^1$ and Teacher answers with a counterexample $\bar v^1$ and so on, until Student finally computes~$y^t\bar w^t$ such that no counterexample exists; then the computation stops with output $y^t\bar w^t$. We can assume that $t$ equals $|q(x)|$ for some $\PV$-term~$q(x)$, independent of Teacher's answers. The whole computation runs in time polynomial in $|x|$, so there is a $\PV$-term $s(x)$ bounding all components of all $y^j\bar w^j$'s, and in particular the $u^j_i$'s.

Now just note that each answer from Teacher can be simulated by $\ell$ 
oracle calls to~$\psi$, namely to get $\bar v^j$ such that $\psi(u^j_i,v^j_i)\big[\alpha/f_i(\cdot,\bar z^j_{i}){=}0  \big]$ for all $i<\ell$ where we write $\bar w^j$ as $\bar x^j \ u^j_0\bar z^j_1\cdots u^j_{\ell-1}\bar z^j_{\ell-1}$ and $\bar v^j$ as $ v^j_0\cdots v^j_{\ell-1}$.
More specifically, we look for $f,g,h\in\PV(\alpha)$ such that the universal closure of \eqref{eq:Tred} is true. The function 
$f$, given $x$ and previous answers, simulates the functions $ f_i(\cdot,\bar z^j_{i})$ for the $\bar z^j_i$ computed by Student; 
$g$ computes the $u^j_i$'s by simulating Student;  $h$ outputs $y^t$ of Student's final candidate $y^t\bar w^t$.

\medskip

The proof of  (b) is similar but simpler. For the forward direction, $\ell$ many $v_0,\ldots,v_{\ell-1}$ can be collected in the tuple  $\langle v_0,\ldots, v_{\ell-1}\rangle$ without need to rely on $\LMIN(\Sigma_1^b(\PV(\alpha)))$. For the converse, $\PV(\alpha)$-provability yields a counterexample computation with constantly many rounds. This follows from the KPT-Theorem~\cite{kpt}, in fact, a simple version of it proved in \cite[Theorem~2.2]{cooktha} by a simple proof that works for $\PVa$.

\medskip

The two forward directions of (c) are clear (recall Lemma~\ref{lem:QE}). The last statement implies the first  
by applying (a) with  $v{=}v$ for $\psi(u,v)$.
\end{proof}

The type 2 NP search problems provably total in universal variants of bounded arithmetics  form a meaningful complexity class in that they are closed under Turing reductions. Again, we state this only for
$\So,\To$ and $\BAo$. Note that by the previous proposition  consequentiality over these theories is implied by Turing reducibility.

\begin{corollary}\label{cor:Tredclosed} Let $\mathsf T$ be  $\So,\To$ or $\BAo$, and let  $\varphi(x,y)$ and  $\psi(u,v)$ be type 2 NP search problems.
If  $\varphi(x,y)$ is a consequence of  $\psi(u,v)$ over $\mathsf T$ and $\mathsf T$ proves $\exists v\psi(u,v)$, then $\mathsf T$ proves $\exists y\varphi(x,y)$.
\end{corollary}

\begin{proof} Assume that $\varphi(x,y)$ is a consequence of  $\psi(u,v)$ over $\mathsf T$  and $\mathsf T$ proves $\exists v\psi(u,v)$. The latter is equivalent to  $\psi(u,v)$ being a consequence of $w{=}w$ over $\mathsf T$. By Proposition~\ref{prop:trans}, $\varphi(x,y)$ is a consequence of $w{=}w$ over $\mathsf T$. Hence $\mathsf T$ proves $\exists y\varphi(x,y)$.
%
\end{proof}

It might be worthwhile to look for complexity theoretic reductions equivalent to consequentiality over  higher levels of the bounded arithmetic hierarchy (cf.~\cite[Section~7]{pudlakbulletin}). Such a notion of reduction is implicit in \cite[Proof of Theorem~8]{bkt}
for the special case of $\mathsf{T}^1_2(\alpha)$ and~$\psi$ the search problem associated to the weak pigeonhole principle~\eqref{eq:wphpformula}. To define this and similar problems we need to agree on how to code finite structures by oracles.

\subsection{Unary and binary codes of structures}\label{sec:coding}

There are at least two common ways how to code  structures by oracles, namely, the unary and the binary encoding.
The unary encoding codes functions by their graphs while the binary encoding uses their bit graphs. Both codings work not only over $\N$ but over certain non-standard models too.

Let $(N,\alpha^N)$ be a model of $\Sa$, and let $L$ be a finite language. For notational simplicity we assume $L\subseteq \N$ and $\N$ is an initial segment of $N$. Recall that $\ar(S)$ denotes the arity of the symbol $S\in L$.
For $n\in N$ we write
$$
[n]:=\{a\in N\mid a<n\}.
$$
Here and  below we omit superscripts as in  $<^N$ for  interpretations of $\PV$-symbols in $N$.

\begin{definition}\label{df:Aalpha} Let $n\in N\setminus\{0\}$. We say $\alpha^N$ is the {\em unary code (in $N$)} of the
 $L$-structure 
 $\str A(L,n,\alpha^N)
 $
with universe $[n]$ if $\alpha^N$ contains exactly the tuples
$\langle S,\bar a\rangle\in\alpha^N$ where $S\in L$ is a relation symbol and  $\bar a\in S^{[n]}$ (the interpretation of $S $ 
in $\str A(L,n,\alpha^N)$), or 
$\langle S,\bar a,b\rangle\in\alpha^N$ where $S\in L$ is a function symbol and  $ S^{[n]}(\bar a)=b$. 
 If such a structure exists, we say $\str A(L,n,\alpha^N)$ {\em is defined (in~$(N,\alpha^N)$)}; otherwise, the notation
$\str A(L,n,\alpha^N)$ is undefined.
\end{definition}

A disadvantage of the unary code is that not every set
 $\alpha^\N\subseteq \N$ is the unary code of an $L$-structure on~$[n]$ because the relations determined for function symbols have to be graphs of functions on~$[n]$.
 Another disadvantage is that function symbols cannot be evaluated in polynomial time given oracle access to the code. 
This is avoided by the binary code:

\begin{definition}\label{df:Balpha} Let $n\in N\setminus\{0\}$. Call an element of $N$ {\em relevant (wrt $L,n$)} 
if it equals either
\begin{enumerate}\itemsep=0pt
\item[--] $\langle S,\bar a\rangle$ for some $\bar a\in [n]^{\ar(S)}$ and $S\in L$ a relation symbol, or
\item[--] $\langle S,\bar a,i\rangle$  for some $i<|n|$ and $\bar a\in [n]^{\ar(S)}$ and $S\in L$  a function symbol.
\end{enumerate}
A set $\alpha^N\subseteq N$ is a
{\em binary code (in $N$)} of the $L$-structure 
$\str B(L,n,\alpha^N)
$ 
with universe $[n]$ if 
\begin{enumerate}\itemsep=0pt
\item[--] every relation symbol $S\in L$ is interpreted in $\str B(L,n,\alpha^N)$ by the set
$S^{[n]}$ of those $\bar a\in [n]^{\ar(S)}$ with $ \langle S,\bar a\rangle\in\alpha^N$;
\item[--] every   function symbol  $S\in L$ is interpreted in $\str B(L,n,\alpha^N)$ by the function
$S^{[n]}$ mapping $\bar a\in [n]^{\ar(S)}$ to $\min\{a,n-1\}$ for
the unique $a\in[n]$ such that for all $i<|n|$ we have $\bit(i,a)$ equal to 1 or 0 depending on whether~$\langle S,\bar a,i \rangle$ is in~$\alpha^N$ or not. 
\end{enumerate}
\end{definition}

\begin{remark} Some comments are in order:
\begin{itemize}\itemsep=0pt
\item[--] Since $(N,\alpha^N)\models\Sa$, 
there exists a unique $a$ as required. Hence, by Lemma~\ref{lem:subst}, $\str B(L,n,A)$ is well-defined for every $n\in N\setminus\{0\}$ and every $\Delta_1^b( \alpha)$-definable $A\subseteq N$.

\item[--] The minimum above is an almost arbitrary convention to ensure the right range. It can be avoided when restricting $n$
to powers of 2 as is frequently done in the context of NP search problems (e.g.~\cite{papa,bj}).

\item[--] Every set $\alpha^N\subseteq N$ such that $(N,\alpha^N)\models\Sa$ is the binary code (in $N$) of a 
unique $L$-structure on~$[n]$. In particular, every set $\alpha^\N\subseteq \N$
is the binary code (in~$\N$) of a unique $L$-structure on~$[n]$. 

\item[--] Functions can be evaluated in polynomial time with oracle access to the binary code: for a, say, unary function symbol $S\in L$ there is $\tilde S(x,y)\in\PV(\alpha)$ such that $\tilde S(a,n)$ is the value of $S$ on $a<n$ in $\str B(L,n,\alpha^\N)$.
\end{itemize}
\end{remark}



The following lemma states for all models  of $\Sa$, that the unary code is in P relative to the binary code, and that the binary code is in $\textup{NP}\cap\textup{coNP}$ relative to the unary code (see~\cite[Section~7.6]{krabuch} for complexity classes in $N$).

\begin{lemma}\label{lem:AB} 
There are a $\Delta_0^b( \alpha)$-formula $\psi_0(u,x)$ and $\Sigma_1^b( \alpha)$-formulas $\psi_1(u,x),\psi_2(u,x)$ independent of $(N,\alpha^N)$
such that for all $n\in N\setminus\{0\}$:
\begin{enumerate}\itemsep=0pt
\item[(a)] If  $A\subseteq N$ denotes the set defined by $\psi_0(u,n)$  in $(N,\alpha^N)$, then
 $\str A(L,n,A)$ is defined  and equals $\str B(L,n,\alpha^N)$.
\item[(b)] If   $\str A(L,n,\alpha^N)$ is defined and  $A\subseteq N$ denotes the set defined by $\psi_1(u,n)$ in $(N,\alpha^N)$, then  $\psi_2(u,n)$ defines
$N\setminus A$ in $(N,\alpha^N)$ and  $\str A(L,n,\alpha^N)=\str B(L,n,A)$.
\end{enumerate}
\end{lemma}

\begin{proof}[Sketch of proof.] We only sketch the definition of $\psi_1(u,x)$. It implements the following procedure:
given $(u,n)$, reject if $u$ is not relevant wrt $L,n$; else, say $u=\langle S,\bar a,i\rangle$ for $\bar a\in[n]^{\ar(S)}$, $i<|n|$ and $S\in L$ a function symbol; guess $b\in [n]$; if $\alpha(\langle S,\bar a,b\rangle)\wedge\bit(i,b){=}1$, accept; else reject.  
\end{proof}

\subsection{NP search problems from finitary combinatorial principles}\label{sec:FOsearch}

Let $L$ be a finite language disjoint from $\PV$.
Following \cite{bj} we use existential first-order $L$-sentences of a syntactically  simple form to define type 2 NP search problems. It is important to allow not only symbols from $L$ but additionally ``built-in'' symbols. For notational simplicity  we only consider built-in symbols from $\PV$:

\begin{definition}\label{df:builtin} An {\em $L$-formula with built-in $\PV$} is a $(\PV\cup L)$-formula.
\end{definition}

The difference is in the semantics: $L$-formulas with built-in $\PV$ are evaluated in $L$-structures with universe $\N$ or $[n]$ for $n\in\N\setminus\{0\}$ (up to isomorphism). On universe $\N$
the evaluation is as usual by considering the expansion interpreting the symbols from $\PV$ as in the standard model. For an $L$-structure on $[n]$ it is usual in finite model theory to consider the expansion by the graphs of $\PV$-function symbols restricted to $[n]$.
We proceed equivalently but avoid the extra symbols for the graphs. Instead 
we require that every atomic formula in which some $\PV$-function symbol $f$ occurs has 
the form $f(\bar t){=}s$ where $\bar t,s$ are $L$-terms. Such an atom expresses that $(\bar t,s)$ is in the graph of $f$.
We omit further details because, in fact, we are only interested in {\em basic} sentences, following Buss and Johnson's \cite{bj} mode of speech:

\begin{definition} \label{df:basic} An $L$-formula with built-in $\PV$ is {\em basic} if it equals \begin{equation}\label{eq:basic} 
\textstyle \exists \bar y\ \bigvee_{i\in I} \bigwedge_{j\in J}\lambda_{ij}, 
\end{equation} where $I,J$ are nonempty index sets and each $\lambda_{ij}$ is
 a literal of the form
$$R(\bar u),\ \neg R(\bar u),\ f(\bar u){=}v,\ \neg u{=}v,\textup{ or } u{=}v,
$$
where $R$ is a relation symbol and $f$ a function symbol from $L\cup\PV$, and $\bar u,u,v$ are variables.
\end{definition}

This slightly deviates from  \cite[Definition~2.9]{bj} in that there relation symbols are forbidden but constant symbols from $\PV$ are allowed within $\bar u,u,v$ above. 

To be precise how such sentences are evaluated in $L$-structures with universe $U=[n]$ or $U=\N$ we stipulate that $u{<}v$ (which is of the form $R(\bar u)$ above) defines the natural order on~$U$, and $f(\bar u){=}v$ for $r$-ary $f(\bar u)\in\PV$ defines $\{(\bar a,b)\in U^{r+1}\mid f^\N(\bar a)=b\}$.


\begin{definition}\label{df:fcp} A {\em finitary combinatorial principle (in the language $L$)} is a basic $L$-sentence with built-in $\PV$ that is {\em valid in the finite}, i.e., true in all finite $L$-structures with universe~$[n]$ for some  $n\in\N\setminus\{0\}$. Being {\em without built-in symbols} means that $\PV$-symbols do not occur. 
\end{definition}

\begin{remark}\label{rem:herbrand} Standard Herbrandization allows to compute from any $L$-formula $\varphi$ with built-in  $\PV$ an equivalid basic $L'$-formula $\varphi'$ with built-in $\PV$ where $L'$ is $L$ plus certain functions symbols. Note that a negative literal $\neg f(\bar u){=}v$ can be eliminated using $\exists y(f(\bar u){=}y\wedge \neg y{=}v)$.
In fact, $\varphi$ is true in all $L$-structures on a given universe ($\N$ or $[n]$) if and only if $\varphi'$ is true in all $L'$-structures on that universe. 
\end{remark}

 Let 
$\exists\bar y\psi(\bar y)$ be a basic  $L$-sentence with built-in $\PV$, and $\bar y=(y_0,\ldots,y_{k-1})$. 
Define
\begin{equation}\label{eq:Amodels}
\q{\str A(L,x,\alpha)\models\psi(y)}
\end{equation}
to be the quantifier free $\PV(\alpha)$-formula obtained from $\psi(\bar y)$ as follows: first, replace $L$-atoms of the form $R(\bar u)$ by $\alpha(\langle R,\bar u\rangle)$ and $ f(\bar u){=}v$ by $\alpha(\langle f, \bar u,v\rangle)$ (note $\PV$-atoms are left untouched); second, letting
$\psi'(y_0,\ldots,y_{k-1})$ denote the resulting formula, define \eqref{eq:Amodels}
to be
$$
0{<}x\to\textstyle\bigwedge_{j{<}k} (y)_j{<}x\wedge\psi'((y)_0,\ldots,(y)_{k-1});
$$
note $(y)_j$ is a $\PV$-term with variable $y$ and a constant $j\in\PV$. 

Even if $\exists \bar y\psi(\bar y)$ is valid in the finite, $\q{\str A(L,x,\alpha)\models\psi(y)}$ might not be a type 2 NP search problem. It can fail to be total (cf.~\eqref{eq:total}) since  $\alpha$ can fail to be the unary code of a some $L$-structure on~$[x]$. One can define a different total search problem: find $y$ such that  $\q{\str A(L,x,\alpha)\models\psi(y)}$ if $\str A(L,x,\alpha)$ is defined, and otherwise $y$ witnesses that~$\alpha$ is not such a code. But this property of $y$ is not verifiable in polynomial time with oracle~$\alpha$, so the search problem is not NP.
These problems disappear when using the binary code.

The formula
\begin{equation}\label{eq:Bmodels}
\q{\str B(L,x,\alpha)\models\psi(y)}
\end{equation}
is similarly defined but replacing $ f(\bar u){=}v$ (not by $\alpha(\langle f, \bar u,v\rangle)$ but instead) by a $\Delta^b_0(\alpha)$-formula
defining the graph of the interpretation of $f$ in $\str B(L,x,\alpha)$.  
The choice of this formula shall not play any further role; for example,
one might take
\begin{equation*}\label{eq:binf}
\begin{array}{rcl}
&&\Big( v{<}x\wedge \forall i{<}|x|\big(\alpha(\langle f,\bar u,i\rangle)\leftrightarrow \bit(v,i){=}1\big)\Big)\\
&&\vee\; \Big( v{=}x{-}1\wedge\exists i{<}|x|\big(  \alpha(\langle f,\bar u,i\rangle) \wedge  \bit(x{-}1,i){=}0 \\
&&\qquad 
\wedge   \
\forall j{<}|x| (  i{<}j\to  (\alpha(\langle f,\bar u,j\rangle)\leftrightarrow \bit(x{-}1,j){=}1    )   ) \big)  \Big).
\end{array} 
\end{equation*}

All formulas have the free variables shown. We employ suggestive notation for substitutions. E.g. 
$\q{\str B(L,n,f_{\bar z}^{-1}(0))\models\psi(a)}$ is obtained by
substituting  $n,a$ for $x,y$ and $f_{\bar z}^{-1}(0)$ for~$\alpha$ (see~\eqref{eq:subfct}).
The following is clear:

\begin{lemma} \label{lem:formAB} Let $(N,\alpha^N)$ be a model of $\Sa$ and $\exists y_0\cdots y_{k-1} \psi(y_0,\ldots,y_{k-1})$ be a basic $L$-formula with built-in $\PV$. Then for all  $(n,a)\in N^2$ with $n\neq 0$:
$$
N\models  \q{\str B(L,n,\alpha)\models\psi(a)}\  \Longleftrightarrow\ \str B(L,n,\alpha^N)\models\psi((a)_0,\ldots, (a)_{k-1}).
$$
If furthermore $\str A(L,n,\alpha^N)$ is defined, then
$$
N\models  \q{\str A(L,n,\alpha)\models\psi(a)}\  \Longleftrightarrow\ \str A(L,n,\alpha^N)\models\psi((a)_0,\ldots, (a)_{k-1}).
$$
\end{lemma}

If $\exists \bar y\psi(\bar y)$ is valid in the finite, then  $\q{\str B(L,x,\alpha)\models\psi(y)}$ is a type 2 NP search problem in the sense of Section~\ref{sec:searchprbl}. Indeed, the above lemma (for $N=\N$) implies totality \eqref{eq:total}, and boundedness is witnessed by  $t(x):=c\#( x\#\cdots\# x)$ with~$k$ iterations of $\#$ and suitable $c\in\N$ (by
 \eqref{eq:seqbound}).
It is the problem, given a natural $n>0$ and access to an oracle~$\alpha^\N\subseteq\N$, to find a satisfying assignment of 
$\psi(y_0,\ldots,y_{k-1})$ in $\str B(L,n,\alpha^\N)$. 

\begin{definition}\label{df:assNPSP} Let $\varphi=\exists\bar y\psi(\bar y)$ be a finitary combinatorial principle in the language $L$. The type 2 NP search problem $Q_\varphi$ {\em associated to} $\varphi$ is $\q{\str B(L,x,\alpha)\models\psi(y)}$.
\end{definition}

Here, and in similar contexts below, we silently assume that the language $L$ is finite and disjoint from $\PV$, and that $\psi(\bar y)$ is quantifier free.



\subsection{Proof of Theorem~\ref{thm:bjstrong}}\label{sec:bj}

Let $\exists\bar y\varphi(\bar y),\exists\bar w\tilde \varphi(\bar w)$ be finitary combinatorial principles in the  languages $L,\tilde L$, respectively.
Hence we have type 
2 NP search problems $\q{\str B(L, x,\alpha)\models\varphi(y)}$ and $\q{\str B(\tilde L, \tilde x,\alpha)\models\tilde\varphi(w)}$.
Let $t(x)$ and $\tilde t(\tilde x)$ be terms  witnessing their boundedness \eqref{eq:tbd}. Using the propositional translation $\langle\cdot\rangle$ of Section~\ref{sec:prop},
the totality of these search problems is naturally expressed by a sequence of propositional tautologies, one for each universe $[n]$ where $n>0$. 
We get two such sequences, one for the unary and one for the binary code of structures.  
There is some recent work~\cite{sergithesis,barny}  comparing the two translations in propositional proof complexity.

\begin{definition} Let $n\in\N\setminus\{0\}$. The {\em binary translation of $\varphi$ on $[n]$} is 
$$
\Big\langle 
  \exists y{<}t(n)\q{\str B(L, n,\alpha)\models\varphi(y)}\Big\rangle.
$$
\end{definition}

The formula $\str A(L,x,\alpha)\textit{ is defined}$ is the conjunction of 
$$
\forall \bar u{<}x\exists v{<}x\ \alpha(\langle S, \bar u,v\rangle)\wedge \forall \bar u,v,v'{<}x
\big(v{=}v'\vee \neg\alpha(\langle S, \bar u,v\rangle)\vee\neg\alpha(\langle S, \bar u,v'\rangle)\big)
$$
for every function symbol $S\in L$. This is a $\Delta_0^b(\alpha)$-formula with free variable $x$. It is satisfied by $n\neq 0$ in a model $(N,\alpha^N)$ of $\Sa$ if and only if 
$ \str A( L, n,\alpha^N)$ is defined in $(N,\alpha^N)$. 

\begin{definition}\label{df:untransl} Let $n\in\N\setminus\{0\}$. The {\em unary translation of $\varphi$ on $[n]$} is 
$$
\Big\langle \str A(L,n,\alpha)\textit{ is defined}\to
  \exists y{<}t(n)\q{\str A(L, n,\alpha)\models\varphi(y)}\Big\rangle.
$$
\end{definition}

`Propositional translation'' in Theorems~\ref{thm:bj} and \ref{thm:bjstrong} refers to the unary tanslation.

\begin{example}\label{ex:wphptransl} A basic sentence expressing the $n^2$ to $n$ weak pigeonhole principle~\eqref{eq:wphpformula} is the existential closure of
\begin{equation*}
\big(f(x,y){=}z\wedge f(x',y'){=}z\wedge \neg x{=}x'\big)\ \vee\ \big(f(x,y){=}z\wedge f(x',y'){=}z\wedge \neg y{=}y'\big).
\end{equation*}
Write $i\in[n^2]$ as  $i=i_0\cdot n+i_1$ for $i_0,i_1\in[n]$. Further write
 $p_{ij}$ for the propositional variable $\langle\alpha(\langle f,i_0,i_1,j\rangle)\rangle$ where $i\in[n^2],j\in[n]$.
The unary translation on $[n]$ has many occurrences of the Boolean constants 0,1. If one eliminates these occurrences by repeatedly replacing subformulas $0\vee F, 1\wedge F$ by $F$ etc., then one gets the  familiar disjunction of 
$$
\begin{array}{lcl}
\textstyle\bigwedge_{j<n}\neg p_{ij}&&i\in[n^2],\\
 p_{ij}\wedge p_{ij'}&&j,j'\in[n],j\neq j',\\
 p_{ij}\wedge p_{i'j}&&j\in[n],i,i'\in[n^2], i\neq i',
\end{array}
$$
with multiple occurrences of the last disjuncts.
\end{example}

 \begin{remark}The unary translation is very similar to the propositional translation used by Buss and Johnson~\cite{bj}. More precisely,  the translation in \cite[Definition~3.2]{bj} produces a sequent $F\Rightarrow G$; if one eliminates Boolean constants as indicated in the example above both in $(\neg F\vee G)$ and in our  unary translation, then one obtains the same formula. 
\end{remark}
 
A  {\em substitution instance} of a propositional formula is  obtained by simultaneously replacing some of its variables by propositional  formulas. The first statement of the following is a slightly more detailed statement of Theorem~\ref{thm:bjstrong}. 

\begin{theorem} If $\q{\str B(L,x,\alpha)\models\varphi(y)}$ is a consequence of $\q{\str B(\tilde L, \tilde x,\alpha)\models\tilde \varphi(w)}$ over $\BAo$, then there are $d,n_0\in\N$ such that for all $n>n_0$ there
are size $2^{|n|^d}$ depth~$d$ Frege proofs of the unary translation of $\varphi$ on $[n]$ from substitution instances of the unary translations 
of $\tilde\varphi$ on~$[\tilde n]$ for all $\tilde n<2^{|n|^d}$.

The same holds for the binary translations of $\varphi$ and $\tilde\varphi$.
\end{theorem}

\begin{proof} 
Assume $\q{\str B(L,x,\alpha)\models\varphi(y)}$ is a consequence of $\q{\str B(\tilde L, \tilde x,\alpha)\models\tilde \varphi(w)}$ over~$\BAo$. 
Recall  $t(x),\tilde t(\tilde x)$ are terms witnessing the boundedness of these search problems.
By compactness there is a finite $\Delta\subseteq\PV(\alpha)$
such that $ \BAo$ proves
\begin{equation}\label{eq:conseq}
\textstyle 
\bigwedge_{f(z,\bar z)\in\Delta}\forall \tilde x \bar z\exists w{<}\tilde t(\tilde x)\ \q{\str B(\tilde L,\tilde x,f^{-1}_{\bar z}(0))\models\tilde\varphi(w)}\ \to\ \exists y{<}t(x)\ \q{\str B(L,x,\alpha)\models\varphi(y)}.
\end{equation}

Let $\psi_0(u,x)$ be the formula from Lemma~\ref{lem:AB}.

 \medskip
 
 \noindent{\em Claim 1.} For every $f(z,\bar z)\in\PV(\alpha)$, $\BAo$ proves 
\begin{equation}\label{eq:ant}
\begin{split}
\textstyle 
  &\Big(\big( \str A(\tilde L,\tilde x,\alpha)\textit{ is defined}\to
  \exists w{<}\tilde t(\tilde x)\q{\str A(\tilde L,\tilde x,\alpha)\models\tilde\varphi(w)}\big)\big[\alpha/\psi_0(\cdot,\tilde x)\big]\Big)\big[\alpha/f^{-1}_{\bar z}(0)\big]\\
 & \to
\exists w{<}\tilde t(\tilde x)\q{\str B(\tilde L,\tilde x,f^{-1}_{\bar z}(0))\models\tilde\varphi(w)}.
\end{split}
\end{equation}

\noindent{\em Proof of Claim 1:} 
By Proposition~\ref{prop:cons}, models of $\BAo$ have the form $\langle M,\alpha^M\rangle$ where $(M,\alpha)\models\BAa$.
Suppose the assignment of $\tilde n,\bar c$ to $\tilde x,\bar z$ falsifies the succedent of~\eqref{eq:ant} in $\langle M,\alpha^M\rangle$, i.e.,
$\str B(\tilde L,\tilde x,f^{-1}_{\bar z}(0))\not\models\exists\bar w\tilde\varphi(\bar w)$ by Lemma~\ref{lem:formAB}.
We have to show that~$\tilde n,\bar c$  falsify the antecedent of~\eqref{eq:ant} in $\langle M,\alpha^M\rangle$.

Let  $A:=\{a\in M\mid f^M(a,\bar c)=0\}$. Then $A$ is $\Delta^b_1(\alpha)$-definable in $(M,\alpha^M)$ by Lemma~\ref{lem:falphaf}, so 
$(M,A)\models\BAa$ by Lemma~\ref{lem:subst}. 
Writing~$B$ for the set defined by $\psi_0(u,\tilde n)$ in $(M,A)$, Lemma~\ref{lem:AB} gives 
that $\str A(\tilde L,m,B)$ is defined (in $(M,A)$) and equals $\str B(\tilde L,\tilde n,A)$. Thus
 $$
 (M,A)\not\models\big( \str A(\tilde L,\tilde n,\alpha)\textit{ is defined}\to
  \exists w{<}\tilde t(\tilde x)\q{\str A(\tilde L,\tilde n,\alpha)\models\tilde\varphi(w)}\big)\big[\alpha/\psi_0(\cdot,\tilde n)\big],
 $$
by Lemma~\ref{lem:formAB}. Then $\tilde n,\bar c$ falsify the antecedent of \eqref{eq:ant}  in $\langle M,\alpha^M\rangle$.
\hfill$\dashv$\medskip

For every $f(z,\bar z)\in\Delta$,
the antecedent of $\eqref{eq:ant}$ is $\PVa$-provably equivalent to a 
$(\PV\cup\{\alpha\})$-formula. This formula is obtained by substituting 
 atoms $f(t,\bar z){=}0$ by suitable  $\Sigma^b_1(\alpha)$-formulas obtained from  $\Sigma^b_1(\alpha)$-definitions of the graph of $f$ (see  Lemma~\ref{lem:falphaf}).  
 Let
  $$
  \chi_0(\tilde x,\bar z_0),\ldots,  \chi_{|\Delta|-1}(\tilde x,\bar z_{|\Delta|-1})
  $$ 
enumerate the $(\PV\cup\{\alpha\})$-formulas thus obtained.
 By conservativity (Proposition~\ref{prop:cons})
 \begin{equation}\label{eq:alphaimpl}
\textstyle 
  \BAa\ \vdash\ \bigwedge_{i< |\Delta|}\forall \tilde x \bar z_{i}\chi_i(\tilde x,\bar z_i)\ \to\ \exists y{<}t(x)\q{\str B(L,x,\alpha)\models\varphi(y)}.
\end{equation}

Let $\psi_1(u,x)$ be the formula  from 
Lemma~\ref{lem:AB}.

\medskip

\noindent{\em Claim 2.} $\BAa$ proves 
\begin{equation}\label{eq:prvA}
\begin{split}&
\textstyle   \bigwedge_{i< |\Delta|}\forall \tilde x \bar z_{i} \chi_i(\tilde x,\bar z_i)\big[ \alpha/\psi_1(\cdot,x)  \big] \\
  &
   \to
  \big( \str A(L,x,\alpha)\textit{ is defined}\to\exists y{<}t(x)\q{\str A(L,x,\alpha)\models\varphi(y)}\big).
  \end{split}
\end{equation}

\noindent{\em Proof of Claim 2:}
Suppose $(M,\alpha^M)\models\BAa$ and $n\in M$ falsifies the succedent of~\eqref{eq:prvA} in~$(M,\alpha^M)$. 
Then $n\neq 0$, $\str A(L,n,\alpha^M)$ is defined and   $\str A(L,n,\alpha^M)\not\models\exists\bar y\varphi(\bar y)$ by Lemma~\ref{lem:formAB}. 
Let $A\subseteq M$ be defined by $\psi_1(u,n)$ in $(M,\alpha^M)$. By Lemma~\ref{lem:subst},
 $(M,A)\models\BAa$.
By Lemma~\ref{lem:AB}, $\str B(L,n,A)$ equals $\str A(L,n,\alpha^M)$, so
 $(M,A)\not\models\exists y{<}t(x)\q{\str B(L,n,\alpha)\models\varphi(y)}$ by Lemma~\ref{lem:formAB}.
  Thus $(M,A)$ falsifies the antecedent of \eqref{eq:alphaimpl}, so $n$ falsifies the antecedent of~\eqref{eq:prvA} in $(M,\alpha^M)$.
\hfill$\dashv$\medskip

Parikh's theorem (see e.g.~\cite[Theorem~1.4.3]{busshand}) allows to bound $\forall \tilde x \bar z_{i}$ in \eqref{eq:prvA} by 
 a $\PV$-term $s(x)$. 
Thereby we get a $\Sigma^b_\infty(\alpha)$-formula and can apply 
Proposition~\ref{prop:simulation}. This yields for every natural $n>0$ a  quasipolynomial (in~$n$) size bounded depth  Frege proof of
the unary translation of $\varphi$ on $[n]$
from  the formulas
$\textstyle 
\left\langle\chi_i(\tilde n,\bar c_i)\big[ \alpha/\psi_1(\cdot,n)  \big]\right\rangle
$
where $ \tilde n,\bar c_i<s(n),i<|\Delta|$. These formulas
are substitution instances of 
the unary translation of $ \tilde\varphi$ on $[\tilde n]$.

\medskip

The proof of the second statement is similar but  simpler:
 from \eqref{eq:conseq} move to a $(\PV\cup\{\alpha\})$-formula by substituting definitions for the graphs of the functions in  $\Delta$. Then bound the quantifiers $\forall\tilde x\bar z$ using Parikh's theorem and apply the simulation (Proposition~\ref{prop:simulation}).
\end{proof}


\section{Finitary combinatorial principles}\label{sec:comb}

From a computational perspective it is natural to view  a finitary combinatorial principle as a search problem as in Definition~\ref{df:assNPSP}. From a more logical perspective one might  think of it as a reasoning rule that allows to infer the existence of certain configurations in finite structures. The interesting case is when the principle fails in some infinite structure, so the rule is sound only in the finite. 
 It is not obvious how to compare the logical strength of such principles in the finite since they all hold in the same (all) finite structures. The crucial observation is that they might behave differently with respect to {\em partial} finite structures, allowing the distinction between {\em weak} and {\em strong} principles. Intuitively, a  principle is weak if seeing only a small fraction of a given structure is already sufficient to verify its truth. We shall verify later that the thus distinguished logical strength of principles implies distinct computational complexities of the associated type 2 NP search problems.

We define partial structures and their logic in Section~\ref{sec:partial}, and their codes by partial oracles in Section~\ref{sec:partialcodes}.
Weak and strong principles are defined in Section~\ref{sec:weakstrong} and examples are discussed in Section~\ref{sec:exas}. Section~\ref{sec:dense} establishes the combinatorial lemmas for the forcing constructions to come.

\subsection{Partial structures}\label{sec:partial}

Let $L$ be a language.
For the sake of exposition, let us agree that the interpretation $S^A$ of a symbol $S\in L$ in an $L$-structure $\str A$ with universe $A$ is a function from $A^{\ar(S)}$ into $A$ or into~$\{0,1\}$ depending on whether~$S$ is a function or a relation symbol. For relation symbols we identify $S^{A}$ with $\{\bar a\in A^{\ar(S) } \mid S^A(\bar a)=1\}$. 

A {\em partial $L$-structure} $\str A$ is similarly explained but allowing value $1/2$ which we read as ``undefined'' and assume to be outside $A$. That is, the interpretation $S^{A}$ for $S\in L$ is a function from $A^{\ar(S)}$ into
$A\ \dot\cup\ \{1/2\}$ or into $\{0,1,1/2\}$ depending on whether~$S$ is a function or a relation symbol. $\str A$ is {\em total}
if $S^A(\bar a)\neq 1/2$ for all $S\in L$ and all $\bar a\in A^{\ar(S)}$.

Let $\str A,\str B$ be partial $L$-structures with universes~$A,B$ respectively.
Then $\str B$ is a {\em partial substructure of} $\str A$ if $B\subseteq A$ and interpretations $S^{B}$ are obtained from $S^A$ by changing
some values to~$1/2$; it is {\em induced} if for every
 $S\in L$ and all $\bar b\in B^{\ar(S)}$ we have $S^B(\bar b)=S^{A}(\bar b)$ except for the case that
$S$ is a function symbol and $S^A(\bar b)\not\in B$; in this case $S^B(\bar b)=1/2$. 
We say $\str A$ {\em extends} a partial substructure $\str B$ if $A=B$. An {\em isomorphism} from $\str A$ onto $\str B$ is a bijection $\pi$ from~$A\cup\{0,1,1/2\}$ onto $B\cup\{0,1,1/2\}$ which is the identity on $\{0,1,1/2\}$ and such that 
$\pi\circ S^A=S^B\circ\pi $ for all $S\in L$; here, we assume $\{0,1,1/2\}\cap(A\cup B)=\emptyset$.
An {\em embedding} from $\str B$ into~$\str A$ is an isomorphism from $\str B$ onto a partial substructure of $\str A$.

\begin{definition} \label{df:size} The {\em size} of $\str A$ is
\begin{equation*}\textstyle
\sum_{S\in L}|\{\bar a\in A^{\ar(S)}\mid S^A(\bar a)\neq 1/2\}|.
\end{equation*}
We let $s_L(n)$ denote the size of a total $L$-structure with a universe of cardinality $n$, that is,
\begin{equation*}
\textstyle
s_L(n):=\sum_{S\in L}n^{\ar(S)}.
\end{equation*}
\end{definition}

We explain how to evaluate formulas in a partial $L$-structure $\str A$. We silently extend all~$S^A$
to domain $(A\cup\{1/2\})^{\ar(S)}$ giving value $1/2$ to all new argument tuples, i.e., $S^A(\bar a):=1/2$ if $\bar a\in(A\cup\{1/2\})^{\ar(S)}\setminus A^{\ar(S)}$. Then the interpretation $t^A$ of a 
 closed $L$-term $t$ (i.e., $t$ has no variables) with parameters from $A$ is defined as 
usual by composition of the interpretation of its function symbols. That is, values of closed terms are computed bottom-up as usual  but  upon encountering the value $1/2$ the computation is aborted with output~$1/2$.

For an $L$-sentence $\varphi$ with parameters from $A$ 
we define the {\em truth value $v^{\str A}(\varphi)\in \{0,1,1/2\}$ of $\varphi$ in $\str A$} in a way
 familiar from 3-valued logic (see e.g.~\cite{partial}):
\begin{itemize}\itemsep=0pt
\item[--] If $\varphi$ has the form $t{=}s$  for closed $L$-terms $t,s$ with parameters from $A$, then
$v^{\str A}(\varphi):=1/2$ if
at least one of $t^A,s^A$ equals $1/2$; otherwise, $v^{\str A}(\varphi)$ is 1 or 0 depending on whether 
$t^A$ equals $s^A$ or not.
\item[--] If $\varphi=S(t_0,\ldots,t_{\ar(R)-1})$ for closed $L$-terms $t_0,\ldots, t_{\ar(S)-1}$ with parameters from $A$ 
and $S\in L$ a relation symbol, then  $v^{\str A}(\varphi):=S^A(t^A_0,\ldots, t^A_{\ar(S)-1})$.
\item[--]  If $\varphi=\neg\psi$, then $v^{\str A}(\varphi):=1-v^{\str A}(\psi)$.
\item[--] If $\varphi=(\psi\wedge\chi)$, then $v^{\str A}(\varphi):=\min\{v^{\str A}(\psi),v^{\str A}(\chi)\}$.
\item[--] If $\varphi=\forall x\psi(x)$, then $v^{\str A}(\varphi):=\min\{v^{\str A}(\psi(a))\mid a\in A\}$.
\end{itemize}
We consider formulas as built from atomic formulas using~$\neg,\wedge,\forall x$ and view $(\varphi\vee\psi)$ and $\exists x\varphi$ as abbreviations of $\neg(\neg \varphi\wedge\neg\psi)$ and $\neg\forall x\neg\varphi$, respectively. Then
\begin{eqnarray*}
v^{\str A}(\varphi\vee\psi)&=&\max\big\{v^{\str A}(\varphi),v^{\str A}(\psi)\big\},\\
v^{\str A}(\exists x\varphi(x))&=&\max\big\{v^{\str A}(\varphi(a))\mid a\in A\big\}.
\end{eqnarray*}

\begin{definition}
A partial structure $\str A$ {\em verifies $\varphi$} if $v^{\str A}(\varphi)=1$; it {\em falsifies} $\varphi$ if it verifies~$\neg\varphi$. 
\end{definition}

Clearly, if a partial structure $\str A$ extends $\str B$, then it verifies every sentence which is verified by $\str B$.
A total structure $\str A$ verifies $\varphi$ if and only if $\str A\models\varphi$.

\begin{lemma}\label{lem:pres} Let $\str A$ be a partial structure and $\str B$ a partial substructure of $\str A$. Then every existential sentence verified by $\str B$ is verified by $\str A$.
\end{lemma}

\begin{proof} 
Call a sentence $\varphi$ with parameters from $B$ {\em good} if $v^{\str B}(\varphi)=1/2$ or $v^{\str B}(\varphi)=v^{\str A}(\varphi)$.
The set of good sentence contains all atomic formulas and is closed under $\wedge$ and $\neg$, so contains all quantifier free sentences with parameters from~$B$. 

If $\str B$ verifies $\exists \bar x\varphi(\bar x)$ for quantifier free $\varphi(\bar x)$, then it verifies $\varphi(\bar b)$ for some tuple $\bar b$ from~$B$. Since~$\varphi(\bar b)$ is good, also $\str A$ verifies $\varphi(\bar b)$ and hence $\exists \bar x\varphi(\bar x)$.
\end{proof}

\subsection{Partial codes}\label{sec:partialcodes}

As structures are coded by oracles, partial structures are coded by ``partial oracles''. 
As in Section~\ref{sec:coding}, we fix finite languages  $L,\tilde L\subseteq\N$. We further fix a model $(N,\alpha^N)$ of $\Sa$, 
so $\langle N,\alpha^N\rangle\models\So$ (Proposition~\ref{prop:cons}). We do not distinguish
 between symbols in~$\PV(\alpha)$ and their interpretations in $\langle N,\alpha^N\rangle$. 
 We also blur the distinction between $p\in N$ and the set 
it codes, namely the set of $a\in N$ with $\bit(p,a)=1$. Recall that {\em relevant} elements of $N$ are those used to code structures (see Definition~\ref{df:Balpha}).

\begin{definition}    Let $n\in N\setminus\{0\}$.
Let $p\in N$ be such that $p=\langle p_0,p_1\rangle$ for certain $p_0,p_1\in N$.
Such a $p$ is a {\em partial $L$-oracle on $[n]$} if $p_0$ and $p_1$ code disjoint sets 
of relevant (wrt~$L,n$) elements
such that for every function symbol~$S\in L$ and $\bar a\in[n]^{\ar(S)}$ either
all or none of $\langle S,\bar a,i\rangle,i<|n|,$ are elements of $p_0\cup p_1$.
Then $p$ codes the following partial structure 
$$
\str B(p)=\str B(L,n,p)
$$ with universe $[n]$:
\begin{enumerate}\itemsep=0pt
\item[--]
 for  a function symbol   $S\in L$
we have $S^{[n]}(\bar a)=1/2$  if $p_0\cup p_1$ does not contain $\langle S,\bar a,i\rangle$ for some  (equivalently all) $i<|n|$; otherwise 
$S^{[n]}(\bar a)=\min\{a,n-1\}$ for the unique $a\in[n]$  with $\langle S,\bar a,i\rangle\in p_{\bit(i,a)}$ for all $i<|n|$; 

\item[--] for a relation symbol $S\in L$
we have $S^{[n]}(\bar a)$ equal to $0$ if $\langle S,\bar a\rangle\in p_0$, equal to 1 if $\langle S,\bar a\rangle\in p_1$, and equal to $1/2$ if $\langle S,\bar a\rangle\not\in p_0\cup p_1$. 
\end{enumerate}
\end{definition}

In the standard $\PV$-model $N=\N$, one might call an $\tilde L$-structure $\str C$ on $[m]$ ``implicitly feasible in'' $\str B(L,n,\alpha^\N)$ if a binary code of $\str C$ is polynomial time Turing reducible to $\alpha^\N$. These are precisely the structures $\str C$ of the form 
$\str B(L,n,f^{-1}(0))$ for some $f\in\PV(\alpha)$. 
 It shall be convenient to work instead with a presentation of such structures
(see Lemma~\ref{lem:family} below)
given by a family of decision trees computing the interpretations of the symbols in $\tilde L$.

Recall, Definition~\ref{df:tree} defines sequences of $\alpha$-answers  to decision trees $t$. The mode of speech for partial oracles is analogous:

\begin{definition}\label{df:pansw} Let $n\in N\setminus\{0\}$,
 $p$ a partial $L$-oracle on $[n]$ and $t(\bar x,z)$ a decision tree in~$N$. Then
$c\in N\setminus\{0\}$ is a {\em sequence of $p$-answers to $t$ on $\bar a$} if
for all $i<|c|-1$ we have $t(\bar a,c_{<i})$ is odd and:
 \begin{enumerate}\itemsep=0pt
\item[--] $\bit(i,c)=1$ and $ \lfloor t(\bar a,c_{<i})/2\rfloor\in p_1$, or 
\item[--] $\bit(i,c)=0$ and $\lfloor t(\bar a,c_{<i})/2\rfloor\in p_0$, or
\item[--] $\bit(i,c)=0$ and $ \lfloor t(\bar x,c_{<i})/2\rfloor$ is not relevant (wrt $L,n$).
\end{enumerate}
It is {\em complete} if $t(\bar a,c)$ is even; it is
{\em maximal} if it is either complete or~$t(\bar a,c)$ is odd and~$\lfloor t(\bar a,  c)/2 \rfloor$ is relevant and outside $p_0\cup p_1$.
\end{definition}

\begin{definition}\label{df:family} For each $\tilde S(\bar x)\in\tilde L$ let $t_{\tilde S}(\bar x,z)$ be a decision tree of height $h_{\tilde S}(\bar x)$ in $N$.
For $m,n\in N\setminus\{0\}$ and a partial $L$-oracle $p$ on $[n]$ we get a partial $\tilde L$-structure 
$$
\str C((t_{\tilde S})_{\tilde S\in\tilde L},m,p)
$$ 
with universe $[m]$ as follows. For 
$\tilde S\in\tilde L$ and $\bar a\in[m]^{\ar(\tilde S)}$ let $\tilde S^{[m]}(\bar a)\neq 1/2$ 
if and only if  there is exactly one  complete
sequence $c$ of $p$-answers to~$t_{\tilde S}$ on $\bar a$; then
$$
\tilde S^{[m]}(\bar a):=\left\{\begin{array}{ll} \min\{ t_{\tilde S}(\bar a,c)/2,m-1\}&\text{if $\tilde S$ is a function symbol,}\\ 
\min\{ t_{\tilde S}(\bar a,c)/2,1\}&\text{if $\tilde S$ is a relation symbol.}\ \end{array}\right.
$$
For $\alpha^N \subseteq\N$ we define $\str C((t_{\tilde S})_{\tilde S\in\tilde L},m,\alpha^N)$
 analogously  using sequences of  $\alpha^N$-answers.
\end{definition}

The minima above are just a convention to ensure the right range. Of course, in the standard $\PV$-model there can only be at most one complete sequence of $p$-answers. In our possibly nonstandard model $N$, this holds if the  decision trees have a sufficiently simple definition like the following.

\begin{definition}\label{def:givenbyterms} A family $(t_{\tilde S})_{\tilde S\in\tilde L}$ of decision trees in~$N$ is {\em given by terms} if 
every $t_{\tilde S}(\bar x),\tilde S\in\tilde L,$ is the interpretation (in $N$) of some $\PV$-term with parameters from $N$, and has height $|h_{\tilde S}(\bar x)|$ for some  $\PV$-term $h_{\tilde S}(\bar x)$ with parameters from $N$.
\end{definition}

\begin{lemma} \label{lem:family}
Let $m\in N\setminus\{ 0\}$, $f(z,\bar z)\in\PV(\alpha)$, and $\bar a$ a tuple from $N$. Then there is a family $(t_{\tilde S})_{\tilde S\in\tilde L}$ of decision trees 
 in~$N$  given by terms such that $$\str B(\tilde L,m,f_{\bar a}^{-1}(0))=\str C((t_{\tilde S})_{\tilde S\in\tilde L},m,\alpha^N).$$
\end{lemma}

\begin{proof} It is easy to see, and also follows from Lemma~\ref{lem:subst}, that $\str B(L,m,f_{\bar a}^{-1}(0))$ is well defined in $\langle N,\alpha^N\rangle$.
Let $\tilde S\in\tilde L$ be a function symbol (the case of a relation symbol is similar).  Consider the following algorithm with oracle $\alpha$ and parameters $\bar z,m$ from~$\N$: on input $\bar x\in[m]^{\ar(\tilde S)}$ compute the length $|m|$ binary string whose  $i$-th bit is 1 or 0 depending on whether
 $f(\langle\tilde S,\bar x,i\rangle, \bar z)=0$ or not; finally output the number with this binary expansion if it is in~$[n]$, otherwise output $n-1$. Now 
choose $t_{\tilde S},h_{\tilde S}$  according Lemma~\ref{lem:falphaf}.
\end{proof}

\subsection{Weak and strong principles}\label{sec:weakstrong}

Let $L$ be a finite language.
We  define simple model-theoretic notions for an $L$-sentence $\varphi$ to be weak or strong. The case of interest is when $\varphi$ is basic (Definition~\ref{df:basic}), valid in the finite and fails in some infinite model.

\begin{definition}\label{df:d}
Let $\varphi$ be an $L$-sentence.
The {\em determinacy of}~$\varphi$ is the function $d:\N\setminus\{0\}\to \N$ such that $d(n)$ is the minimal
$m\in\N$ such that every partial $L$-structure with universe of cardinality $n$ and size at least $ m$ verifies $\varphi$.
%
If $s_L(n)\ge n^{\Omega(1)}\cdot d(n)$, then we say $\varphi$ is {\em weak}.
\end{definition}

Observe that $d(n)>0$ because there is no sentence verified by the completely undefined structure. We have $d(n)\le s_L(n)$ if and only if $\varphi$ is valid in the finite, and otherwise $d(n)=s_L(n)+1$. Intuitively,  the smaller the determinacy  the weaker the principle (i.e., the claim that it has no finite models).

\begin{remark}\label{rem:weakPV} The same definitions apply to $L$-sentences with built-in $\PV$ understanding verification as follows:  a  partial $L$-structure $\str B$ with universe $[n]$ for some $n\in\N\setminus\{0\}$ {\em verifies} a $(\PV\cup L)$-sentence if and only if so does the partial $(\PV\cup L)$-structure that interprets the symbols from $L$ as $\str B$ and the symbols from $\PV$ as  the partial substructure induced on~$[n]$ in the standard $\PV$-structure $\N$.
It is easy to check that for basic sentences verification coincides with truth as explained in Section~\ref{sec:FOsearch} (after Definition~\ref{df:basic}).
\end{remark}


\begin{definition}\label{df:g}
Let  $g:\N\setminus\{0\}\to\N$, and $\str B$ be an infinite (total) $L$-structure with universe~$B$.
An induced partial substructure $\str B_0$ of $\str B$ with finite universe $B_0$ is {\em $g$-large} if there exists a subset  $V\subseteq B\setminus B_0$ of  size at most $ g(|B_0|)$ such that for every function symbol $S\in L$ 
the interpretation~$S^B$ of~$S$ in $\str B$ maps $B_0^{\ar(S)}$ into $B_0\cup V$.
The structure $\str B$ is {\em $g$-large} if every finite partial substructure of $\str B$ embeds into a $g$-large partial substructure of $\str B$ with a universe of the same cardinality.

An $L$-sentence  is {\em strong} if its negation has an infinite $n^{o(1)}$-large model.
 \end{definition}

Assume $\str B\not\models\varphi$ where $\varphi$ is basic and valid in the finite.
Then no finite subset of $B$ is closed under the interpretations of the function symbols in $\str B$. 
 Definition~\ref{df:g} quantifies how many function values are outside a given finite subuniverse. Intuitively, the smaller~$g$, the closer $\varphi$ is to be satisfiable in the finite; hence, the smaller $g$, the stronger the principle.
\medskip

Our aim is to verify these intuitions to some extent, namely in the sense of Theorem~\ref{thm:main}.
The proof requires some fair amount of work, and before getting there we consider 

\subsection{Examples}\label{sec:exas}

We start with common pigeonhole principles.

\begin{example}\label{ex:php1} Let $L:=\{f,c\}$ for a unary function symbol $f$ and a constant $0$. The {\em ($n$ to $n{-}1$) pigeonhole principle $\PHP$} is the existential closure of
$$
(f(x){=}u\wedge f(y){=}u\wedge \neg x{=}y)\vee (f(x){=}u\wedge c{=}u).
$$
This is a basic (Definition~\ref{df:basic}) variant of $(f(x){=}f(y)\wedge \neg x{=}y)\vee f(x){=}c$.
It has maximal determinacy $d(n)=s_L(n)=n+1$. It is not weak and it is strong, indeed, its negation has a  1-large model.
\end{example}

\begin{proof}
To prove the second statement, let $\str A$ have universe $A:=\N$ and interpret $c$ by $0$ and $f$ by the successor function. 
Let $\str A_0$ be a partial substructure of $\str A$ with universe~$A_0$ of cardinality $n$. Map the minimal element of $A_0$ to 0, the second largest element of $A_0$ to 1 and so on. This embeds $\str A_0$ into the partial substructure induced on $[n]$ in $\str A$. This partial substructure is $1$-large witnessed by
$V:=\{n\}$. The first statement follows noting that this partial substructure
has size $s_L(n)-1=n$ and, of course, does not verify the principle. 
\end{proof}

For readability we write our principles from now on not  in basic form as in Definition~\ref{df:basic} but  allowing ourselves atoms with more than one symbol of the language.

\begin{example}\label{ex:bphp} Let $L=\{f,g,c\}$  for unary function symbols $f,g$ and a constant $c$. Following~\cite{bj}, let the {\em onto pigeonhole principle} $\OPHP$ be the existential closure of
$$
\neg g(c){=}c\ \vee\ f(x){=}c \ \vee\ \neg g(f(x)){=}x\ \vee\ \big(\neg c{=}x\wedge\neg f(g(x)){=}x\big),
$$
and the  {\em left pigeonhole principle} $\LPHP$ is the same with the last disjunct deleted. 

Both  principles have maximal determinacy $d(n)=s_L(n)=2n+1$, are not weak and are strong, indeed, their negations have 1-large models.
\end{example}

\begin{proof} Expand the structure  $\str A$ of the previous example letting  $g^A$ be the predecessor function (understanding $g^A(0)=0$). Then argue as there. 
\end{proof}

\begin{example}\label{ex:php3}  Let $L:=\{f\}$ for a binary function symbol $f$. The {\em $n^2$ to $n$ weak pigeonhole principle} 
$\WPHP$ is defined in Example~\ref{ex:wphptransl}.

It has determinacy $d(n)=\sqrt{s_L(n)}+1=n+1$. It is weak and not strong.
\end{example}

\begin{proof} Note $s_L(n)=n^2$. It is clear that once a structure on a universe with cardinality $n$ has $n+1$ values distinct from $1/2$, then there is a collision and the principle is verified. It is also clear that there are partial structures of size $n$ not verifying the principle. 

To see $\WPHP$ is not strong, let  $\str A$ be a model of its negation. Restricted on a set $A_0$ of~$n$ points $f^A$ takes at least $n^2-n$ many values outside $A_0$. Hence, a set $V$ from the definition of a large partial substructure must have cardinality at least $n^2-n$.
\end{proof}

\begin{example}\label{ex:php2} Let $L:=\{f,g\}$ for  unary function symbols $f,g$. The {\em $2n$ to $n$ weak pigeonhole principle} $\WPHP'$ is the existential closure of
$$
(f(x){=}f(y)\wedge \neg x{=}y)\ \vee\ (g(x){=}g(y)\wedge \neg x{=}y)\ \vee\ f(x){=}g(y).
$$
It has determinacy $d(n)=s_L(n)/2+1=n+1$. It is neither weak nor strong.
\end{example}

\begin{proof}
Note $s_L(n)=2n$. In any partial $L$-structure with $n+1$ elements where $f$ or $g$ is defined (i.e., value $\neq 1/2$) two of these values are equal, so the principle is verified. Hence $d(n)\le n+1$. But $d(n)>n$ because there are size $n$ partial structures on $[n]$ that do not verify the principle, e.g., interpret $f$ by a permutation and let $g$ be completely undefined.

That $\WPHP'$ is not strong can be seen as in the previous example.
\end{proof}

\begin{example} \label{ex:rphp} The provably total (type 1) NP search problems of Je\v r\'abek's \cite{japx} theory of approximate counting $\mathsf{APC}_1$ are many-one reducible to the {\em $n$ to $n^2$ retraction pigeonhole principle} $\rPHP$ (see~\cite[Proposition~1.14]{jdual}) for $\PV$-functions. To express it by  a first-order formula over universe $[n]$ we take $L=\{g,f_0,f_1\}$ for a  binary function symbol $g$ and two unary function symbols $f_0,f_1$ and state that $g$ does not witness  that $x\mapsto (f_0(x),f_1(x))$ is a surjection from~$[n]$ onto~$[n]^2$: the existential closure of
\begin{equation*}\label{eq:rphp}
\neg f_0(g(x,y)){=}x\ \vee\ \neg f_1(g(x,y)){=}y.
\end{equation*}
It is neither weak nor strong.\end{example}

\begin{proof} A partial structure on $[n]$ that interprets $g$ by an arbitrary binary function and 
has~$f_0,f_1$ completely undefined does not verify $\rPHP$ and has size $n^2$. Since $s_L(n)=n^2+2n$, this shows that $\rPHP$ is not weak. 
That it is not strong is seen as in Example~\ref{ex:php3}.
 \end{proof}

We turn to other principles.

\begin{example}\label{ex:par} Let $L=\{f\}$ for a unary function symbol $f$. The {\em parity principle} $\PAR$ states that involutions have fixed points and is valid in structures of odd finite size: the existential closure of
\begin{equation*}\label{eq:par}
\neg x{=}f(f(x))\ \vee \  x{=}f(x).
\end{equation*}
It has determinacy 
$$
d(n)=\begin{cases}n&\text{if $n$ is odd,}\\0&\text{else.}\end{cases}
$$
It is not weak and it is strong, indeed, its negation has a 1-large model.
\end{example}

\begin{proof} 
Let $\str A$ have universe $A=\N$ and let $f^\N$ map even $n$ to $n+1$, and odd $n$ to $n-1$.
Then $\PAR$ fails in $\str A$. It is easy to see that any partial substructure of $\str A$ of 
size~$n$ embeds into the partial substructure induced in $\str A$ on $[n]$. For even $n$ this substructure is total, and for odd $n$, only the last point $n-1$ is mapped to something outside. Our claims follow.
\end{proof}

\begin{example}\label{ex:hop} Let $L:=\{f,\prec\}$ for a unary function symbol $f$ and a binary relation symbol~$\prec$ (with infix notation) and constants $0,1$. The {\em Herbrandized ordering principle}~$\HOP$ negates the Skolemized infinity axiom
stating ``$\prec$ is a partial order without a minimal element'':
the existential closure of
\begin{equation*}\label{eq:hop}
x{\prec} x\ \vee\ (x{\prec} y\wedge y{\prec}z\wedge \neg x{\prec} z) \ \vee \ \neg f(x){\prec} x.
\end{equation*}
It has maximal determinacy $d(n)=s_L(n)=n^2+n$. It is not weak and it is strong, indeed, its negation has a 1-large model.
\end{example}

\begin{proof} To prove the second statement, let $\str A$ have universe $A:=\N$, interpret $\prec$ by the inverse natural order, i.e., $\prec^A:=\{(i,j)\mid j<i\}$, and $f$ by the successor function. Every partial substructure with universe of cardinality $n$ embeds into the partial substructure induced on~$[n]$ which is $1$-large witnessed by $V:=\{n\}$. 

The claim about determinacy follows noting that the described 1-large partial substructure of $\str A$ has size $n^2+n-1$, namely, it has only one value $1/2$ (taken by $f^{A}$ on $n-1$). 
\end{proof}

\begin{remark}\label{rem:hopvar} The same reasoning applies to weaker variants \cite{bkt,atstha} of $\HOP$ adding disjuncts saying that $\prec$ is not linear and/or $f$ does not map all points to immediate predecessors:
$(\neg x{\prec} y\wedge\neg y {\prec} x\wedge \neg x{=}y)$ and/or $(f(x){\prec} y\wedge y{\prec} x)$.
\end{remark}

\begin{example} Let $L:=\{P,s,\prec,\textit{min},\textit{max}\}$ where $P$ is a unary and $\prec$ a binary relation symbol, $s$ is a unary function symbol, and $\textit{min},\textit{max}$ are constants. The {\em Induction principle} \IND\ states induction for the predicate $P$ on a discrete linear order $\prec$ with minimum $\textit{min}$, maximum $\textit{max}$ and successor $s$: the existential closure of
\begin{eqnarray*}
&&x{\prec} x\ \vee\ (x{\prec} y\wedge y{\prec}z\wedge \neg x{\prec} z) \ \vee\ (\neg x{\prec}y\wedge\neg y{\prec}x\wedge \neg x{=}y) \\
&& \vee \ x{\prec}\textit{min}\ \vee\ \textit{max}{\prec}x  \ \vee\ ( x{\prec} y\wedge y{\prec} s(x))\ \vee\ (\neg \textit{max}{=}x\wedge\neg x{\prec} s(x))  \\
&& \vee \  \neg P(\textit{min})\ \vee\ P(\textit{max})\ \vee\ (P(x)\wedge \neg P(s(x))).
\end{eqnarray*}
It has maximal determinacy $d(n)=n^2+2n+2$ for $n>1$. It is not weak and it is strong, indeed, its negation has a $2$-large model.
\end{example}

\begin{proof}  
For the first claim, consider the natural partial structure on $[n]$ that interprets $P$ by~$[n]$ and leaves $\textit{max}$ undefined. 
For the second claim consider the structure $\str A$ on $A:=\N\ \dot\cup\ \{\infty\}$ that interprets $\prec$ by the natural order extended by declaring $\infty$ larger than all natural numbers, $P$ by $\N$, $s$ by the natural successor extended by $s^{A}(\infty):=\infty$, and $\textit{min},\textit{max}$ by $0,\infty$.
Clearly, $\str A$ falsifies $\IND$. To see it has 2-bounded overflow, let $\str A_0$ be a partial substructure on a universe $A_0$ of size~$n$ and distinguish two cases. If $\infty\notin A_0$, then map $A_0$ order preserving onto $[n]$; the partial structure induced on $[n]$ has only $\textit{max}$ and $s(n-1)$ undefined, so is $2$-large. If $\infty\in A_0$ then map $A_0\setminus\{\infty\}$ onto $[n-1]$ as above and note that the induced partial substructure on $[n-1]\cup\{\infty\}$ is 1-large. 
\end{proof}

\begin{example}\label{ex:ba} Let $L:=\{\sqcup,\sqcap,\sim,f,0,1\}$ where $\sqcup,\sqcap$ are binary function symbols (in infix notation), $\sim,f$ are unary function symbols and $0,1$ are constants. Recall Boolean algebras are axiomatized by a finite set $E$ of equations in the language $L\setminus\{f\}$. The {\em Herbrandized atomicity principle} $\HAP$ negates the Skolemized infinity axiom stating ``here is an atomless Boolean algebra'': the existential closure of
$$\textstyle
\bigvee_{\zeta\in E} \neg \zeta\ \vee\ \neg f(0){=}0\ \vee\ \neg f(x){\sqcap} x{=}f(x)\ \vee\ (f(x){=}x\wedge \neg x{=}0).
$$
For $n$ a power of 2, its determinacy is $> s_L(n)-\log n$. It is neither weak nor strong.
\end{example}

\begin{proof} On a universe of cardinality  $n=2^k$, take a Boolean algebra with $k$ atoms and interpret~$f$ to map the interpretation of $0$ to itself, any other non-atom to an atom below it, and declare it undefined on all atoms. This shows the claim about the determinacy and that $\HAP$ is not weak.
To see $\HAP$ is not strong, let $\str A$ falsify $\HAP$. Any partial substructure of $\str A$ that contains~$n$ pairwise disjoint non-empty elements (in the sense of $\str A$), has $\sqcup^A$ completely undefined. This gives ${n\choose 2}$ many values of $\sqcup^A$ outside its universe.
\end{proof}

\begin{example}\label{ex:dlo} Take $L:=\{\prec,b,0,1\}$ for a binary relation symbol $\prec$, a binary function symbol $b$ (for ``between'') and constants $0,1$.
The {\em Herbrandized discreteness principle} $\HDP$ negates  the Skolemized infinity axiom stating ``here is a dense non-empty partial order'': the existential closure~of
\begin{eqnarray*}
&& x{\prec} x\ \vee\ (x{\prec} y\wedge y{\prec} z\wedge \neg x{\prec} z) \ 
 \vee\ (x{\prec} y\wedge \neg b(x,y){\prec} y)\ \vee\ (x{\prec} y\wedge \neg x{\prec}b(x,y))\ \vee\ \neg 0{\prec}1.
\end{eqnarray*}
The last disjunct ensures that the partial order is non-empty and thus $\HDP$ is valid in the finite.
$\HDP$ has determinacy $d(n)>2n^2 -2n$ for $n>1$. It is neither weak nor  strong.
\end{example}

\begin{proof} Note $s_L(n)=2n^2+2$. Consider a partial structure on  $[n]$ for $n>1$ that 
interprets~$\prec$ by the natural order, $0,1$ by themselves, and 
$b$ by some function that maps  $(i,j)$ with $|i-j|>1$ to some point between $i$ and $j$, and maps $(i,i)$ to 0, 
and is undefined on the $2(n-1)$ many pairs $(i,j)$ with $|i-j|=1$. This  does not verify $\HDP$, and has size $s_L(n)-2(n-1)$.

To see $\HDP$ is not strong, let $\str A$  falsify $\HDP$ and 
consider a linearly ordered subset $A_0$ of size $n$. Then $b^A$ takes a value outside~$A_0$ on each pair of $\prec^A$-consecutive points in $A_0$; this gives at least $n-1$ pairwise distinct values outside $A_0$. 
\end{proof}

\begin{remark}\label{rem:papa} Following \cite{beame}, every basic $L$-sentence $\varphi$  valid in the finite defines a complexity class, namely the type 1 NP search problems many-one reducible to $Q_\varphi$
(see Definition~\ref{df:assNPSP}). The classes associated to $\PHP,\textit{Onto-PHP},\textit{Left-PHP}$ and $\PAR$ 
are Papadimitriou's classes PPP, PPAD, PPADS and PPA~\cite{papa}. 
Papadimitriou~\cite{papa} showed that his classes contain many natural search problems of independent interest.\footnote{A minor difference is that usually the problems are defined only for structures with a universe of the form $[2^n]$ while we allow any $[n]$. The principle~$\PAR$, then, has to be slightly changed so as to be valid in even instead of odd structures.}  
\end{remark}

Finally, we mention an important example with built-in $\PV$:

\begin{example}\label{ex:iter} Let $\ITER$ be  $\exists y\ITER(y)$ (cf.~\cite{bk}) where $\ITER(y)$ is the following formula with a unary function symbol $f$ and built-in order $<$ and constant 0:
$$
f(0){=}0\ \vee \ f(y){<} y \ \vee \ \big(y{<}f(y)\wedge f(y){=}f(f(y))\big).
$$
It has maximal determinacy $d(n)=n$, so is not weak. 
\end{example}

\begin{proof} Interpret $f$ on $[n]$ by the successor, undefined on $n-1$.
\end{proof}

\begin{remark} \label{rem:iter}
The complexity class associated to $\ITER$ is the complexity class PLS  from~\cite{pls}.
Built-in symbols are necessary to characterize PLS. More precisely, assume that not all PLS problems are solvable in polynomial time. Then there does not exist a finitary combinatorial principle  without built-in symbols whose associated  class would equal PLS. 
\end{remark}

\begin{proof} Let $\varphi=\exists \bar y\psi(\bar y)$ be  such a principle, say in language $L$. If $\varphi$ fails in some infinite $L$-structure, then by Theorem~\ref{thm:morioka}, $Q_\varphi$ is not Turing reducible to $Q_\ITER$.
Otherwise $\varphi$ is valid. By Proposition~\ref{prop:Tred}~(3) it suffices to show that $\Sa$ proves $\exists y \q{\str B(L,x,\alpha)\models\psi(y)}$. But,
if $(N,\alpha^N)\models\Sa$ and $n\in N\setminus \{0\}$, then $\str B(L,n,\alpha^N)\models\varphi$ since $\varphi$ is valid, so $(N,\alpha^N)\models \exists y \q{\str B(L,n,\alpha^N)\models\psi(y)}$ by Lemma~\ref{lem:formAB}. 
\end{proof}


\subsection{Density arguments} \label{sec:dense}

We now establish the combinatorics needed for the forcing proofs of Theorems~\ref{thm:Briis} and~\ref{thm:main}. The sense of the forcing set-up in \cite{am} is to reduce independence questions for bounded arithmetics  to questions in finite combinatorics. Consequently, the combinatorics in this section are carried out in the standard model $\N$.
 For the rest of this section we let
\begin{enumerate}\itemsep=0pt
\item[--] $L$ and $\tilde L$ be finite languages;
\item[--] $r_{L}:=1+\max_{S\in L}\ar(S)$ and $r_{\tilde L}:=1+\max_{\tilde S\in \tilde L}\ar(\tilde S)$;
\item[--] $\tilde \varphi$ be a finitary combinatorial principle in the language $\tilde L$ as in Definition~\ref{df:fcp}, hence possibly with built-in $\PV$.
\end{enumerate}


 \begin{definition} \label{def:extend} Let $n\in\N\setminus\{0\}$ and
$p,q$ be partial $L$-oracles $p,q$ on $[n]$. The {\em size} $\|p\|$ of~$p$ is the size of $\str B(p)$ (as a partial structure, see Definition~\ref{df:size}).

We say $p$ {\em extends} $q$ if $\str B(p)$ extends~$\str B(q)$, in other words, if
$q_0\subseteq p_0$ and $q_1\subseteq p_1$; if additionally $b\in \N$ and $\|p\|\le \|q\|+b$ we call~$p$ a 
{\em $b$-extension} of $q$.

Call $a\in[n]$  {\em active in $\str B(p)$} if there are $S\in L$ and $\bar a\in[n]^{\ar(S)}$  such that $S^{[n]}(\bar a)\neq 1/2$ in~$\str B(p)$ 
and $a$ appears in $\bar a$ or $S$ is a function symbol and  $a=S^{[n]}(\bar a)$  in~$\str B(p)$.  
\end{definition}

\begin{lemma}\label{lem:dense1} Let $n\in\N\setminus\{0\}$, $\str B$ an $L$-structure, 
$p$ a partial $L$-oracle on $[n]$ such that $\str B(p)$ embeds into~$\str B$, $S\in L$ and $\bar a\in[n]^{\ar(S)}$. If
\begin{equation}\label{eq:psmall}
n>\|p\|\cdot r_L,
\end{equation}
then there is a $1$-extension $q$ of $p$ such that $\str B(q)$ embeds into $\str B$ and $S^{[n]}(\bar a)\not=1/2$ in $\str B(q)$.
\end{lemma}

\begin{proof}  
Write $p=\langle p_0,p_1\rangle$ and let $W\subseteq [n]$ be the set of $a\in[n]$ that are active in $\str B(p)$.
Note that $|W|\le \|p\|\cdot r_L$. 
Let $e$ be the embedding of $\str B(p)$ into $\str B$.

If $S$ is a relation symbol, obtain $q$ from $p$ by adding $\langle S,\bar a\rangle$ to $p_{b}$ where $b=S^{B}(e(\bar a))$ in $\str B$.

If $S$ is a function symbol and $v:=S^B(e(\bar a))$ is in the image of $e$, obtain 
$q$ from $p$ by adding $\langle S,\bar a,i\rangle$ to $p_{\bit(e^{-1}(v),i)}$ for all $i<|n|$.
If $v$ is not in the image of $e$, note that by \eqref{eq:psmall} there is $a\in[n]\setminus W$. Then change $e$ by mapping $a$ to $v$ 
and proceed as before.
\end{proof}

The following two lemmas show extendibility of a partial oracle to ensure that a partial $\tilde L$-structure of the form
$\str C((t_{\tilde s})_{\tilde s\in\tilde L},m,q)$ verifies $\tilde \varphi$. The first is simple and  useful for small $m$, and the second is useful for large~$m$ and the 
combinatorial core of the proof of Theorem~\ref{thm:main}.

\begin{lemma}\label{lem:dense3} Let $\str B$ be an $L$-structure,  $n,m,b_0\in\N\setminus\{0\}$,  $p$ a partial $L$-oracle on $[n]$ such that~$\str B(p)$ embeds into~$\str B$, and $(t_{\tilde S})_{\tilde S\in\tilde L}$ a family of decision trees of height at most $ b_0$.
If
\begin{equation}\label{eq:msmall}
n> r_L\cdot (\|p\|+b_0|\tilde L|m^{r_{\tilde L}-1}),
\end{equation}
then there exists a $b_0|\tilde L|m^{r_{\tilde L}-1}$-extension $q$ of $p$ 
such that  $\str C((t_{\tilde s})_{\tilde s\in\tilde L},m,q)$ verifies $\tilde \varphi$ and
$\str B(q)$ embeds into $\str B$.
\end{lemma}

\begin{proof} For $\tilde S\in \tilde L$ and $\bar a\in[m]^{\ar(\tilde S)}$ let $z_{\tilde S,\bar a}$ be a maximal sequence of 
$p$-answers to~$t_{\tilde S}$ on~$\bar a$. Note there are at most $ |\tilde L|m^{r_{\tilde L}-1}$ many such sequences and 
each has length at most $b_0$.
 If all these sequences are complete, then $\str C((t_{\tilde s})_{\tilde s\in\tilde L},m,q)$ is total and thus verifies $\tilde\varphi$ (being valid in the finite). Otherwise choose a 1-extension $q$ of $p$ that prolongues at least one of the answer sequences.
This is possible by the previous lemma if $r_L\|p\|< n$. By \eqref{eq:msmall} we can repeat this step until  all sequences are complete.
\end{proof}

By the {\em size} $|\varphi|$ of a formula $\varphi$, we mean the size (number of nodes) of the formula tree, that is, the number of occurrences of atomic subformulas and logical symbols $\wedge,\vee,\neg,\exists,\forall$.

\begin{lemma}[Core Lemma]\label{lem:dense2} 
Suppose the assumptions of the previous lemma hold and additionally
\begin{enumerate}\itemsep=0pt
\item[(i)] $\str B$ is $g$-large where $g:\N\setminus\{0\}\to\N$ is some function;
\item[(ii)] $n\ge (2b_0^2r_L+1)\cdot g(n)+r_L\cdot\|p\|$; 
\item[(iii)] $s_{\tilde L}(m)\ge 2b_0\tilde d(m)$ where $\tilde d$ is the determinacy of $\tilde \varphi$.
\end{enumerate}
Then there exists a $b_0|\tilde\varphi|$-extension $q$ of  $p$ 
such that $\str C((t_{\tilde s})_{\tilde s\in\tilde L},m,q)$ verifies $\tilde \varphi$
and $\str B(q)$ embeds into $\str B$. 
\end{lemma}

\begin{proof} We claim that it suffices to find $q$ as desired but neglecting the size bound, i.e., such that
 $q$ extends $p=\langle p_0,p_1\rangle$,  $\str C((t_{\tilde s})_{\tilde s\in\tilde L},m,q)$ verifies $\tilde \varphi$, and $\str B(q)$ embeds into $\str B$.

Given such $q=\langle q_0,q_1\rangle$, we have to find some $q'=\langle q'_0,q'_1\rangle$ with the same properties and of size $\|q'\|\le \|p\|+b_0|\tilde\varphi|$.
Recall  that $\tilde \varphi$ has the form \eqref{eq:basic} from Definition~\ref{df:basic}.
That $\str C:=\str C((t_{\tilde s})_{\tilde s\in\tilde L},m,q)$ verifies~$\tilde \varphi$ means that
there are a tuple~$\bar b$ from $[m]$
and  $i\in I$ such that~$\str C$ verifies~$\lambda_{ij}(\bar b)$ for all $j\in J$.
The literals $\lambda_{ij}(\bar b), j\in J,$ are verified in a partial substructure~$\str C'$ of $\str C$ of size at most 
$|J|<|\tilde\varphi|$.
For every $\tilde S\in \tilde L$ and $\bar b\in[m]^{\ar(\tilde S)}$ such that $\tilde S^{[m]}(\bar b)\neq 1/2$ in~$\str C'$ choose
a complete sequence $z_{\tilde S,\bar b}$ of $q$-answers to~$t_{\tilde S}$ on $\bar b$.
Consider the set  $Q$ of relevant (wrt $L,n$) queries $q$ needs to answer in these sequences. More precisely, this is  the set of
all relevant $\floor{f_{\tilde S}(\bar b,(z_{\tilde S,\bar b})_{<i})/2}$ where $i<|z_{\tilde S,\bar b}|-1$.
Then $|Q|<b_0|\tilde \varphi|$. 
By deleting certain elements from $q_0$ and $q_1$ we get a partial $L$-oracle $q'$ extending $p$ (and extended by $q$)  
of size at most $ \|p\|+b_0|\tilde\varphi|$ such that $\str C((t_{\tilde S})_{\tilde S\in\tilde L},m,q')$ extends~$\str C'$, so verifies $\tilde\varphi$. Namely, obtain $Q'$  from~$Q$ by adding $\langle S,\bar a,i\rangle$ whenever 
$\langle S,\bar a,i'\rangle\in Q$ for some $i'$ (here, $S\in L$, $\bar a\in[n]^{\ar(S)}$ and $i,i'$ range over $[|n|]$), and define $q'_0:=q_0\cap(p_0\cup Q')$, and similarly $q'_1$.
This proves the claim.

\medskip

For the sake of contradiction, assume that $q$ as in the claim does not exist.

\medskip

Consider a pair $(X,q)$ where $X$ is a set of pairs $(\tilde S,\bar b)$ with $\tilde S\in \tilde L$ and $\bar b\in[m]^{\ar(\tilde S)}$,
and~$q$ is a partial $L$-oracle on $[n]$ that extends $p$ and such that $\str B(q)$ embeds into $\str B$. From $(X,q)$ we compute another such pair $(X',q')$ as follows.

Choose an embedding $e$ of the partial structure $\str B(q)$ coded by $q$ into a $g(n)$-large partial 
substructure $\str B^*$ of $\str B$ with universe $B^*$ of size $n$. Note $\str B(q)$ extends $\str B(p)$ and $e$ embeds $\str B(p)$ into $\str B^*$. Let $\str B^*_n$ be the partial structure on $[n]$ which is isomorphic under $e$ to $\str B^*$ and let
$q^*$ be the partial oracle coding it.
Then $q^*$ extends~$q$. 

Choose $V\subseteq B$ witnessing that $\str B^*$ is $g(n)$-large. 
Let $W_n$ be the set of $a\in [n]$ active in~$\str B(p)$.
Let $R_n:=[n]\setminus W_n$ and let $R,W$ be the images of~$R_n,W_n$ under $e$.
Note 
$$
|R_n|\ge n-\|p\|\cdot r_L.
$$

For $(\tilde S,\bar b)\in X$ choose a maximal sequence $z_{\tilde S,\bar b}$ of $q^*$-answers to~$t_{\tilde S}$ on $\bar b$.
Let~$Y$ be obtained from $X$ by deleting all $(\tilde S,\bar b)$ such that $z_{\tilde S,\bar b}$ is complete. 
Then
$$
|Y|> |X|-\tilde d(m).
$$
Indeed, if at least $ \tilde d(m) $ many $z_{\tilde S,\bar b}$
are  complete, then $\str C((t_{\tilde S})_{\tilde S\in\tilde L},m,q^*)$ has
size at least $ \tilde d(m)$, and thus verifies $\tilde \varphi$. But this contradicts our assumption. 

Say $(\tilde S,\bar b)$ {\em touches} $a\in[n]$ if
there are  $j\le|z_{\tilde S,\bar b}|$ and $S\in L$ and $\bar a\in[n]^{\ar(S)}$ such that
$\floor{t_{\tilde S}(\bar b,z_{\tilde S,\bar b})_{<j}/2}$ equals $\langle S,\bar a\rangle$ or $\langle S,\bar a,i\rangle$ 
for some $i<|n|$, and such that $a$ appears in $\bar a$ or ($S$ is a function symbol and) $e(a)=S^{B^*}(e(\bar a))$ in $\str B^*$.

Note that any $(\tilde S,\bar b)\in Y$ touches at most $b_0\cdot r_L$ many $a\in[n]$. 
By averaging, there exists $r_0\in R_n$ which is touched by at most $ |Y|\cdot b_0\cdot r_L/|R_n|$ many pairs in $Y$.
Similarly, there exists $r_1\in R_1\setminus\{r_0\}$ touched by   at most $|Y|\cdot b_0\cdot r_L/(|R_n|-1)$ many pairs in $Y$.
Continuing like this we find pairwise distinct $r_0,\ldots, r_{|V|-1}$ such that at most
$$
|V|\cdot |Y|\cdot b_0\cdot r_L/(|R_n|-|V|)\le  b_0\cdot \frac{g(n)\cdot s_{\tilde L}(m)\cdot r_L}{n-\|p\|\cdot r_L-g(n)}
$$
many pairs in $Y$ touch any of them. Observe that (ii) implies that the denominators are positive.
Define $X'$ by deleting all these pairs from $Y$ and note
\begin{equation}\label{eq:shrink}
|X'|> |X|-\tilde d(m)- b_0\cdot \frac{g(n)\cdot s_{\tilde L}(m)\cdot r_L}{n-\|p\|\cdot r_L-g(n)}.
\end{equation}

Let $e'$ map $r_0,\ldots, r_{|V|-1}\in[n]$ bijectively onto $V$ and otherwise agree with $e$. Let $\str B'$ be the induced partial substructure of
 $\str B$ whose universe $B'$ is  the image of $e'$, 
 and let 
$q'$ be the partial $L$-oracle on $[n]$ such that $e':\str B(q')\cong\str B'$.
Then $q'$ extends $p$ since $e'$ equals $e$ on~$W_n$.


For $( \tilde S,\bar b)\in X'$ let $z'_{\tilde S,\bar b}$ be a maximal sequence of $q'$-answers 
to~$t_{\tilde S}$ on~$\bar b$. Then, as strings of bits, $z_{\tilde s,\bar b} $ is an initial segment of  $z'_{\tilde s,\bar b}$, i.e.\
$\bit(i,z'_{\tilde S,\bar b})=\bit(i,z_{\tilde S,\bar b})$ for all $i<|z_{\tilde S,\bar b}|-1$. We claim
$$
|z'_{\tilde S,\bar b}|>|z_{\tilde s,\bar b}|.
$$


Indeed, $t_{\tilde S}(\bar b,z_{\tilde S,\bar b})$ is odd and  
$\floor{t_{\tilde S}(\bar b,z_{\tilde S,\bar b})/2}$ equals $\langle S,\bar a,i\rangle$ for some function symbol $S\in L,\bar a\in[n]^{\ar(S)}$ and $i<|n|$ such that $S^{A}(e(\bar a))\in V$ in~$\str B$ (if $\floor{t_{\tilde S}(\bar b,z_{\tilde S,\bar b})/2}$ would not have this form, then 
$z_{\tilde S,\bar b}$ could be prolongued).
 Since all components of $\bar a$ are touched by 
$(\tilde S,\bar b)$ and $(\tilde S,\bar b)\in X'$, we know $e'(\bar a)=e(\bar a)$ and $S^{B'}(e'(\bar a))=e'(r_j)$ in $\str B'$ for some $j<|V|$.
As $z'_{\tilde S,\bar b}$ is maximal, it has a length $\ge |z_{\tilde S,\bar b}|+1$ 
with $\bit(|z_{\tilde S,\bar b}|-1,z'_{\tilde S,\bar b})=\bit(i,r_j)$.

Consider $(X_0,p_0)$ for $p_0:=p$ and $X_0$  the  set of all pairs $(\tilde S,\bar b)$ with 
$\tilde S\in \tilde L$ and $\bar b\in[m]^{\ar(\tilde S)}$.  Define a sequence $(X_0,p_0), (X_1,p_1),\ldots$ by iterating the function 
$(X,q)\mapsto (X',q')$. This gives a sequence $X_0\supseteq X_1\supseteq \cdots$ and a sequence of partial oracles $p=p_0,p_1,\ldots$ each extending~$p$. The maximal sequence of $p_i$-answers  to $t_{\tilde S}$ on $\bar b$ for pairs 
$(\tilde S,\bar b)\in X_i$ is prolonged in each step, and the pair gets deleted once the sequence is completed (recall the definition of $Y$ above). As the decision trees have height at most $b_0$, we conclude that $X_{b_0}$ is empty. On the other hand, the sets $X_i$ shrink per step as estimated in~\eqref{eq:shrink}. At the start $|X_0|=s_{\tilde L}(m)$, so
$$
0=|X_{b_0}|> s_{\tilde L}(m)-b_0\cdot \tilde d(m)- b_0^2\cdot\frac{g(n)\cdot s_{\tilde L}(m)\cdot r_L}{n-\|p\|\cdot r_L-g(n)},
$$
hence (recall $\tilde d(m)>0$)
$$
b_0> \frac{s_{\tilde L}(m)}{\tilde d(m)}\cdot \big( 1-b_0^2\cdot \frac{g(n)\cdot r_L}{n-\|p\|\cdot r_L-g(n)} \big).
$$
By (ii), the r.h.s. is $\ge (s_{\tilde L}(m)/\tilde d(m))\cdot 1/2$, a contradiction to (iii).
\end{proof}


\section{Typical forcing}

This section gives  a general method to construct models of $\To$ by forcing. We define {\em typical}  forcings with {\em typical graded} forcing frames that encompass many forcing type arguments in bounded arithmetic~\cite{pw,riisbrics,ajtai}. Theorem~\ref{thm:T12} states that such forcings  produce models of $\To$ if they satisfy a series of simple technical conditions. 
We give an application in the next section but believe general result is of independent interest.
The proof follows the set-up from \cite{am}, a simplified form of which is recalled in
Section~\ref{sec:basics}. Section~\ref{sec:definability} proves Theorem~\ref{thm:T12}. 
Throughout this section we fix
\begin{enumerate}\itemsep=0pt
\item[--] a countable language $L$ containing $\PV$;
\item[--]  a unary relation symbol $\alpha\notin L$;
\item[--]  an $L$-expansion $\N$ of the standard $\PV$-model;
\item[--] a countable proper elementary extension $M$ of $\N$. 
\end{enumerate}
\subsection{Forcing basics}\label{sec:basics}

We recall some standard forcing terminology. New notions are highlighted as definitions.

\medskip

A (countable) {\em  forcing frame} is a triple $(P,\fle,\mathcal D)$ where $(P,\fle)$ is a countable partial order 
with elements called {\em conditions} and $p\fle q$ reads as $p$ {\em extends} $q$, and $\mathcal D$ is a countable family of dense 
subsets of $P$. A subset of $P$ is  {\em dense (below $p$)} if every condition ($\fle p$) has an extension in it. Conditions
$p,q$ are {\em compatible}, written $p\|q$, if they  have a common extension.

\begin{definition} \label{def:graded} A {\em graded} forcing frame has additionally a non-increasing function $\| \cdot\|$ from~$P$ into $M$, that is,
$\|q\|\le^M \|p\|$ for all $p,q\in P$ with $p\fle q$. We say $p$ is a $b$-extension of $q$ if~$p\fle q$ and $M\models \|p\|\le\|q\|+b$.
A graded forcing frame is {\em typical} if $P\subseteq M$ and there are formulas $\q{x\fle y}$ and $\q{x\| y}$ and $\q{\|x\|=y}$ such that
for all $p,q\in P$ and $b\in M$:
\[
\begin{array}{llll}
M\models& \q{p\fle q}& \Longleftrightarrow & p\fle q;\\
M\models& \q{p\| q}&\Longleftrightarrow& p\| q;\\
M\models& \q{\|p\|=b}&\Longleftrightarrow& \|p\| =b.
\end{array}
\]
\end{definition}

Since this mode of speech does not depend on $\mathcal D$ we shall also refer to $(P,\fle,\|\cdot\|)$ as typical.
The {\em forcing language} is $L\cup\{\alpha\}$ together with the elements of $M$ as constants.
A {\em (universal) pre-forcing} is a binary  relation $\Vdash$ between conditions and sentences of the forcing language satisfying the following:
\begin{equation}\label{eq:recurrence}
\begin{split}
p\Vdash(\varphi\wedge\psi)&\ \Longleftrightarrow\  p\Vdash\varphi\text{ and }p\Vdash\psi;\\
p\Vdash\neg\varphi &\ \Longleftrightarrow\  q\not\Vdash\varphi \text{ for all } q\fle p;\\
p\Vdash\forall x\varphi(x)&\ \Longleftrightarrow\  p\Vdash\varphi(a)\text{ for all } a\in M.
\end{split}
\end{equation}
We write formulas with $\wedge,\neg,\forall$ and view $(\varphi\vee\psi)$ and $\exists x\varphi$ as abbreviations of  the classical dualities $\neg(\neg \varphi\wedge\neg\psi)$ and $\neg\forall x\neg\varphi$. Then
\begin{equation}\label{eq:recurrenceOR}
\begin{split}
p\Vdash(\varphi\vee\psi)&\ \Longleftrightarrow\ \{p\in P\mid p\Vdash\varphi\}\cup\{p\in P\mid p\Vdash\psi\}  \text{ is dense below }p;\\
p\Vdash\exists x\varphi(x)&\ \Longleftrightarrow\  \textstyle \bigcup_{a\in M} \{p\in P\mid p\Vdash\varphi(a)\} \text{ is dense below }p.
\end{split}
\end{equation}
 Also note that  $p\Vdash\neg\neg\varphi$ if and only if $\{p\in P\mid p\Vdash\varphi\}$ is dense below~$p$; for typcial forcings, defined next, this is equivalent to $p\Vdash\varphi$ (see Lemma~\ref{lem:forcing}~(e) below).

 \begin{definition}\label{def:typical}
A {\em typical forcing} is a pre-forcing that satisfies the following for $p,q\in P$ and all atomic sentences $\varphi$ and closed terms $s,t$  of the forcing language:
\begin{eqnarray*}
\text{(Extension)}&&\text{if $q\fle p\Vdash\varphi$, then $q\Vdash\varphi$};\\
\text{(Stability)}&&\text{if the set of conditions forcing $\varphi$ is dense below $p$, then $p\Vdash\varphi$};\\
\text{(Conservativity)}&&\text{if $\varphi$ does not mention $\alpha$, then: } p\Vdash\varphi\Longleftrightarrow M\models\varphi;\\
\text{(Extensionality)}&&\text{if $M\models s{=}t$, then: } 
p\Vdash\alpha(t)\Longleftrightarrow p\Vdash\alpha(s).
\end{eqnarray*}
\end{definition}

A {\em filter} $G$ is set of conditions that contains a common extension of any two  $p,q\in G$, and that contains 
any condition of which it contains an extension.
 A {\em generic} filter is one that intersects ``sufficiently many'' dense sets including those in $\mathcal D$. We refer to \cite[Definition~2.9]{am} for a definition, and just recall the standard lemma that every condition is contained in some generic filter (\cite[Lemma~2.12]{am}).
For such a filter $G$ \cite[Definition 2.16]{am} defines a structure~$M[G]$ interpreting the forcing language. 
We skip the definition as we only need the following genuine properties:

\begin{lemma}[Forcing Lemma]\label{lem:forcing} Assume $(P,\fle,\mathcal D)$ is a forcing frame and $\Vdash$ is a typical forcing. 
Then for every generic filter $G$, sentence $\varphi$ of the forcing language, and  $p\in P$:
\begin{enumerate}\itemsep=0pt
\item[(a)]  There is $\alpha^M_G\subseteq M$ such that $M[G]\cong (M,\alpha^M_G)$ as structures interpreting the forcing language ($(M,\alpha^M_G)$ interprets each constant $a\in M$ by $a$ itself).
\item[(b)] {\em (Truth Lemma)}  $M[G]\models \varphi$ if and only if $q\Vdash \varphi$ for some $q\in G$.
\item[(c)] {\em (Forcing Completeness)} $p\Vdash\varphi$ if and only if $M[H]\models\varphi$ for all generic filters $H$ containing $p$.
\item[(d)]  The set of sentences forced by $p$ is closed under logical consequence.
\item[(e)] (Extension), (Stability) and (Conservativity) hold for all sentences $\varphi$ of the forcing language.
\end{enumerate}
\end{lemma}

\begin{proof} This is proved in \cite{am}, we give precise references. First observe that, in the sense of \cite[Definition~2.16]{am}, $M[G]$ {\em is defined} for all generic filters $G$.
Thus,
(a)-(d) are \cite[Proposition~2.26]{am}, \cite[Theorem~2.19]{am}, \cite[Corollary~2.20~(2)]{am} and \cite[Corollary~2.20~(3)]{am}, respectively. 
In (e), (Extension) and (Stability) are \cite[Lemma~2.6~(1),(2)]{am}, and (Conservativity) is implied by (a) and (c).
\end{proof}


We remark that typical forcings behave nicely with bounded quantifiers, namely:
\begin{equation}\label{eq:bdqu}
p\Vdash\forall y{<}t\ \varphi(y)\ \Longleftrightarrow\ p\Vdash\varphi(a)\text{ for all $a \in M$ with }M\models a{<}t.
\end{equation}

\subsection{Partially definable  forcing}\label{sec:definability}

A condition $p$ is {\em compatible} with a  sentence $\varphi$ of the forcing language, written $p\|\varphi$,
if some extension of $p$ forces $\varphi$. Compatibility is dual to forcing in the sense that $p\|\varphi$ if and 
only if $p\not\Vdash\neg\varphi$, and, $p{\not\hspace*{-0.3ex}\|} \neg\varphi$ if and only if $p\Vdash\varphi$.

\begin{theorem}\label{thm:principal} Let $\Phi$ be a set of formulas of the forcing language. 
Under the assumptions of the previous lemma, suppose $\Vdash$ is {\em definable for $\Phi$}, i.e., for all $p\in P$
and $\varphi(\bar x)\in\Phi$, the set of tuples $\bar a$ from $M$ such that $p\|\varphi(\bar a)$ is definable in $M$. 
Then~$M[G]\models\LNP(\exists \Phi)$.
\end{theorem}

\begin{proof} This follows from  \cite[Theorem~3.5]{am} and \cite[Lemma~3.9~(1)]{am}.
\end{proof}

\begin{definition} Let $b_0\in N\subseteq M$. A {\em $\Delta^{b_0}_0(\alpha)$-formula with parameters from $N$} is a
$\PV\cup\{\alpha\}$-formula with parameters from $N$ all of whose quantifiers are {\em $b_0$-bounded}, i.e., of the form 
$\forall x{<}b_0$ and $\forall x{<}b_0$. Closing these formulas under
positive Boolean combinations, $b_0$-bounded quantifiers and bounded existential quantifiers $\exists x{<}t$ (where $t$ is a $\PV$-term without $x$ and possibly with parameters from $N$) yields the set of
{\em $\Sigma^{b_0}_1(\alpha)$-formulas with parameters from~$N$}. 
 \end{definition}

Recall (Section~\ref{sec:dense}) the size $|\varphi|$ of a formula $\varphi$ is the size of its formula tree. 

\begin{lemma}[Definability Lemma] \label{lem:definability} Let $(P,\fle,\mathcal D,\|\cdot\|)$ be a  typical graded forcing frame,~$\Vdash$ a typical forcing, and $b_0\in M\setminus\{0,1\}$.
Suppose
\begin{enumerate}\itemsep=0pt
\item[(a)] 
for every $r\in\N$ and $p\in P$ the set $\{ q\in P \mid q\text{ is a $b_0^r$-extension of } p\}$ is definable in $M$;
\item[(b)] for every  literal sentence $\varphi$ of the forcing language and  all $p^*,p\in P$ with $p\fge p^*\Vdash\varphi$ there 
exists a $b_0$-extension $q$ of $p$ that is compatible with $p^*$ and forces $\varphi$;
\item[(c)] for every atomic formula $\varphi(\bar x)$ of the forcing language and $p\in P$ the set of tuples $\bar a$ from $M$ such that
$p\Vdash\varphi(\bar a)$ is definable in $M$.
\end{enumerate}
Then $\Vdash$ is definable for $\Delta^{b_0}_0(\alpha)$-formulas with parameters from $M$.
\end{lemma}

Intuitively, conditions (a)-(c) are not much to ask for after a suitable choice for $b_0$, and this choice is mainly restricted by condition (a). Consider the usual case that $P$ has a minimum, is undefinable in $M$ and there is an upper bound $s\in M$
on $\|p\|,p\in P$. Then (a) implies $b_0^r\le^M s$ for all $r\in\N$, equivalently, $b_0$ is bounded by an infinitesimal power of $s$.

\begin{proof}[Proof of Lemma~\ref{lem:definability}]
By the Forcing Lemma~\ref{lem:forcing}~(d) we can restrict attention to $\Delta_0^{b_0}(\alpha)$-formulas in negation normal form (NNF), i.e., formulas built from 
literals by $\wedge,\vee$ and  $b_0$-bounded quantification $\exists x{<}b_0,\forall x{<}b_0$.
For  $\varphi$ in NNF let~$\varphi\neg$ be the formula in NNF obtained from $\neg\varphi$ by pushing the negation inside, that is, by swapping $\forall/\exists$ and $\wedge/\vee$ and literals with their complementary version. 
Let $k_\varphi$ denote the number of occurrences of  $\forall,\exists,\wedge,\vee$ in $\varphi$.

We show by induction on $k_\varphi$ that, if $\varphi$ has quantifier rank at most $r$, then:

\begin{enumerate}\itemsep=0pt
\item[(i)] for all tuples $\bar a$ from $M$ and all conditions
$p,p^*\in P$ with $p\fge p^*\Vdash \varphi(\bar a)$ there exists a $|\varphi|\cdot b_0^{r+1}$-extension $q$ of $p$ with $p^*\| q$ and  
$q\Vdash\varphi(\bar a)$;
\item[(ii)] there is a formula $\hat\varphi(z,\bar x)$ such that for all $p\in P$ the formula $\hat\varphi(p,\bar x)$ defines the set
 $\{\bar a\mid p\|\varphi(\bar a)\}$  in $M$;
 \item[(iii)] there is a formula $\tilde\varphi(z,\bar x)$ such that for all $p\in P$ the formula $\tilde\varphi(p,\bar x)$ defines the set
 $\{\bar a\mid p\Vdash\varphi(\bar a)\}$  in $M$.
\end{enumerate}

For $k_\varphi=0$, $\varphi$ is a literal. If $\varphi$ does not mention $\alpha$, then (i)-(iii) are trivial. 
If $\varphi(\bar x)$ is $\alpha(t(\bar x))$ for some term $t(\bar x)$, then (i) and (iii) hold by (b) and (c), respectively. For (ii), note that by~(b) we have that $p\|\alpha(t(\bar a))$ if and only if there is a 
$b_0$-extension $q$ of $p$ that forces~$\alpha(t(\bar a))$; this is easy to express using (a) and (c).

If $\varphi(\bar x)$ is $\neg\alpha(t(\bar x))$ for some term $t(\bar x)$, then (i) holds by (b).  For (ii), using (Stability), set
$\hat{\varphi}(z,\bar x):=\neg\widetilde{\alpha(t)}(z,\bar x)$. For (iii) set 
$\tilde{\varphi}(z,\bar x):=\neg\widehat{\alpha(t)}(z,\bar x)$.

\medskip

For the induction step we distinguish four cases whether $\varphi(\bar x)$ is obtained by $\wedge,\vee,\forall x{<}b_0$ or $\exists x{<}b_0$ from formulas $\psi$ with $k_\psi<k_\varphi$. 

\begin{enumerate}\itemsep=0pt

\item
Suppose $\varphi(\bar x)= (\varphi_0(\bar x)\wedge\varphi_1(\bar x))$. 
For (i) let $\bar a$ be a tuple from $M$ and suppose 
$$
p\fge p^*\Vdash(\varphi_0(\bar a)\wedge\varphi_1(\bar a)).
$$
Then $p\fge p^*\Vdash\varphi_0(\bar a)$.
By induction there is 
a $|\varphi_0|b_0^{r+1}$-exten\-sion~$q^0$ of~$p$ which is compatible with~$p^*$ and forces 
$\varphi_0(\bar a)$. 
Choose $q^*$ extending both $p^*$ and~$q^0$. Then $q^0\fge q^*\Vdash\varphi_1(\bar a)$. By induction there is 
a $|\varphi_1|b_0^{r+1}$-extension $q$ of $q^0$ which is compatible with~$q^*$ and forces~$\varphi_1(\bar a)$. 
Then $q$ is a $|\varphi_0|b_0^{r+1}+|\varphi_1|b_0^{r+1}< |\varphi|b_0^{r+1}$-extension of $p$ and compatible with~$p^*$. It forces  $\varphi_1(\bar a)$ by choice and  $\varphi_0(\bar a)$
as it extends~$q^0$, so $q\Vdash(\varphi_0(\bar a)\wedge\varphi_1(\bar a))$.

For (ii) observe we just showed that $p\|\varphi(\bar a)$ if and only if there is a $|\varphi|b_0^{r+1}$-extension~$q$ of~$p$ that forces both $\varphi_0(\bar a)$ and $\varphi_1(\bar a)$. This can be expressed using (a) and (iii) for~$\varphi_0,\varphi_1$.

For (iii) set $\tilde\varphi(z,\bar x):=\tilde\varphi_0(z,\bar x)\wedge\tilde\varphi_1(z,\bar x)$.

\item
Suppose $\varphi(\bar x)= (\varphi_0(\bar x)\vee\varphi_1(\bar x))$. For (i) let $\bar a$ be a tuple from $M$ and suppose 
$$
p\fge p^*\Vdash(\varphi_0(\bar a)\vee\varphi_1(\bar a)).
$$ 
Then there are $b\in\{0,1\}$ and 
$\tilde p$ such that $p^*\fge\tilde p\Vdash\varphi_b(\bar a)$  (recall~\eqref{eq:recurrenceOR}). Then $p\fge \tilde p$ and induction gives  
a $|\varphi_b|b_0^{r+1}$-extension $q$ of $p$ which is compatible with~$\tilde p$, and hence also 
with~$p^*$, and forces $\varphi_b(\bar a)$,
and hence also $\varphi(\bar a)$. 

For (ii) set $\hat\varphi(z,\bar x):=\hat\varphi_0(z,\bar x)\vee\hat\varphi_1(z,\bar x)$.

For (iii) set $\tilde\varphi(z,\bar x):=\neg\widehat{\varphi\neg}(z,\bar x)$; note $\varphi\neg$ is a conjunction with $k_{\varphi\neg}=k_{\varphi}$, so $\widehat{\varphi\neg}$ has been defined in the previous case.

\item
Suppose $\varphi(\bar x)=\forall y{<} b_0 \psi(y,\bar x)$.
Let $\bar a$ be a tuple from $M$ and suppose 
$$
p\fge p^*\Vdash\forall y{<} b_0\ \psi(y,\bar a).
$$ 
We claim that for every~$b\le^M b_0$  there
 is a $b\cdot |\psi|\cdot b_0^{r}$-extension~$q^b$ of $p$ such that $q^b\|p^* $ and $q^b\Vdash \psi(c,\bar a)$ for all~$c<^Mb$. 
 
 This is an $M$-definable property of $b$. Indeed, using (a) and the definability of $\|\cdot\|$ (Definition~\ref{def:graded}), the set of $b\cdot |\psi|\cdot b_0^{r}$-extensions of $p$ is definable in $M$ (with parameter~$b$), forcing $\psi(c,\bar a)$ for all~$c<^Mb$ is expressed using $\hat\psi(z,y,\bar x)$, and compatibility with $p^*$ is expressed using $\q{x\| y}$.

Since $M$ is an elementary extension of $\N$, it satisfies  induction  for all formulas in its language. 
We  can thus prove our claim by induction on $b$ in $M$. Then (i) will follow, witnessed by $q^{b_0}$ (recall~\eqref{eq:bdqu}).

For~$b=0$ take~$q^0:=p$. Assume that~$b<^Mb_0$ and we found $q^b$ as desired. Let $q^*$ be a common extension of $q^b$
and $p^*$. Then $q^b\fge q^*\Vdash \psi(b,\bar a)$. Note $\psi$ has quantifier rank at most $r-1$. Applying (i) for $\psi$ gives a $|\psi|b_0^{r}$-extension 
$q^{b+1}$ of~$q^b$ that forces $\psi(b,\bar a)$ and is compatible with $q^*$ and hence 
with~$p^*$; since~$q^{b+1}$ extends $q^b$ it forces $\psi(c,\bar a)$ 
for all $c$ with $M\models c{<}b{+}1$. 

To see (ii), note we showed that $p\|\varphi(\bar a)$ if and only if there exists a  $|\psi|\cdot b_0^{r+1}$-extension of~$p$ forcing $\psi(c,\bar a)$ for all $c$ with   $M\models c{<}b_0$. This is easily expressed using (a) and the formula~$\tilde\psi(z,y,\bar x)$.

For (iii) set $\tilde\varphi(z,\bar x):=\forall y{<}b_0\tilde{\psi}(z,y,\bar x)$ (recall \eqref{eq:bdqu}).

\item
Suppose $\varphi(\bar x)= \exists y{<} b_0\psi(y,\bar x)$.  For~(i) let
$\bar a$ be a tuple from $M$ and suppose $$p\fge p^*\Vdash\exists y{<} b_0\ \psi(y,\bar a).
$$ 
Then there are $b\in M$ and 
$\tilde p$ such that $p^*\fge\tilde p\Vdash (b{<}b_0\wedge \psi(b,\bar a))$  (recall~\eqref{eq:recurrenceOR}). By (Conservativity), $b<^M b_0$ and
$\tilde p\Vdash\psi(b,\bar a)$.
As $\psi(y,\bar x)$ has quantifier rank $\le r-1$, induction gives 
a $|\psi|b_0^{r}$-extension $q$ of $p$ which is compatible with~$\tilde p$, and hence with $p^*$, and forces~$\psi(b,\bar a)$
and hence $\exists y{<} b_0\ \psi(y,\bar a)$.

For (ii), note we just saw that $p\| \exists y{<} b_0\psi(y,\bar a)$ if and only if $p\| \psi(b,\bar a)$ for some $b<^Mb_0$. We thus set
$\hat\varphi(z,\bar x):=\exists y{<}b_0\hat\psi(z,y,\bar x)$.

For (iii), set $\tilde\varphi(z,\bar x):=\neg\widehat{\varphi\neg}(z,\bar x)$; note $\varphi\neg$ starts with $\forall y{<}b_0$ and has $k_{\varphi\neg}=k_{\varphi}$, so $\widehat{\varphi\neg}$ has been defined in the previous case.
\end{enumerate}
This finishes the proof of the Definability Lemma.
\end{proof}

We are ready to prove the main result in this section, a  general method to produce models of $\To$ by typical forcings. 

\begin{definition}
A  {\em $\PV$-cut} in $M$ is a substructure $N$ of the $\PV$-reduct of~$M$ such that $a <^Mb\in N$ implies $a\in N$ for all $a,b\in M$. 
\end{definition}

Recall the notation $\alpha^M_G$ from the Forcing Lemma~\ref{lem:forcing}~(a).

\begin{theorem}\label{thm:T12}  Assume the forcing frame $(P,\fle,\mathcal D,\|\cdot\|),$ the forcing $\Vdash$ 
 and $b_0\in M$ satisfy the assumption of the previous lemma, and
let $G$ be a generic filter.

 Assume further that $N$
is a $\PV$-cut of $M$  such that $b_0\in N$ and $b_0$ {\em bounds lengths in~$N$}, i.e., $N\models \forall x\ |x|{<}b_0$. 
Set 
$$
\alpha^N:=\alpha^M_G\cap N.
$$
Then $(N,\alpha^N)$ has a unique expansion to a model of $\To$.
 \end{theorem}

\begin{proof} By Lemma~\ref{lem:definability} and Theorem~\ref{thm:principal} we have 
$(M,\alpha^M_G)\models\LNP(\exists\Delta_0^{b_0}(\alpha))$.
We claim that 
\begin{equation*}\label{eq:s12}
(N,\alpha^N)\models\LNP(\Sigma_1^{b}(\alpha)).
\end{equation*}
For contradiction, assume $\varphi(x)$ is a $\Sigma_1^{b}(\alpha)$-formula with parameters from~$N$ that defines in~$(N,\alpha^N)$ a non-empty set without  minimum. Since~$b_0$ bounds lengths in $N$,
$\varphi(x)$ is in $(N,\alpha^N)$ equivalent to a $\Sigma_1^{b_0}(\alpha)$-formula~$\varphi'(x)$ with parameters from $N$.
Since $N$ is a $\PV$-cut in $M$, $\varphi'(x)$ defines also in~$(M,\alpha^M_G)$ a non-empty set without minimum. But a standard collection argument  
(see e.g.~\cite[Proof of Theorem~4.3]{am}) shows $\varphi'(x)$ is in $(M,\alpha^M_G)$ equivalent to a $\exists\Delta_0^{b_0}(\alpha)$-formula. We thus get a contradiction to $\LNP(\exists\Delta_0^{b_0}(\alpha))$. 

Clearly, $N\models\forall\PV$, so the theorem follows  by Lemma~\ref{lem:PVpa} and 
Proposition~\ref{prop:cons}.
\end{proof}


\section{Riis' theorem and extensions}\label{sec:riis}

We define a forcing whose conditions are partial oracles on $[n]$ coding partial structures that do not verify a given $\varphi$. The oracle in the generic expansion then codes a total structure on~$[n]$ where $\varphi$ fails. It is a routine task to verify that our forcing has various desirable properties (typical, graded,  etc.). We shall give the details in Section~\ref{sec:fwps}. Sections~\ref{sec:Briis} and \ref{sec:riisext} then prove certain stronger variants of Theorems~\ref{thm:Briis} and \ref{thm:main} as an application of Theorem~\ref{thm:T12}. We view Theorem~\ref{thm:main} as an extension of Theorems~\ref{thm:Briis} because the proof of the latter is not much more than the former plus an additional application of the Core Lemma~\ref{lem:dense2}. 
The proof exemplifies the role of forcing in bounded arithmetic, as viewed in~\cite{am}, to reduce independence to finite combinatorics, here, density arguments.

\subsection{Forcing with partial structures}\label{sec:fwps}

We define a notion of forcing in the following situation:
\begin{enumerate}\itemsep=0pt
\item[--] $L$ is a finite language and $\varphi$ is a basic $L$-sentence (Definition~\ref{df:basic}); 
\item[--] $\N$  is an expansion of the standard $\PV$-model interpreting a countable  language including ($\PV$ and) $L$; 
\item[--] $\str B\not\models\varphi$ where $\str B$ is the $L$-reduct of $\N$;
\item[--] $M$ is a countable proper elementary extension of $\N$;
\item[--] $b_0,n\in M\setminus \N$ such that $b_0^k<^M n$ for all $k\in\N$.
\end{enumerate}
Hence the role of $L$ in Section~\ref{sec:basics} is played by the language of $M$ here. 
\medskip

There is a $(\PV\cup L)$-formula $\PaOr(y,x)$ that defines in $\N$ the pairs $(m,p)$ such that $m>0$ and
$p$ is a code of a partial $L$-oracle on $[m]$ and  the partial $L$-structure $\str B(L,m,p)=\str B(p)$  is embeddable into $\str B$. The size $\|p\|$ of such~$p$ does not depend on $m$ and is definable in $\N$. We  have a $(\PV\cup L)$-formula ``$x$~is relevant (wrt~$L,y$)'' defining in $\N$ the set of pairs $(a,m)$ such that $a$ is relevant (wrt $L,m$).
A {\em partial $L$-oracle on $[n]$ in $M$} is an element satisfying $\PaOr(n,x)$ in $M$, and $a\in M$ is {\em relevant (wrt $L,n$)} if it satisfies ``$x$ is relevant (wrt $L,n$)'' in $M$. We do not distinguish a  partial $L$-oracle $p$ notationally from the pair of sets it codes. We write
$p=\langle p_0,p_1\rangle$ (in $M$); formally, 
$p_0,p_1$ are $(p)_0,(p)_1$ calculated in $M$. 
Since $M$ is an elementary extension of $\N$, the function $\|\cdot\|$ extends to $M$.

Let $P\subseteq M$ be the set of partial $L$-oracles $p$ (on $[n]$)
 in $M$ such that  
 $$
 M\models \|p\|\le b_0^k\text{ for some }k\in\N.$$
 We let $p,q,\ldots$ range over~$P$. Note that $P$ is not definable in $M$. We set $p\fle q$ if and only if~$p$ extends $q$ in the sense of Definition~\ref{def:extend} (applied in $M$). 

\begin{lemma} \label{lem:typ}
$(P,\fle,\|\cdot\|)$ is a typical graded forcing frame.
\end{lemma}

\begin{proof} 
Clearly,  $\|\cdot\|$ is non-increasing. For typicality, we already noted the formula $\q{\|x\|{=}y}$ and set (recall $x\in y$ is $\bit(y,x){=}1$)
\begin{eqnarray*}
\q{x{\fle} y}&:=&\forall z(z\in (y)_0\to z\in (x)_0)\wedge \forall z(z\in (y)_1\to z\in (x)_1);\\
\q{x\| y}&:=&\exists z(\PaOr(n,z)\wedge \q{z{\fle} x}\wedge \q{z{\fle} y}).
\end{eqnarray*}
A pair of conditions $(p,q)$ satisfies $\q{x\| y}$ in $M$ if and only if  $p$ and $q$ have a common extension {\em in $P$}. Indeed, if there is a partial $L$-oracle extending both $p$ and $q$, then there is one of size at most $ \|p\|+\|q\|$ which hence is in $P$.
\end{proof}

This completes the definition of the forcing frame up to the choice of $\mathcal D$. This choice will be based on the following corollaries to Section~\ref{sec:dense}, explaining the title of that section.

\begin{corollary}\label{cor:D1}
For every relevant $a\in M$, the  set $
D(a):=\{q\in P\mid a\in q_0\cup q_1\}
$ is dense. 
\end{corollary}
\begin{proof} Let $p\in P$, and  $a\in M$ be relevant, i.e., in $M$ of
the form $\langle S,\bar a,i\rangle$ or $\langle S,\bar a\rangle$ for $S\in L$ and~$\bar a\in [n]^{\ar(S)}$ and $i<|n|$ in $M$. 
Lemma~\ref{lem:dense1} formalizes as a sentence which is true in $\N$, and hence in $M$. The assumption~\eqref{eq:psmall} of this lemma holds in $M$ for all $p\in P$. Hence its conclusion gives in $ M$ a 1-extension $q$ of $p$ in $D(a)$. Clearly, $q\in P$. 
\end{proof}

The  following corollary is proved by a case distinction as to whether $m$ is small or large and then applies Lemma~\ref{lem:dense3} or \ref{lem:dense2}. It is not needed in the proof of Riis' theorem given in the next section.

\begin{corollary}\label{cor:D2} Assume $\str B$ is $n^{o(1)}$-large and $\tilde \varphi$ is a  weak finitary combinatorial principle in the  language $\tilde L$.
Then  for every $m\in M\setminus\{0\}$ and every family of decision trees  $(t_{\tilde S})_{\tilde S\in\tilde L}$
in~$M$ of height at most $b_0$  the following set is dense:
$$
D((t_{\tilde S})_{\tilde S\in\tilde L},m):=\big\{q\in P\mid \str C((t_{\tilde S})_{\tilde S\in\tilde L},m,q)\text{ verifies }\tilde \varphi\big\}.
$$
\end{corollary}

\begin{proof}  
There is a definable function in $\N$ that maps $n,p$ to the (natural numbers coding the)
 partial structure $\str B(L,n,p)=\str B(p)$. Similarly, 
  $\str C((t_{\tilde S})_{\tilde S\in\tilde L},m,q)$ is the value of a definable (in~$\N$) function 
 on $m,q$ and the parameters in the definitions of the decision trees $t_{\tilde S}, \tilde S\in\tilde L$. 
 The size function $s_{\tilde L}$, the determinacy $\tilde d$ of $\tilde\varphi$ are clearly definable in $\N$, and so is some function $g(n)\le n^{o(1)}$ witnessing that $\str B$ is $g$-large. Since $M$ is 
 an elementary extension of~$\N$ these functions extend to $M$, and we denote the extensions again by by $s_{\tilde L},\tilde d$ and $g$. We have
Lemmas~\ref{lem:dense3} and \ref{lem:dense2} for $M$ instead~$\N$. 
Let $p\in P$ be given. We distinguish two cases.

Assume first that $m$ satisfies  (iii) of Lemma~\ref{lem:dense2}. We have assumptions~(i) and (ii) of this lemma. For (ii), observe that
overspill gives  $t\in M\setminus \N$ such that 
$M\models g(n) {<} n^{1/t}$; hence
the r.h.s.\ of~(ii) is at most $n^{1/t'}$ for some~$t'\in M\setminus \N$. 
The conclusion of Lemma~\ref{lem:dense2} gives in~$M$ a $b_0|\tilde\varphi|$-extension 
$q$ of $p$ in $D((t_{\tilde S})_{\tilde S\in\tilde L},m)$. Note $q\in P$ because $\|q\|\le \|p\|+|\tilde\varphi|b_0\le b_0^k$ for suitable standard $k\in \N$.

Now assume that $m$ violates (iii) of Lemma~\ref{lem:dense2}, i.e., $s_{\tilde L}(m)< 2b_0\tilde d(m)$ in $M$. 
As $\tilde \varphi$ is weak, $s_{\tilde L}(m)\ge m^{1/\ell}\cdot \tilde d(m)$ for some $\ell\in\N\setminus\{0\}$. As $\tilde\varphi$ is valid in the finite, $\tilde d(m)>0$ in~$M$.  It follows that $m< (2b_0)^\ell$ in~$M$. 
 But then the assumption \eqref{eq:msmall} of Lemma~\ref{lem:dense3} holds true in $M$: the r.h.s.\ is bounded by $b_0^{k}$ for some standard $k\in \N$ and $b_0^k<n$ in $M$.
The conclusion of this lemma gives in~$M$ some $q\in D((t_{\tilde S})_{\tilde S\in\tilde L},m)$ extending~$p$; indeed $q\in P$ because $\|q\|\le \|p\|+b_0|\tilde L|m^{r_{\tilde L}-1}< b_0^{k}$ in $M$.
\end{proof}

We next define a typical forcing $p\Vdash\varphi$ for $p\in P$ and $\varphi$  a sentence in the forcing language.
One might be tempted  to define $p\Vdash\alpha(t)$ if and only if $t^M\in p_1$; recall $t^M$ is the value of the closed term $t$ of the forcing language in $M$ (treating its constants from $M$ as parameters).
This, however, does not work: assume $t^M=\langle S,\bar a\rangle\notin p_1$ with $S^{[n]}(\bar a)=1/2$ in~$\str B(p)$; 
it might be that 
every partial substructure of $\str B$ containing an isomorphic copy of $\str B(p)$ is such that the copy of $\bar a$ is mapped to 1 by 
$S^A$ in $\str B$. In this case, $t^M\in q_1$ for all extensions $q$ of $p$ with $S^{[n]}(\bar a)\neq 1/2$ in~$\str B(q)$. Then Forcing Completeness (Lemma~\ref{lem:forcing}~(c)) 
fails:~$\alpha(t)$ is not forced by $p$ but holds in all generic expansions built by filters containing $p$. 

The issue is sidestepped using a weaker and slightly more technical definition:

\begin{lemma}\label{lem:defass} There is exactly one typical forcing $\Vdash$ satisfying for all closed terms $t$ of the forcing language and all $p\in P$:
 \begin{equation}\label{eq:forcedef}
 p\Vdash\alpha(t)\ \Longleftrightarrow
\  t^{M}\textup{ is relevant and $t^{M}\notin q_0$ for every 1-extension $q$ of $p$}. 
 \end{equation}
Moreover,
 $(P,\fle,\mathcal D,\|\cdot\|),\Vdash$ 
 and $b_0$ satisfy the assumption of the Definability Lemma~\ref{lem:definability}.
\end{lemma}

\begin{proof} We define $p\Vdash\varphi$ for atomic formulas $\varphi$ without $\alpha$ according
to (Conservativity) and use the recurrence~\eqref{eq:recurrence} to define it on more 
complex formulas. Uniqueness being clear, we check this defines a typical forcing. The rest being obvious we have to check (Extension) and (Stability) for atoms of the form $\alpha(t)$ where $t$ is a closed term of the forcing language.

For (Extension) assume $p\fge q\not\Vdash\alpha(t)$. We  show $p\not\Vdash \alpha(t)$. This is clear if~$t^{M}$ is not relevant. Otherwise there is a 1-extension $q'$ of $q$ with
$t^{M}\in q'_0$. Deleting some elements from $q_0',q'_1$ gives a 1-extension $p'$ of $p$ with $t^{M}\in p'_0$, so $p\not\Vdash\alpha(t)$.

For (Stability) assume $p\not\Vdash\alpha(t)$. We have to find
some extension $q$ of $p$ that does not have an extension forcing $\alpha(t)$. If 
$t^{M}$ is not relevant, we take $q:=p$. Otherwise there is a 1-extension $q$ of $p$ with $t^{M}\in q_0$. Clearly, no extension of $q$ 
forces $\alpha(t)$.

We now verify the assumptions of the Definability Lemma~\ref{lem:definability}.
Assumptions (a) and (c) being clear, we prove (b). 
Let $\varphi$ be a literal sentence of the forcing language and suppose 
 $p\fge p^*\Vdash \varphi$. 
We can assume $\varphi$ mentions $\alpha$
(otherwise  take $q:=p$), so equals $\alpha(t)$ or $\neg\alpha(t)$ for some closed term $t$. Assume the former (the latter case is similar).
Then $t^{M}$ is relevant, so Corollary~\ref{cor:D1} gives
 $r\fle p^*$ with $r\in D(t^M)$. Then $t^{M}\in r_1$ because $r\Vdash\alpha(t)$. From $r$ get a 1-extension $q$ of~$p$
with $t^{M}\in q_1$  by deleting some elements from $r_0,r_1$. 
Clearly,~$q$ is compatible with $p^*$ and forces $ \alpha(t)$. 
\end{proof}

Finally, 
we observe that the generic $\alpha_G^M$ from the Forcing Lemma~\ref{lem:forcing}~(a) is as expected:

\begin{lemma}\label{lem:alphagood}
For every relevant $a\in M$:
\begin{eqnarray*}\label{eq:p1}
\textstyle
a\in \alpha^{M}_G& \Longleftrightarrow & a\in p_1\text{ for some } p\in G\\\label{eq:p0}
& \Longleftrightarrow & a\not\in p_0\text{ for all } p\in G.
\end{eqnarray*}
\end{lemma}
\begin{proof}
If
$a\in \alpha^{M}_G$, then there is
$p\in G$ forcing~$\alpha(a)$ by  Lemma~\ref{lem:forcing}~(a),(b). By genericity there is $q\in G\cap D(a)$. Then $a\notin q_0$ as otherwise $q\Vdash\neg \alpha(a)$ and  then~$p,q\in G$ would not be compatible. Hence $a\in q_1$. 
Conversely, if $a\in p_1$ for some $p\in G$, then $p\Vdash\alpha(a)$. Then  $a\in \alpha^{M}_G$ by Lemma~\ref{lem:forcing}~(a),(b).
This shows the first equivalence. The second is similar.
\end{proof}

\subsection{Proof of Theorem~\ref{thm:Briis}}\label{sec:Briis}

We prove the following stronger version of Theorem~\ref{thm:Briis}. 
A function $f:\N\to\N$ is  {\em subexponential} if $f(n)\le 2^{n^{o(1)}}$. 
If $f$ is definable (in the standard $\PV$-structure $\N$), then it
has an extension $f^M$  to any elementary extension $M$ of $\N$. 
Call a $\PV$-cut $N$ of $M$ {\em subexponential in~$n$} if $f^M(n)\in N$ for  all definable subexponential  functions $f:\N\to\N$.

To be clear about the notation in the following statement, recall that by Proposition~\ref{prop:cons} every model of $\To$ has the form $\langle N,\alpha^N\rangle$ for $N\models\forall\PV$ and $\alpha^N\subseteq N$.

\begin{theorem}\label{thm:strongmain} Let $L$ be a finite language and $\varphi$ a basic $L$-sentence without built-in symbols
that fails in some infinite model. 

Then there exists a model $\langle N,\alpha^N\rangle$ of $\To$ and $n\in N\setminus\{0\}$ such that 
$$
\str B(L,n,\alpha^N)\not\models\varphi.
$$

Moreover, if $\psi(x)$ is a $\PV$-formula that defines an unbounded set in~$\N$, then $N,n$ can be chosen such that $N$ is a $\PV$-cut in an elementary extension $M$ of $\N$ such that $N$ is subexponential in $n$ and   $M\models\psi(n)$. 
\end{theorem}

\begin{remark} If $\varphi$ is not valid in the finite, the  first statement is trivial
but the second is not. An interesting case is that the spectrum of $\neg\varphi$ is co-infinite and belongs to  the polynomial hierarchy, or equivalently, the set of $n>0$ such that $\varphi$ is valid in structures of size $[n]$ is infinite and definable by a bounded $\PV$-formula $ \q{\varphi \textit{ is valid on } [x]}$. Then we get $\str B(L,n,\alpha^N)\not\models\varphi$ and $N\models\q{\varphi \textit{ is valid on } [n]}$ (since this is bounded and true in $M$).
\end{remark}

\begin{proof}[Proof of Theorem~\ref{thm:strongmain}] The proof consists mainly in putting the pieces together. 
Let $\str B$ be an infinite model of~$\neg\varphi$. We can assume it has universe $B=\N$. We let $\N$ be the $\PV\cup L$-structure whose $\PV$-reduct is the standard model and whose $L$-reduct is $\str B$. Let $f_0,f_1,\ldots$ enumerate the definable subexponential functions. For every $k\in\N$ the formula
\begin{equation}\label{eq:nb}
\textstyle
\psi(x)\wedge y^k{<}x\wedge\bigwedge_{i<k}|f_i(x)|{<}y
\end{equation}
is satisfiable in $\N$. Thus there exists a countable elementary extension $M$ of $\N$ and $n,b_0\in M$ such that assigning $n$ to $x$ and $b_0$ to $y$ satisfies  \eqref{eq:nb} for all $k\in\N$. Clearly, $n,b_0\in M\setminus\N$.

Let $\mathcal D$ be the family of dense sets $D(a), a\in M,$ from Corollary~\ref{cor:D1}. The previous section gives a typical graded forcing frame $(P,\fle,\mathcal D,\|\cdot\|)$ and a typical forcing $\Vdash$ satisfying the assumptions of the Definability Lemma~\ref{lem:definability} (see Lemma~\ref{lem:defass}).
Let $G$ be a generic filter and 
$$
\textstyle N:=\bigcup_{k\in\N}\big\{a\in M\mid a\le^M f^M_k(n)  \big\}.
$$
This is a $\PV$-cut in $M$ and $b_0$ bounds lengths in $N$. By Theorem~\ref{thm:T12}, $(N,\alpha^N)$ has an expansion $\langle N,\alpha^N\rangle\models\To$ where $\alpha^N:=\alpha^M_G\cap N$. Note that $\alpha^N=\alpha^M_G$ since $\alpha^M_G$ contains only relevant elements and these are in $N$. 

The ``moreover'' part is obvious. To  verify $\str B(L,n,\alpha^N)\not\models\varphi$, we first observe that $\str B(L,n,\alpha^N)$ is the union of the partial structures $\str B(p), p\in G$. More precisely and first, every $\str B(p), p\in G,$ is a partial substructure of $\str B(L,n,\alpha^N)$ because, by
 Lemma~\ref{lem:alphagood}, sequences of $p$-answers are sequences of $\alpha^N=\alpha^M_G$-answers.
Second, assume
$S^{[n]}(\bar a)=b$ in $\str B(L,n,\alpha^N)$ for some~$S\in L$ and $\bar a,b$ from $[n]$. We claim that $S^{[n]}(\bar a)=b$ in $\str B(p)$ for some $p\in G$.
Say, $S$ is a function symbol (the case of a relation symbol is similar). Choose 
$p\in G\cap D(\langle S, \bar a,0\rangle)$. Then $\langle S, \bar a,i\rangle\in p_0\cup p_1$ for all $i<|n|$, so  $S^{[n]}(\bar a)\neq 1/2$ in~$\str B(p)$. Since $\str B(p)$ is a partial substructure of $\str B(L,n,\alpha^N)$ we have
$S^{[n]}(\bar a)=b$ in $\str B(p)$.

We now verify $\str B(L,n,\alpha^N)\not\models\varphi$. Assume otherwise and recall $\varphi$ has the form~\eqref{eq:basic} (Definition~\ref{df:basic}). 
 Choose~$i\in I$ and a  tuple $\bar a$ 
from $[n]$ such that $\str B(L,n,\alpha^N)$  
 verifies  $\lambda_{ij}(\bar a)$ for all $j\in J$. The literals $\lambda_{ij},j\in J,$ are verified in a partial substructure $\str C$ of $\str B(L,n,\alpha^N)$ of size at most $|J|$. Let 
$
(S_0,\bar a_0),\ldots, (S_{|J|-1},\bar a_{|J|-1})
$
list all pairs $(S,\bar a)$ with $S\in L,\bar a\in[n]^{\ar(S)}$ and $S^{[n]}(\bar a)\neq 1/2$ in~$\str C$. 
As observed above,
for every $j<|J|$ there is $p^j\in G$ such that the value $S_j^{[n]}(\bar a_j)$ in~$\str B(p^j)$ is equal to this value in~$\str C$.
Since $G$ is a filter, there is $p\in G$ extending all~$p^j,j<|J|$. 
Then~$\str C$ is a partial substructure of $\str B(p)$, so $\str B(p)$ verifies $\varphi$. 
As $p\in P$ we have that $\str B(p)$ embeds into (the $L$-reduct of) $M$. Hence
$M\models\varphi$ by  Lemma~\ref{lem:pres},  so $\str B\models\varphi$ by elementarity -- a contradiction. 
\end{proof}

The following is repeated from the Introduction and  strengthens of Buresh-Oppenheim and Morioka's Theorem~\ref{thm:morioka}. Recall Definitions~\ref{df:conseq} and~\ref{df:assNPSP} and Example~\ref{ex:iter}.

\begin{corollary}  If $\varphi$ is  finitary combinatorial principle without built-in symbols that fails in some  infinite model, then  $Q_\varphi$ is independent  from $Q_\ITER$ over $\To$ .
\end{corollary}

\begin{proof} Assume  $\varphi$ satisfies the hypothesis, write it  as $\exists\bar y\psi(\bar y)$ for $\psi(\bar y)$ quantifier free and say it has language $L$. 
By the previous theorem and Lemma~\ref{lem:formAB}, $\To$ does not prove
$\exists y\q{\str B(L,x,\alpha)\models\psi(y)}$. But it is not hard to see that $\To$ proves $\exists y\ \q{\str B(L,x,\alpha)\models\ITER(y)}$
(cf.~\cite[Theorem 4.4]{bk}). Now apply Corollary~\ref{cor:Tredclosed}.
\end{proof}

\subsection{Proof of Theorem~\ref{thm:main}}\label{sec:riisext}

We prove the following stronger version of Theorem~\ref{thm:main}. Its proof is an extension of the previous one. For readability 
statement (b) blurs the distinction between the symbol $f\in\PV(\alpha)$ and its interpretation in $\langle N,\alpha^N\rangle$.

\begin{theorem}\label{thm:1} 
Let $L$ be a finite language and $\varphi$ a strong basic $L$-sentence without built-in symbols.  Further, let $\tilde \varphi$ be  a weak finitary combinatorial principle in the language $\tilde L$. 

Then there exists a model $\langle N,\alpha^N\rangle$ of $\To$  such that
\begin{enumerate}\itemsep=0pt
\item[(a)] $
\str B(L,n,\alpha^N)\not\models\varphi
$ for some  $n\in N\setminus\{0\}$;
\item[(b)] $
\str B(\tilde L,m,f_{\bar a}^{-1}(0))\models\tilde\varphi
$ for all $m\in N\setminus\{0\}$, $f(x,\bar z)\in \PV(\alpha)$ and  tuples $\bar a$ from $N$.
\end{enumerate}

Moreover, if $\psi(x)$ is a $\PV$-formula that defines an unbounded set in~$\N$, then $N,n$ can be chosen such that $N$ is a $\PV$-cut in an elementary extension $M$ of $\N$ such that $N$ is subexponential in $n$ and   $M\models\psi(n)$. 
\end{theorem}

\begin{proof} Proceed as in the previous proof with two changes. First, since $\varphi$ is strong, we can additionally assume that the structure $\str B$ chosen in the beginning is $n^{o(1)}$-large. 
This ensures the assumptions of Corollary~\ref{cor:D2}. Second,
we let $\str D$ include additionally the countably many sets $D((t_{\tilde S})_{\tilde S\in\tilde L},m)$ from this corollary, where
$m$ runs over $M\setminus\{0\}$ and $(t_{\tilde S})_{\tilde S\in\tilde L}$ runs over families of decision trees  
 of height at most $b_0$ in $M$.
We are left to verify (b).

Let
 $f(x,\bar z)\in \PV(\alpha)$, $m\in N\setminus\{0\}$ and a tuple $\bar a$ from $N$ be given.
Choose $(t_{\tilde S})_{\tilde S\in\tilde L}$ according to Lemma~\ref{lem:family}. We show
$\str C((t_{\tilde S})_{\tilde S\in\tilde L},m,\alpha^N)\models\tilde\varphi$.												
We can assume that every~$t_{\tilde S}$ outputs~$0$ on arguments outside $[m]$ (otherwise modify $t_{\tilde S}$ adding $m$ to its parameters). Then every $t_{\tilde S}$ is a decision tree also in~$M$. 
As~$b_0$ bounds  lengths in $N$, the trees~$t_{\tilde S}$ have
height at most $b_0$. By genericity, there is $p\in G\cap D((t_{\tilde S})_{\tilde S\in\tilde L},m)$, so
$\str C((t_{\tilde S})_{\tilde S\in\tilde L},m,p)$ verifies~$\tilde\varphi$. By Lemma~\ref{lem:alphagood}, $\str C((t_{\tilde S})_{\tilde S\in\tilde L},m,\alpha^N)$ extends
$\str C((t_{\tilde S})_{\tilde S\in\tilde L},m,p)$ and hence verifies 
 $\tilde\varphi$ too.
\end{proof}


\section{Discussion}\label{sec:disc}

We discuss the applicability of Theorem~\ref{thm:main} using the examples from Section~\ref{sec:exas}. There we saw many strong finitary combinatorial principles and also that $\WPHP$ is weak. To these principles Theorem~\ref{thm:main} applies directly and thus, as stated in the Introduction, gives a simple and general criterion for independence from $Q_\WPHP$ over $\To$. The main limitation of the applicability of 
Theorem~\ref{thm:main} is that $\WPHP$ is our only natural example of a weak principle.
Despite its naturality, weakness seems to be a surprisingly restrictive condition. We are unable to offer any sort of explanation for this.

However, one can get independence from principles that are not weak
via Theorem~\ref{thm:main}:  

\begin{corollary}\label{cor:appl} $Q_\varphi$ is independent from $Q_{\tilde\varphi}$ over $\To$ for
\begin{eqnarray*}
\tilde\varphi&\in&\big\{\WPHP,\WPHP',\rPHP\big\},\\
\varphi&\in&\big\{\PHP,\LPHP,\OPHP,\PAR,\HOP,\IND\big\}.
\end{eqnarray*}
\end{corollary}

\begin{proof} For $\tilde\varphi=\WPHP$ this follows directly from Theorem~\ref{thm:main} because $\WPHP$ is weak and all listed choices for $\varphi$ are strong. The principles $\WPHP'$ and $\rPHP$ are not weak but both $Q_{\WPHP'}$ and $Q_\rPHP$ are consequences of $Q_\WPHP$ over $\PVa$, so our claim follows by Proposition~\ref{prop:trans}. For $\WPHP'$ this is well known (see \cite{jwphp} for this and other comparisons of various pigeonhole principles over $\PV(\alpha)$ and $\So$). For $\rPHP$ note that $Q_\rPHP$ is many-one reducible to $Q_\WPHP$ and apply Proposition~\ref{prop:Tred}~(b).
\end{proof}

Some of these independence results are known to hold in a much stronger form
following Ajtai's work: $Q_\OPHP$ is not provably total in $\BAo$ \cite{ajtai,kpwbip} while  $Q_\WPHP$ is provably total in $\mathsf T^2_2(\PV(\alpha))$~\cite{mpw}. Further,~$Q_\PAR$ is independent from  $Q_\PHP$ over~$\BAo$: this follows from Theorem~\ref{thm:bjstrong} and the exponential lower bound on  bounded depth Frege proofs~\cite{bp}. We refer to \cite{beameriis} and the references therein for more on counting principles. 

As mentioned in Example~\ref{ex:rphp}, the choice $\tilde\varphi=\rPHP$ implies that $Q_\varphi$ for $\varphi$ as in Corollary~\ref{cor:appl} is independent from $\So$ plus the surjective weak pigeonhole principle for $\PV(\alpha)$-functions. For $\varphi=\HOP$ this is known~\cite{atstha} even for  $\To$ instead $\So$.

\medskip

As in the proof of the corollary above we see that $Q_\psi$ is independent from $Q_{\tilde\varphi}$ over $\To$ if 
$Q_\varphi$ is and $Q_\varphi\le^m_pQ_\psi$; here, $\tilde\varphi,\varphi,\psi$ are arbitrary finitary combinatorial principles and
 $\le^m_p$ denotes (polynomial time) many-one reducibility.
 In this sense
all independence results in Corollary~\ref{cor:appl} follow from the ones for  
 $\OPHP$ and $\IND$:
 $$
  \begin{array}{ccccccc}
  &&&&&&\\
 &&\PHP&&&&\HOP\\[1ex] 
 &&\uparrow&\nwarrow&&\nearrow&\uparrow\\[1ex]
 \PAR&&\LPHP&&\IND&&\HDP\\[1ex]
 &\nwarrow&\uparrow&&&&\uparrow\\[1ex]
 &&\OPHP&&&&\HAP\\
 &&&&&&
 \end{array}
$$
In this figure, e.g.\ the arrow from $\HAP$ to $\HDP$ indicates $Q_\HAP\le^m_pQ_\HDP$. Recalling that $\PAR$ is not total, by  $Q_\OPHP\le^m_p O_\PAR$ we mean a many-one reduction $f,g,h$ as in \eqref{eq:manyone} of Section~\ref{sec:searchprbl} with the additional property that $g$ has only odd values. 

We give the reductions involving $\IND,\HDP$ and $\HAP$ below, the others are well-known.

\begin{remark}
The principles $\HDP$ and $\HAP$ are not well studied in proof 
complexity and Theorem~\ref{thm:main} does not seem to shed any light on their complexity. 
Their propositional proof complexity is  low:
the negations of their unary translations have polynomial size refutations in Res$(k)$ for some constant $k\in\N$. This follows from our proof of Proposition~\ref{prop:exareductions} below. There we give  quantifier free definitions of $\HAP$ and $\HDP$ in $\HOP$ in the sense of \cite[p.57f]{sergithesis}, and  this allows \cite[Lemma~15]{sergithesis} to translate well known short Resolution refutations of the negation of the unary translation of $\HOP$ \cite{stalmark} into short Res$(k)$ refutations as claimed.
\end{remark} 

\begin{proposition} \label{prop:exareductions}\
\begin{enumerate}\itemsep=0pt
\item[(a)] $Q_\HAP\le^m_pQ_\HDP$.
\item[(b)] $Q_\HDP\le^m_pQ_\HOP$.
\item[(c)] $Q_\IND\le^m_pQ_\HOP$.
\item[(d)] $Q_\IND\le^m_pQ_\PHP$.
\end{enumerate}
\end{proposition}

The proof will be easy based on the following ad hoc  lemma:

\begin{lemma}\label{lem:fored} Let $\tilde\varphi,\varphi$ be finitary combinatorial principles without built-in symbols in finite languages $\tilde L,L$ respectively. 
Assume there is a family $I:=(\delta_S)_{S\in L}$ 
of quantifier free $\tilde L$-formulas such that:
\begin{enumerate}\itemsep=0pt
\item[(i)] if $S\in L$ is a relation symbol, then 
$\delta_S$ has $\ar(S)$ many free variables;
\item[(ii)] if $S\in L$ is a function symbol, then 
$\delta_S$ has $\ar(S)+1$ many free variables and defines in every $\tilde L$-structure  the graph of some $\ar(S)$-ary function;
\item[(iii)] for every $\tilde L$-structure $\str B$ falsifying $\tilde\varphi$, the $L$-structure $I(\str B)$ falsifies $\varphi$; this structure has the same universe as $\str B$ and interprets $S\in L$ by the set defined by $\delta_S$ in $\str B$.
\end{enumerate}
Then $Q_{\tilde \varphi}\le^m_pQ_{\varphi}$. 
\end{lemma}

\begin{proof} Let $\tilde \varphi=\exists \bar y\tilde\psi(\bar y)$ and $\varphi=\exists \bar w\psi(\bar w)$ for quantifier free $\tilde \psi,\psi$ and recall $Q_{\tilde\varphi}$ and $Q_\varphi$ are $\q{\str B(\tilde L,x,\alpha)\models\tilde \psi(y)}$ and $\q{\str B(L,x,\alpha)\models\psi(w)}$ respectively. 

\medskip

\noindent\textit{Claim:} There exists $\hat{I}(u,v)\in\PV(\alpha)$ such that in every model $\langle N,\alpha^N\rangle$ of $\So$ and every $n\in N\setminus\{0\}$ we have
\begin{equation}\label{eq:equI}
\str B(\tilde L,n,\hat{ I}^{-1}_n(\alpha^N))= I\big(\str B(\tilde L,n,\alpha^N)\big).
\end{equation}

\noindent\textit{Proof of the Claim:}
We show that 
for every $S\in L$ there is $f_S(u,\bar x)\in\PV(\alpha)$ such that $f_S(n,\bar a)=S^{[n]}(\bar a)$ in $\langle\N,\alpha^\N\rangle$  for every $\alpha^\N\subseteq\N$
and $\bar a\in[n]^{\ar(S)}$ and $n\in\N\setminus\{0\}$; here, $S^{[n]}$ denotes the interpretation of $S$ in $I\big(\str B(\tilde L,n,\alpha^\N)\big)$.

This is clear for relation symbols. For a function symbol $S\in L$ observe that the empty theory 
proves $\exists y\delta_S(\bar x,y)$ by (ii). Hence,
Herbrand's theorem gives finitely many  $\tilde L$-terms $ t_0(\bar x),\ldots, t_{\ell-1}(\bar x)$ such that $\bigvee_{i<\ell}\delta_S(\bar x,t_i(\bar x))$ is valid.
Then $S^{[n]}(\bar a)$ can be computed in polynomial time with oracle $\alpha^\N$ by testing which of $ t_0(\bar a),\ldots, t_{\ell-1}(\bar a)$ satisfies $\delta_S(\bar a,y)$ in~$\str B(\tilde L,n,\alpha^\N)$. 

The function $\hat I$ is easily constructed from the functions $f_S,S\in L$, so that \eqref{eq:equI} holds in $\langle \N,\alpha^\N\rangle$ for all $\alpha^\N\subseteq\N$ and all $n\in\N\setminus\{0\}$. 
To see \eqref{eq:equI} holds in $\langle N,\alpha^N\rangle$ let $S\in L$ be a unary function symbol; other symbols are treated similarly. We have to show that 
$$
\str B(\tilde L,n,\hat{ I}^{-1}_n(\alpha^N))\models S(a){=}b\; \Longleftrightarrow \;
I(\str B(\tilde L,n,\alpha^N)\models S(a){=}b. 
$$
The l.h.s.\ is equivalent to $f_S(a)=b$ (in $\langle N,\alpha^N\rangle$) because this equivalence is expressed by a $\Delta_0^b(\PV(\alpha))$ sentence, so proved by $\PVa$ by Lemma~\ref{lem:QE}. The r.h.s. too is equivalent to $f_S(a)=b$. Indeed, let $\delta'_S(u,v,\bar u)$ be a simple $\tilde L$-formula such that $\exists \bar u\delta_S(u,v,\bar u)$ is  logically equivalent to $\delta_S(u,v)$; intuitively, the variables $\bar u$ collect values of (sub)terms appearing in~$\delta_S$. Then  $\So$ proves (recall Lemma~\ref{lem:formAB})
$$
 u{<}x\wedge v{<}x \to \Big(\exists z \Big( 
\q{\str B(\tilde L,x,\alpha)\models \delta'_S(z) }
\wedge (z)_0{=}u\wedge(z)_1{=}v\Big)\leftrightarrow f_S(u){=}v\Big).
$$
This implies the claim.\hfill$\dashv$\medskip

Recalling Lemma~\ref{lem:formAB},  the Claim and (iii) imply that $\So$ proves
$$
\q{\str B(L,x,\hat{ I}^{-1}_x(\alpha))\models\psi(w)}\to \exists y\q{\str B(\tilde L,x,\alpha)\models\tilde \psi(y)}
$$
By witnessing, there is $h(x,w)\in \PV(\alpha)$ witnessing $y$. This implies $Q_{\tilde \varphi}\le^m_pQ_{\varphi}$ (using the identity function and $\hat I$ for $g$ and $f$ in \eqref{eq:manyone} of Section~\ref{sec:searchprbl}).
\end{proof}

\begin{proof}[Proof of Proposition~\ref{prop:exareductions}] 
For (a), given a Boolean algebra $\str B$ falsifying $\HAP$ we falsify $\HDP$ taking for $\prec$ the proper subset relation (in the sense of $\str B)$; a point between $a$ and a proper superset $b$ is obtained adding to $a$ a 
proper non-empty subset of $b\setminus a$ (in the sense of $\str B$); such a subset is found by $f^B$. More precisely, we apply the previous lemma with~$I$ collecting the following formulas:
\begin{eqnarray*}
\delta_\prec(x_0,x_1)&:=& x_0{\sqcap} x_1{=}x_0\wedge \neg x_0{=}x_1,\\
\delta_b(x_0,x_1,y)&:=& (y{=} x_0{\sqcup}f(x_1{\sqcap}{\sim}x_0)\wedge \delta_\prec(x_0,x_1))\vee(y{=} 0\wedge \neg\delta_\prec(x_0,x_1)),\\
\delta_0(y)&:=& y{=}0,\\
\delta_1(y)&:=& y{=}1.
\end{eqnarray*}

For (b) we use a variant of \cite[Example~2, p.65]{sergithesis}: given $\str B$ violating $\HDP$ we find a $\{\prec,f\}$-structure falsifying $\HOP$ by taking the $\prec^B$-interval $[0,1]$, with regressive function $b^B(0,\cdot)$ and declaring everything outside $[0,1]$ to be pairwise incomparable and bigger than~1. 
More precisely, writing $\q{x\in [0,1]}$ for
$
 x{=}0 \vee x{=}1 \vee (0{\prec}x\wedge x{\prec} 1)
$
and $\q{x\notin[0,1]}$ for its negation, 
\begin{eqnarray*}
\delta_\prec(x_0,x_1)&:=& (\q{x_0\in [0,1]}\wedge \q{x_1\in[0,1]}\wedge x_0{\prec}x_1)\ \vee \ 
(\q{x_1\notin[0,1]}\wedge \q{x_0\in [0,1]}),
\\
\delta_f(x,y)&:=& (\q{x\in[0,1]}\wedge y{=}b(0,x))\vee (\q{x\notin[0,1]}\wedge y{=}1).
\end{eqnarray*}

For (c), given $\str B$ falsifying $\IND$ we get a structure falsifying $\HOP$ by taking the inverse of the order of $\str B$ restricted to $P^B$, declaring everything outside $P^B$ to be pairwise incompatible and bigger than $\textit{min}^B$, and taking $s^B$ as regressive function. More precisely,
\begin{eqnarray*}
\delta_\prec(x_0,x_1)&:=& (P(x_0)\wedge P(x_1)\wedge x_1{\prec} x_0)\ \vee\
(\neg P(x_1)\wedge P(x_0)),
\\
\delta_f(x,y)&:=& y{=}s(x).
\end{eqnarray*}

For (d), take $ y{=}\textit{min}$ for $\delta_c(y)$, and 
$
(P(x)\wedge y{=}s(x))\vee(\neg P(x)\wedge y{=}x)
$
for $\delta_f(x,y)$.
\end{proof}

\subsubsection*{Acknowledgements}

I thank the referee for detailed comments.
I thank Neil Thapen and Emil Je\v r\'abek for their help understanding the material in Section~\ref{sec:aux} during a visit to Prague supported by the ERC advanced grant 339691 (FEALORA).

\end{document}